\documentclass[12pt]{amsart}
\usepackage{amsmath,amssymb}

\theoremstyle{plain}
\newtheorem{thm}{Theorem}[section]
\newtheorem{theorem}[thm]{Theorem}

\newtheorem{lemma}[thm]{Lemma}
\newtheorem{corollary}[thm]{Corollary}
\newtheorem{proposition}[thm]{Proposition}
\theoremstyle{definition}
\newtheorem{remark}[thm]{Remark}

\newtheorem{definition}[thm]{Definition}

\newtheorem{example}[thm]{Example}

\numberwithin{equation}{section}
\newcommand{\0}{{\mathcal O}}

\newcommand{\sC}{{\mathcal C}}
\newcommand{\sD}{{\mathcal D}}

\newcommand{\sF}{{\mathcal F}}
\newcommand{\sG}{{\mathcal G}}
\newcommand{\sH}{{\mathcal H}}

\newcommand{\sK}{{\mathcal K}}

\newcommand{\sM}{{\mathcal M}}
\newcommand{\sN}{{\mathcal N}}
\newcommand{\sO}{{\mathcal O}}

\newcommand{\sQ}{{\mathcal Q}}

\newcommand{\sS}{{\mathcal S}}

\newcommand{\sU}{{\mathcal U}}

\newcommand{\sW}{{\mathcal W}}

\newcommand{\BO}{{\mathbb O}}

\newcommand{\C}{{\mathbb C}}

\newcommand{\BP}{{\mathbb P}}
\newcommand{\pit}{{\mathbb P}}
\newcommand{\Q}{{\mathbb Q}}

\newcommand{\Z}{{\mathbb Z}}
\newcommand{\End}{{\rm End}}

\def\Sec{\mathop{\rm Sec}\nolimits}

\newcommand{\fg}{{\mathfrak g}}
\newcommand{\fsl}{{\mathfrak s}{\mathfrak l}}
\newcommand{\fgl}{{\mathfrak g}{\mathfrak l}}

\newcommand{\fp}{{\mathfrak p}}

\newcommand{\fso}{{\mathfrak s}{\mathfrak o}}
\newcommand{\fco}{{\mathfrak c}{\mathfrak o}}
\newcommand{\fsp}{{\mathfrak s}{\mathfrak p}}
\newcommand{\aut}{{\mathfrak a}{\mathfrak u}{\mathfrak t}}
\newcommand{\fcsp}{{\mathfrak c}{\mathfrak s}{\mathfrak p}}
\newcommand\Aut{\rm Aut}
\newcommand\sd{\!>\!\!\!\! \lhd \:}

\def\RatCurves{\mathop{\rm RatCurves}\nolimits}
\def\Gr{\mathop{\rm Gr}\nolimits}

\def\Lag{\mathop{\rm Lag}\nolimits}
\def\Sym{\mathop{\rm Sym}\nolimits}
\def\Pic{\mathop{\rm Pic}\nolimits}
\def\Hom{\mathop{\rm Hom}\nolimits}
\def\Im{\mathop{\rm Im}\nolimits}
\def\Ker{\mathop{\rm Ker}\nolimits}

\def\Bl{\mathop{\rm Bl}\nolimits}

\def\SYM{\mathop{\rm Sym}\nolimits}

\title[Projective manifolds with non-zero prolongations]{Classification of non-degenerate projective varieties with
non-zero prolongation and application to target rigidity}
\author{Baohua Fu and Jun-Muk Hwang}
\thanks{Baohua Fu is supported by National Scientific Foundation of
China and the KIAS Scholar Program.  Jun-Muk Hwang is supported by
National Researcher Program 2010-0020413 of NRF and MEST}
\begin{document}
\maketitle  \setcounter{tocdepth}{1}
\begin{abstract}
The prolongation $\fg^{(k)}$ of a linear Lie algebra $\fg \subset
\fgl(V)$ plays an important role in the study of symmetries of
$G$-structures. Cartan and Kobayashi-Nagano have given a complete
classification of irreducible linear Lie algebras $\fg \subset
\fgl(V)$ with non-zero prolongations.

If $\fg$ is the Lie algebra $\aut(\hat{S})$ of infinitesimal
linear automorphisms of a projective variety $S \subset \BP V$,
its prolongation $\fg^{(k)}$ is related to the symmetries of  cone
structures, an important example of which is the variety of
minimal rational tangents in the study of uniruled projective
manifolds. From this perspective, understanding  the prolongation
$\aut(\hat{S})^{(k)}$ is useful in questions related to the
automorphism groups of uniruled projective manifolds. Our main
result is  a complete classification of irreducible non-degenerate
nonsingular variety $S \subset \BP V$ with $\aut(\hat{S})^{(k)}
\neq 0$, which can be viewed as a generalization of the result of
Cartan and Kobayashi-Nagano.   As an application, we show that
when $S$ is linearly normal and $\Sec(S) \neq \BP V$, the blow-up
${\rm Bl}_S(\BP V)$ has the target rigidity property, i.e., any
deformation of a surjective morphism $f:Y \to {\rm Bl}_S(\BP V)$
comes from the automorphisms of ${\rm Bl}_S(\BP V)$.
\end{abstract}

 \tableofcontents
\section{Introduction}
For a linear algebraic group $G \subset {\rm GL}(V)$, a
$G$-structure on a complex manifold $M$ with $\dim M = \dim V$ is
a $G$-subbundle of the frame bundle on $M$. Many classical
geometric structures in differential geometry are $G$-structures
for various choices of $G$. For this reason, the
(self)-equivalence problem for $G$-structures has been studied
extensively. It turns out that the graded pieces (under a natural
filtration) of the Lie algebra of infinitesimal symmetries of
$G$-structure are contained in the prolongations  $\fg^{(i)}, i
\geq 1$, of the Lie algebra $\fg \subset \fgl(V)$ (cf. Definition
\ref{d.prolongation} for a precise definition) and in fact, equal
to the prolongations when the $G$-structure is flat (cf.
Proposition \ref{p.HM05}). In other words, an essential
information of the symmetries of $G$-structures is encoded in
$\fg^{(i)}$. A fundamental result in the study of the
prolongations is the following result of E. Cartan, S. Kobayashi
and T. Nagano.

\begin{theorem}\label{t.CKN}
Let $\fg \subset \fgl(V)$ be an irreducible representation of a
Lie algebra $\fg$.
\begin{enumerate} \item[(1)] If $\fg^{(2)} \neq 0$, then
$\fg= \fgl(V), \fsl(V), \fsp(V)$ or $\fcsp(V)$ where $\dim V$ is
even for the last two cases. \item[(2)] If $\fg^{(2)} =0$, but
$\fg^{(1)} \neq 0$, then $\fg \subset \fgl(V)$ is isomorphic to
the isotropy representation on the tangent space at a base point
of an irreducible Hermitian symmetric space of compact type,
different from the projective space.
\end{enumerate} \end{theorem}

 The proof of Theorem
\ref{t.CKN} as given in \cite{KN} is purely algebraic, and depends
heavily on the theory of semi-simple Lie algebras and their
representations. For that reason, there is little hope of
generalizing it to non-reductive Lie algebras.

In \cite{HM05}, motivated by algebro-geometric questions, the
prolongation of $\fg \subset \fgl(V)$ associated to a projective
variety $S \subset \BP V$ was studied. More precisely,
  for a projective subvariety $S
\subsetneq \BP(V)$, consider the Lie algebra $\aut(\hat{S})
\subset \fgl(V)$ of infinitesimal linear automorphisms of the
affine cone $\hat{S}$. \cite{HM05} shows that one can study
prolongations $\aut(\hat{S})^{(k)}$ using projective geometry of
$S \subset \BP V$ and the deformation theory of rational curves on
$S$. Combining these two geometric tools, the following
generalization of Theorem \ref{t.CKN} (1) is proved in Theorem
1.1.2 of \cite{HM05}.

\begin{theorem}\label{t.HM112}
Let $S \subset \BP V$ be an irreducible nonsingular non-degenerate
projective variety. If $\aut(\hat{S})^{(2)} \neq 0$, then $S = \BP
V$. \end{theorem}

It is easy to derive Theorem \ref{t.CKN} (1) from Theorem
\ref{t.HM112}. On the other hand, the latter is stronger than the
former, because there is no a priori reason that $\aut(\hat{S})$
is reductive in Theorem \ref{t.HM112}. For example, for the
deformation rigidity studied in \cite{HM05}, it is essential to
have this stronger result.

It is natural to ask the generalization of Theorem \ref{t.CKN} (2)
in the form of Theorem \ref{t.HM112}. Some partial results in this
direction was obtained in \cite{HM05} (e.g. Theorem \ref{t.HM}
below).  The goal of this paper is to give a complete answer to
this question in the following form.

\vspace{0.5cm}
 {\bf Main Theorem}. Let $S \subsetneq \BP V$ be an irreducible nonsingular
  non-degenerate
 variety such that $\aut(\hat{S})^{(1)} \neq 0$. Then $S \subset
\BP V$ is projectively  equivalent to one of the following:

\begin{enumerate} \item[(A1)] the second Veronese embedding
$v_2(\BP^n) \subset \BP^{\frac{1}{2}(n^2+3n)} $ for $n \geq 2$;
\item[(A2)] Segre embedding $\BP^a \times \BP^b \subset
\BP^{ab+a+b} $ for $a,b \geq 2$; \item[(A3)] a natural embedding
$\BP(\sO_{\BP^{k}}(-1)^m \oplus \sO_{\BP^{k}}(-2)) \subset
\BP^{m(k+1) + \frac{1}{2}(k+2)(k+1) -1}$ for $k\geq 2,m \geq 1$;
\item[(B1)] odd-dimensional hyperquadrics $\Q^1, \; \Q^3,
\;\ldots, \; \Q^{2 \ell -1}, \ldots$; \item[(B2)] even-dimensional
hyperquadrics $\Q^2, \; \Q^4, \; \ldots, \; \Q^{2\ell}, \ldots$;
\item[(B3)]
 Segre embedding  $\BP^1
\times \BP^m \subset \BP^{2m+1}$ and Pl\"ucker embedding  $\Gr(2,
\C^{m+3}) \subset \BP^{\frac{1}{2}(m^2+5m+4)}$ for $m \geq 3$;
\item[(B4)] Segre embedding  $\BP^1 \times \BP^2 \subset \BP^5$,
Pl\"ucker embedding  $\Gr(2,\C^5) \subset \BP^9$, spinor embedding
$\mathbb{S}_5 \subset \BP^{15}$ and the $E_6$-Severi embedding
$\BO \BP^2
 \subset \BP^{26}$; \item[(B5)] general hyperplane sections of the first three in (B4), i.e., $(\BP^1 \times \BP^2) \cap H_0 \subset \BP^{4}, \; \Gr(2,\C^5) \cap
H_1 \subset \BP^8, \; \mathbb{S}_5 \cap H_2 \subset \BP^{14}$;
\item[(B6)] general hyperplane section of the first in (B3), i.e.
$(\BP^1 \times \BP^m) \cap H \subset \BP^{2m}$, for $m \geq 3$.
 \item[(C)]
some biregular projections of (A1), (A2), (A3) and $\Gr(2,
\C^{m+3})$ in (B3).
\end{enumerate}\vspace{0.2cm}

The varieties in (A1)-(A3) and (C) satisfy $\Sec(S) \neq \BP V$
while the first entries of (B1)-(B6) verify $\Sec(S)=\BP V$. Note
that the varieties in (B1)-(B5) are listed as  sequences $S_0,
S_1, S_2, \ldots$. The reason behind this  way of listing the
varieties will become clear in the course of the proof of Main
Theorem. In fact, the variety $S_i$ is the VMRT (cf. Definition
\ref{d.VMRT}) of the variety $S_{i+1}$, a crucial fact in the
proof of Main Theorem. More detailed description of the varieties
(A1)-(B6) and the explicit computation of the prolongation
$\aut(\hat{S})^{(1)}$ for each of them are given in Section
\ref{s.examples}. Which biregular projections in (C) have non-zero
prolongations will be described completely in Section
\ref{s.prolong}. One may have the impression that compared with
the linearly normal cases of (A1)-(B6), the projections in (C) are
mere technicalities. This is not the case. In fact, in the
induction process of the proof of Main Theorem, it is crucial to
understand the cases in (C). In other words, imposing the
additional condition of linear normality on $S \subset \BP V$ in
Main Theorem would not make the proof any simpler, and it is
essential to include varieties which are not necessarily linearly
normal to carry out the proof of Main Theorem.

As we will  explain in Section \ref{s.G-structure},
$\aut(\hat{S})^{(1)}$ is an essential part of the symmetries of
cone structures, in particular, the structure coming from the
varieties of minimal rational tangents, which is an important tool
in the study of uniruled projective varieties. In this respect,
 Main Theorem will be useful in algebraic geometric questions involving  automorphism groups
 of uniruled varieties.
As an example, we will give a direct application of Main Theorem
in Section \ref{s.target}, in the proof of the target rigidity for
the blow-up of  $\BP V$ along $S$.  More precisely, we shall show
(cf. Corollary \ref{c.comp}) that if $S \subset \BP V$ is an
irreducible nonsingular non-degenerate linearly normal variety
such that $\Sec(S) \neq \BP V$, then any deformation  $f_t: Y \to
\Bl_S(\BP V)$ of a surjective morphism  $f_0: Y \to \Bl_S(\BP V)$
comes from automorphisms of $\Bl_S(\BP V)$.

 Turing to the proof of Main Theorem, the main strategy is
  to carry out an induction on
 VMRT. In fact, by the partial result in \cite{HM05} and
 the work of Ionescu-Russo
 \cite{IR}, the question is quickly reduced to the case when $S
 \subset \BP V$ has Picard number 1 and covered by lines. In this
 setting, we show in Proposition \ref{p.VMRTCC} and Theorem \ref{t.induction} that the VMRT of $S$ at a general point, say,
 $S' \subset \BP V'$, is again an irreducible nonsingular non-degenerate projective variety
 with $\aut(\hat{S}')^{(1)} \neq 0$. By induction, we have a
 classification of $S' \subset \BP V'$. From the information on
 $S' \subset \BP V'$, we can recover $S \subset \BP V$ by
 Cartan-Fubini type extension theorem as explained in Corollary
 \ref{c.reconstruct}. An essential ingredient in this induction process
 is the local flatness of the associated cone structure, or equivalently,
  G-structure. For that purpose, we develop some general theory of
  these differential geometric machinery in Section
  \ref{s.G-structure}.

 The induction process enables us to prove Main Theorem, modulo the
 termination of the sequence of varieties in (B3)-(B6).
 Among these, the termination of (B3) and (B4) is an easy consequence of
 the condition on the secant varieties, via a result from \cite{HK}.
 The termination of (B5) and (B6) is more complicated and
 technically demanding. It will be proved in Section
 \ref{s.non1} and Section \ref{s.non2}. Many of the
 geometric ideas in these two sections are borrowed from
 Section 6 and Section 8 of \cite{HM05}. However, the main line of
 arguments and details of the proof are rather different from \cite{HM05},
  except
 Propositions
 \ref{p.submanifold} and \ref{p.submanifoldQ} whose proofs are essentially contained in those of Proposition 6.3.4 and
 Proposition 8.3.4 of \cite{HM05}, respectively.

\section{Prolongation of a projective variety: basic properties}

\begin{definition} \label{d.prolongation}
Let $V$ be  a complex vector space and $\fg \subset {\rm End}(V)$
a Lie subalgebra. The {\em $k$-th prolongation} (denoted by
$\fg^{(k)}$) of $\fg$ is the space of symmetric multi-linear
homomorphisms $A: \Sym^{k+1}V \to V$ such that for any fixed $v_1,
\cdots, v_k \in V$, the endomorphism $A_{v_1, \ldots, v_k}: V \to
V$ defined by $$v\in V \mapsto A_{v_1, \ldots, v_k, v} := A(v,
v_1, \cdots, v_k) \in V$$ is in $\fg$. In other words, $\fg^{(k)}
= \Hom(\Sym^{k+1}V, V) \cap \Hom(\Sym^kV, \fg)$.
\end{definition}

It is immediate from the definition that $\fg^{(0)} = \fg$ and if
$\fg^{(k)}=0$, then $\fg^{(k+1)} = 0$.

In this paper, we are interested in the case where $\fg$ arises
from geometric situations. First recall some basic definitions.

\begin{definition}\label{d.basic}
Let $S \subset \BP V$ be an irreducible subvariety.
\begin{enumerate}
 \item[(i)]
$S$  is said to be {\em non-degenerate} (resp. {\em linearly
normal}) if the restriction map $H^0(\BP V, \sO_{\BP V}(1)) \to
H^0(S, \sO_{S}(1))$ is injective (resp. surjective). \item[(ii)]
$S$ is said to be
 {\em conic-connected} if through two general
points of $S$, there passes an irreducible conic contained in $S$.
\item[(iii)] The {\em secant variety} $\Sec(S) \subset \BP V$ of
$S$ is the closure of the union of lines through two points of
$S$. \item[(iv)] The projective automorphism group of
 $S \subset \BP V$ is  $$\Aut(S):=\{g \in {\rm PGL}(V) |
g S = S\}.$$  Its Lie algebra will be denoted by $\aut(S)$.
\item[(v)] Denote by $\hat{S} \subset V$ the affine cone of $S$
and by $T_{\alpha}(\hat{S}) \subset V$ the tangent space at a
smooth point $\alpha \in \hat{S}$. The Lie algebra of
infinitesimal linear automorphisms of $\hat{S}$ is
$$\aut(\hat{S}):=\{g \in \End(V)| g(\alpha) \in T_{\alpha}(\hat{S})
\mbox{ for any smooth point } \alpha \in \hat{S} \}.$$ Its
prolongation  $\aut(\hat{S})^{(k)}$ will be called the $k$-th {\em
prolongation} of $S \subset \BP V$.
\end{enumerate}
\end{definition}

We have the following vanishing result.
\begin{theorem}[\cite{HM05}, Theorem 1.1.2]\label{t.HMk=2}
Let $S \subsetneq \BP V$ be an irreducible nonsingular
non-degenerate subvariety. Then $\aut(\hat{S})^{(k)} = 0$ for all
$k \geq 2$.
\end{theorem}

However, there are several examples of $S$ with non-zero first
prolongation $\aut(\hat{S})^{(1)}$. In \cite{HM05}, some partial
results on the structure of such varieties were obtained. Here we
collect them with some immediate improvements.

The following theorem is  essentially proved  in  Theorem 1.1.3 of
\cite{HM05} for linearly normal $S \subset \BP V$. We will explain
how the proof in \cite{HM05} can be modified to give the general
result.

\begin{theorem}\label{t.HM}
Let $S \subsetneq \BP V$ be a nonsingular non-degenerate
projective subvariety. Let $\sO(1)$ be the hyperplane line bundle
and $S \to \BP H^0(S, \sO(1))^*$ be the linearly normal embedding
inducing the inclusion $\iota: \hat{S} \to H^0(S, \sO(1))^*$. We
have a natural projection $p: H^0(S, \sO(1))^* \to V$ satisfying
$p (\iota (\hat{S})) = \hat{S}$.
 Then the following holds.
\begin{enumerate} \item[(i)] If $\aut(\hat{S})^{(1)} \neq 0,$ then $S$ is
  conic-connected. \item[(ii)]  For each non-zero $A \in   \aut(\hat{S})^{(1)}, $ there exists
a non-zero $\lambda_A \in H^0(S, \sO(1))$ such that for each
$\alpha \in \hat{S}$, $A_{\alpha \alpha} =
\lambda_A(\iota(\alpha)) \alpha$.
%For any $\alpha, \beta \in \hat{S}$ general, we have $A_{\alpha
%\beta} \in T_\alpha(\hat{S}) \cap T_\beta(\hat{S})$.
  \item[(iii)] In the notation of (ii), for any $\alpha \in \hat{S}$ and $\alpha' \in T_{\alpha}(\hat{S})$,
  $$ \lambda_A(\iota(\alpha)) \alpha' + \lambda_A(\iota({\alpha'})) \alpha = 2
A_{\alpha \alpha'}.$$ In particular, the endomorphism $A_{\alpha}$
acts on the tangent space of $S$
$$T_{[\alpha]}(S) = {\rm Hom}(\C \alpha, T_{\alpha}(\hat{S})/\C
\alpha)$$ as the scalar multiplication by $\frac{1}{2}
\lambda_A(\iota(\alpha))$.
  \item[(iv)] Suppose $\aut(\hat{S})^{(1)} \neq 0$.
  Then for a general point $s \in S$, there exists an element $E \in \aut(\hat{S})$
   which generates a $\C^\times$-action on $S$ with an isolated fixed point at $s$
    such that the isotropy action on $T_s(S)$ is the scalar multiplication by $\C^\times$. \end{enumerate} \end{theorem}

\begin{proof}
 By a similar proof as that of Theorem 1.1.3 (ii)
\cite{HM05}, there exists a point $s_o \in S$ such that
  $S$ is covered by conics passing through $s_o$. By Lemma 1 in \cite{Fu}, this implies that
$S$ is conic-connected, proving (i).  The proof of (ii) is the
same as that of Proposition 2.3.1 in \cite{HM05}. The proof of
(iii) and (iv) is a modification of that of Theorem 1.1.3 (iii) in
\cite{HM05}.  In fact, fix $\alpha \in \hat{S}$ outside the zero
locus of $\lambda_A$ and pick any $\alpha' \in
T_{\alpha}(\hat{S})$. By the natural identification of $\hat{S}$
and $\iota(\hat{S}),$ we can regard $\alpha$ and $\alpha'$ as
vectors in $H^0(S, \sO(1))^*$. Then choosing a holomorphic arc in
$\hat{S}$ passing through $\alpha$ with tangent $\alpha'$,
$$\lambda_A(\alpha + t \alpha' + \cdots)\; (\alpha + t \alpha' +
\cdots) = A_{\alpha + t \alpha' + \cdots, \alpha + t \alpha' +
\cdots} =A_{\alpha \alpha} + 2t A_{\alpha \alpha'} + \cdots$$
where $(\cdots)$ stands for terms containing $t^2$-factor. From
this, the equation in (iii)  follows.
 Thus as in pp. 606--607 of  \cite{HM05}, we can say that
the endomorphism $A_{\alpha}$ acts on the tangent space
$$T_{[\alpha]}(S) = {\rm Hom}(\C \alpha, T_{\alpha}(\hat{S})/\C
\alpha)$$ as the scalar multiplication by $\frac{1}{2}
\lambda_A(\iota(\alpha))$. If we choose $\alpha$ outside the zero locus
of $\lambda_A,$  then the semi-simple part of $A_{\alpha}$ in
$\aut(\hat{S})$ generates the required $\C^\times$-action.
\end{proof}

We have the following variation of Lemma 2.3.3 in \cite{HM05}. The
proof there works verbatim.
\begin{lemma}\label{l.sumihiro}
Let $ S \subsetneq \BP V$ be a nonsingular non-degenerate
projective subvariety which is not biregular to a projective
space. Let $A \in \aut(\hat{S})^{(1)}$. Suppose for some $\alpha
\in V$ and a subspace $H \subset  V$ of codimension 1, the
endomorphism $A_\alpha$ satisfies $A_{\alpha \beta} = 0$ for all
$\beta \in H \cap \hat{S}$. Then $A_{\alpha} = 0$.
\end{lemma}

 Theorem \ref{t.HM}  has the
following consequences.

\begin{proposition}\label{p.secondorder}
In the setting of Theorem \ref{t.HM} (ii), assume $S$ is not
biregular to a projective space.
 Choose a general point
of the hyperplane section  $$\alpha \in \hat{S} \cap (\lambda_A =
0).$$ Then the vector field on $S$ induced by $A_{\alpha}$  is not
identically zero and vanishes at $\bar{\alpha} \in S$ to second
order.
\end{proposition}
\begin{proof}
Suppose that for any general point $\alpha \in \hat{S} \cap
(\lambda_A =0)$, the vector field on $S$ induced by $A_{\alpha}$
is identically zero on $S$, i.e., for each $\beta \in \hat{S}$,
$A_{\alpha \beta}$ is proportional to $\beta$. Then it is
proportional to $\alpha$ by symmetry. We conclude that
$A_{\alpha}$ is identically zero on $\hat{S}$, thus on $V$. By
symmetry, for each $\gamma \in V$, $A_{\gamma}$ vanishes on
$\hat{S} \cap (\lambda_A =0)$. Then by Lemma \ref{l.sumihiro},
$A_{\gamma}$ is identically zero, a contradiction to $A \neq 0$.
This shows that the vector field on $S$ induced by $A_{\alpha}$ is
not identically zero.

Now Theorem \ref{t.HM} (iii) says that this vector field on $S$
vanishes to second order at $\alpha\in \hat{S}$ with
$\lambda_A(\alpha) =0.$
\end{proof}

\begin{proposition} \label{p.aut=1}
Let $S \subset \BP V$ be a nonsingular non-degenerate subvariety.
Then
$$\dim \aut(\hat{S})^{(1)} \neq 1.$$
\end{proposition}

\begin{proof}
Let us write $\fg = \aut(\hat{S}) \subset \fgl(V)$. Assuming that
$\dim \fg^{(1)} =1,$ we will derive a contradiction.  Let $G
\subset {\rm GL}(V)$ be the connected component of the linear
automorphism group of the cone $\hat{S} \subset V$ whose Lie
algebra is $\fg$.
  Note that $G$ contains the
central subgroup $\C^\times \cdot {\rm Id}$.  The natural
$G$-action on $\Hom(\Sym^2V, V)$ induces a $G$-action
$$ \chi: G \to {\rm GL}(\fg^{(1)}) \cong \C^\times,$$ which is a
character of $G$. Let $G' \subset G$ be the kernel of $\chi$ and
$\fg' \subset \fg$ be its Lie algebra.
 Since $\fg^{(1)} \subset {\rm
Hom}(\SYM^2V, V)$, the central subgroup $\C^\times \cdot {\rm Id}$
acts non-trivially on $\fg^{(1)}$ and the normal subgroup $G'
\subset G$ is complementary to $\C^\times \cdot {\rm Id}.$ Thus we
have a direct sum decomposition of the Lie algebra $\fg= \C \cdot
{\rm Id} \oplus \fg'.$ Let $\bar{G} \subset {\rm PGL}(V)$ be the
image of $G$, under the projection ${\rm GL}(V) \to {\rm PGL}(V)$.
Then $\bar{G}$ is the identity component of the projective
automorphism group of $S$. The homomorphism $G'\to \bar{G}$ has
finite kernel and the Lie algebra $\fg'$ is isomorphic to the Lie
algebra of $\bar{G}$.

From $\fg^{(1)} \neq 0$ and Theorem \ref{t.HM} (iv), for each
general point $x \in S$, we have
 a
$\C^\times$-subgroup $G_x \subset G'$ which acts as the
multiplication by $\C^\times$ on the tangent space $T_x(S)$. Let
$T_x(\hat{S})$ be the affine tangent space at $x$. Since $G_x$ has
weight 1 on $T_x(S)$, it has exactly two distinct weights on
$T_x(\hat{S})$. In fact, from $T_x(S) = {\rm Hom}(\hat{x},
T_x(\hat{S})/\hat{x})$, if it has weight $k$ on $\hat{x}$, the
other weight on $T_x(\hat{S})/\hat{x}$ must be $k+1$.

Pick a vector $\alpha \in \hat{x}$ and $\alpha' \in
T_x(\hat{S})/\hat{x}.$ For $A \in \fg^{(1)}$ and $g \in G'$,
$g\cdot A = A$ implies that $$A_{g\alpha, g\alpha'} = g \cdot
A_{\alpha, \alpha'}.$$ On the other hand,  by Theorem \ref{t.HM}
(iii) we have
$$2 A_{\alpha, \alpha'} = \lambda(\iota(\alpha)) \alpha' +
\lambda(\iota(\alpha')) \alpha,$$ for some non-zero $\lambda \in
H^0(S, \0(1))$ and an injection $\iota: \hat{S} \to H^0(S,
\sO(1))^*.$
 Thus for any $t \in G_x \cong \C^\times,$
$$2 t\cdot A_{\alpha, \alpha'}= t^{k+1} \lambda(\iota(\alpha)) \alpha' +
t^k \lambda(\iota(\alpha')) \alpha,$$ while
$$ 2 A_{t \cdot \alpha, t\cdot \alpha'} = 2 A_{t^k \alpha, t^{k+1}
\alpha'} = 2 t^{2k +1} A_{\alpha, \alpha'} =  t^{2k+1}
\lambda(\iota(\alpha)) \alpha' +  t^{2k+1} \lambda(\iota(\alpha'))
\alpha.$$ Thus either $\lambda(\iota(\alpha)) = 0$ or
$\lambda(\iota(\alpha')) =0$. Since the set of  such
$\iota(\alpha)$ or $\iota(\alpha')$ spans the vector space $H^0(S,
\sO(1))^*$, we get $\lambda=0$, a contradiction.
\end{proof}

\section{Examples of linearly normal varieties with non-zero first prolongation }\label{s.examples}

 In this section, we will list
examples of linearly normal $S \subset \BP V$ with non-zero first
prolongation. Before we give
 these examples, it is convenient to recall the notion of VMRT, because our examples arise as VMRT of
 some uniruled manifolds.

\begin{definition}\label{d.VMRT}
Let $X$ be a uniruled projective manifold. An irreducible
component $\sK$ of the space $\RatCurves^n(X)$ of rational curves
on $X$ is called {\em a minimal rational component} if the
subvariety $\sK_x$ of $\sK$ parameterizing curves passing through
a general point $x \in X$ is non-empty and proper. Curves
parameterized by $\sK$ will be called {\em minimal rational
curves}. Let $\rho: \sU \to \sK$ be the universal family and $\mu:
\sU \to X$ the evaluation map. The tangent map $\tau: \sU
\dasharrow \BP T(X)$ is defined by $\tau(u) = [T_{\mu(u)}
(\mu(\rho^{-1} \rho(u)))] \in \BP T_{\mu(u)}(X)$. The closure $\sC
\subset \BP T(X)$ of its image is the total space of {\em  variety
of minimal rational tangents}. The natural projection $\sC \to X$
is a proper surjective morphism and a general fiber $\sC_x \subset
\BP T_x(X)$ is called the {\em variety of minimal rational
tangents} (VMRT for short) at the point $x \in X$.
\end{definition}
The following is well-known (cf. Proposition 1.5 in \cite{Hw01}).

\begin{proposition}\label{p.nonsingular}
Let $X \subset \BP^N$ be a nonsingular projective variety  covered
by lines. A component of family of lines covering $X$ is a minimal
rational component and the VMRT $\sC_x \subset \BP T_x(X)$ at a
general point $x \in X$ is nonsingular. \end{proposition}

The following is immediate.

\begin{lemma}\label{l.VMRTsection}
In the setting of Proposition \ref{p.nonsingular}, let $X \cap H$
be a general hyperplane section.  If $\sC_x \subset \BP T_x(X)$ is
the VMRT of $X$ at a general point $x \in X\cap H$ and $\dim \sC_x
\geq 1$, then the VMRT associated to
 a family of lines covering $X \cap H$ is $$\sC_x \cap \BP T_x(H) \subset \BP T_x(X \cap H).$$ \end{lemma}

\subsection{VMRT of an irreducible Hermitian symmetric space
of compact type}\label{e.IHSS}

An irreducible Hermitian symmetric space of compact type is a
homogeneous space $M= G/P$ with a simple Lie group $G$ and a
maximal parabolic subgroup $P$ such that the isotropy
representation of $P$ on $T_x(M)$ at a base point $x \in M$ is
irreducible.
 The highest weight orbit of the isotropy action on $\BP T_x(M)$
is exactly the VMRT at $x$.

The following table collects some well-known facts on irreducible
Hermitian symmetric spaces of compact type (see e.g. \cite{HK}
Section 6.2).

\begin{center}
%\begin{tabular*}{0.908\textwidth}{|c| c| c| c| c| c|}
\begin{tabular}{|c| c| c| c| c| c|}
\hline Type & I.H.S.S. $M$   &  VMRT  $S$ &  $S \subset \BP
T_x(M)$ & $\dim \BP T_x(M)$  & $\dim \Sec(S)$
\\ \hline I    &  $ \Gr(a, a+b) $ & $\pit^{a-1} \times \pit^{b-1}$
& Segre & $ab-1$ & $2a+2b-5$ \\  \hline II  & $\mathbb{S}_{n}$ & $\Gr(2, n)$  & Pl\"ucker & $\frac{1}{2}(n^2-n-2)$ & $4n-11$ \\
\hline III & $ \Lag(2n)$ & $\pit^{n-1}$ & Veronese &
$\frac{1}{2}(n^2+n-2)$ & $2n-2$\\  \hline
 IV & $\Q^n$ & $\Q^{n-2}$ & Hyperquadric & $n-1$ & $n-1$ \\  \hline
 V & $\mathbb{O}\pit^2$ & $\mathbb{S}_{5}$ & Spinor & 15 & $15$ \\  \hline
 VI & $E_7/(E_6 \times U(1)) $ & $\mathbb{O}\pit^2$ & Severi & 26
 & $25$ \\
 \hline
\end{tabular}
\end{center}

Here $\Gr(a, a+b)$ is the Grassmannian of $a$-dimensional
subspaces in an $(a+b)$-dimensional vector space, $\mathbb{S}_{n}$
is the spinor variety, i.e. the variety parameterizing
$n$-dimensional isotropic linear subspaces in an orthogonal vector
space of dimension $2n$. $\Lag(2n)$ is the Lagrangian
Grassmannian, which parameterizes Lagrangian subspaces in a
$2n$-dimensional symplectic vector space. $\Q^n$ denotes the
$n$-dimensional hyperquadric. $\mathbb{O}\pit^2$ is the Cayley
plane, which is of dimension 16 and homogeneous under the action
of $E_6$.

In the following, we always assume that $M$ is not a projective
space. Let  $o \in M=G/P$ be the point with isotropy $P$ and set
$V= T_o(M)$. Let $S \subset \BP V$ be the VMRT of $M$. There
exists a depth 1 decomposition  of the Lie algebra $\fg$ of $G$
(cf. \cite{HM05} Section (4.1)): if we denote by $\alpha_k$ the
simple root corresponding to the maximal parabolic subgroup $P$,
and $\Phi_i$ the set of roots whose coefficient in $\alpha_k$
equals to $i$, then $\Phi_i$ is not empty exactly for $i=0, \pm
1$. Let $\fg_i = \oplus_{\alpha \in \Phi_i} \fg_\alpha$, then we
get the decomposition $\fg= \fg_{-1} \oplus \fg_0 \oplus \fg_1$
satisfying $[\fg_i, \fg_j] \subset \fg_{i+j}$ for all $i, j$.
There exist natural isomorphisms $\fg_{-1} \cong T_o(G/P) = V$,
$\fg_1 \cong T_o^*(G/P) = V^*$ and $\fg_0 \cong \aut(\hat{S})$.
When $G$ is of classical type, the gradation of $\fg$ can be found
in Section 4.4 of \cite{Ya}.

We have a natural injective map:
$$
\phi: \fg_1 \to \Hom(\fg_{-1}, \fg_0), \mbox{ given by } \phi_X(Y)
= [X, Y], \forall X \in \fg_1, Y \in \fg_{-1}.
$$
By Theorem 5.2 of \cite{Ya}, the image $\Im(\phi)$ is exactly the
prolongation $\fg_0^{(1)} = \aut(\hat{S})^{(1)}.$ This gives
\begin{proposition}\label{p.IHSSProlong}
Let $S \subsetneq \BP V$ be the VMRT of an irreducible Hermitian
symmetric space $M$ of compact type. Then
 $\aut(\hat{S})^{(1)} \cong \fg_1
\cong V^*$.

%{\em Type I.} $S = \BP A \times \BP B \hookrightarrow \BP A
%\otimes B$, and $\aut(\hat{S})^{(1)} \cong \Hom(B, A)$.

%{\em Type II.} $S = \Gr(2, W) \hookrightarrow \BP(\wedge^2 W)$ and
%$\aut(\hat{S})^{(1)} \cong \wedge^2 W$
\end{proposition}

As explained in Corollary 1.1.5 of \cite{HM05}, Theorem \ref{t.HM}
implies the following result of \cite{KN}.

\begin{theorem}\label{t.KN}
Let $S \subsetneq \BP V$ be the highest weight variety of an
irreducible representation. Then $\aut(\hat{S})^{(1)} \neq 0$ if
and only if $S$ is the VMRT of an irreducible Hermitian symmetric
space of compact type. \end{theorem}

\subsection{VMRT of symplectic Grassmanians} \label{e.SymGrass}
Let $\Sigma$ be an $n$-dimensional vector space endowed with a
skew-symmetric 2-form $\omega$ of maximal rank. We denote by
$\Gr_\omega(k, \Sigma)$ the variety of all $k$-dimensional
isotropic subspaces of $\Sigma$. When $n$ is even, this is the
usual symplectic Grassmanian, which is homogeneous under the
action of ${\rm Sp}(\Sigma)$.  When $n$ is odd, $\Gr_\omega(k,
\Sigma)$ is the odd symplectic Grassmanian, which is not
homogeneous and it has two orbits under the action of its
automorphism group ${\rm PSp}(\Sigma):=\{g\in {\rm
PGL}(\Sigma)|g^*\omega=\omega\}.$

  Let $W$ and $Q$ be vector spaces of
dimensions $k \geq 2$ and $m$ respectively. Let $ v_2: \pit(W
\oplus Q) \hookrightarrow \pit(\SYM^2(W \oplus Q))$ be the second
Veronese embedding. Let $$ U:=(W \otimes Q) \oplus \SYM^2W \subset
\SYM^2(W \oplus Q).$$ Let $p_{\SYM^2Q}: \pit(\SYM^2(W \oplus Q))
\dasharrow \pit U$ be the projection from $\pit(\SYM^2Q)$. We
denote by $Z$ the proper image of $v_2(\BP(W \oplus Q))$ under the
projection $p_{\SYM^2Q}$. Then $Z$ is isomorphic to the projective
bundle $\pit((Q \otimes \mathbf{t}) \oplus \mathbf{t}^{\otimes
2})$ over $\pit W$, where $\mathbf{t}$ is the tautological line
bundle over $\pit W$. The embedding $Z \subset \pit U$ is given by
the complete linear system $$H^0(\pit W, (Q \otimes \mathbf{t}^*)
\oplus (\mathbf{t}^*)^{\otimes 2}) = (W \otimes Q)^* \oplus \SYM^2
W^* = U^*.$$

%The following diagram summarizes the notations.
%\begin{CD}
%W = \pit(U \oplus Q) @>>{Veronese}> \pit(\SYM^2(U \oplus Q)) \\
%\end{CD}

The following lemma was proved in Proposition 3.2.1 [HM05] for the
case of (even) symplectic Grassmanians. The proof there works also
for odd symplectic Grassmannians.
\begin{lemma}
 The linearly normal embedding $Z
\hookrightarrow \pit U$ is isomorphic to the VMRT of the
symplectic Grassmannian $\Gr_\omega(k, \Sigma)$ (with $\dim \Sigma
= m+2k)$.
\end{lemma}

We also have
\begin{lemma}\label{l.sectionGr}
If $k=2$, then  $Z \hookrightarrow \pit U$ is the VMRT of a
general hyperplane section of the Pl\"ucker embedding of $\Gr(2,
m+4)$.
\end{lemma}

\begin{proof}
Let $\Sigma$ be a vector space of dimension $m+4$, then we have
the Pl\"ucker embedding $\Gr(2, \Sigma) \hookrightarrow
\pit(\wedge^2 \Sigma)$. Let $\omega$ be a general element of
$\wedge^2 \Sigma^*$, i.e.,
 a skew-symmetric 2-form on $\Sigma$ with maximal rank.
 Let $H \subset \wedge^2 \Sigma$ be
the kernel of $\omega \in \wedge^2 \Sigma^*$. Then we get $\Gr(2,
\Sigma) \cap H = \Gr_{\omega}(2, \Sigma)$, the latter being the
symplectic Grassmannian.
\end{proof}

The following will be proved after Proposition \ref{p.autSymp}.

\begin{proposition}\label{p.autZ} $\aut(\hat{Z}) = (W^* \otimes Q) \sd
(\fgl(W) \oplus \fgl(Q))$ and $\aut(\hat{Z})^{(1)} \cong \SYM^2
W^*.$  \end{proposition}

\subsection{Hyperplane section of $\mathbb{S}_{5}$}\label{e.S_5}

Let $Q$ be a 7-dimensional orthogonal vector space  and let $W$ be
the 8-dimensional spin representation of $\fso(Q) = \fso(7)$.
There exists a ${\rm Spin}(7)$-stable  $9$-dimensional Fano
manifold $Z$ of Picard number 1 with an embedding $Z \subset \BP
(W \oplus Q)$ which is isomorphic to a general hyperplane section
of the 10-dimensional spinor variety (cf. Section 7 in \cite{HM05}
where it is denoted by $\sC_o$). In fact, as explained in Section
7 of \cite{HM05}, $Z \subset \BP (W \oplus Q)$ is isomorphic to
the VMRT of a 15-dimensional $F_4$-homogeneous space.
 The variety $Z$ is biregular  to the
horospherical Fano manifold of Picard number 1, the case 4 in
Theorem 1.7 of \cite{Pa}. The next proposition follows from
Theorem 1.11 of \cite{Pa}.

\begin{proposition}\label{p.autZQ} $\aut(\hat{Z}) = \C \oplus W \sd (\fso(Q) \oplus \C).$ \end{proposition}

Here the center $\C$ corresponds to the scalar multiplication on
$W \oplus Q$, while the second $\C$ acts with  weight 1 on $W$ and
0 on $Q$. The action of $W$ on $W \oplus Q$ is annihilating $W$
and given by $W \subset {\rm Hom} (Q,W)$ induced from  the natural
inclusion of $W$ as an irreducible $\fso(7)$-factor of ${\rm
Hom}(Q, W)$. The inclusion $\aut(\hat{Z}) \subset \End(W \oplus
Q)$ can be represented as follows:
$$
\begin{pmatrix}
\mathbb{C} \subset \End(W)  &  0 \\
W \subset \Hom(Q, W) & \mathfrak{co}(Q) \subset \End(Q)
\end{pmatrix}.
$$

 The following is from Proposition 7.2.3 of
\cite{HM05}. We give a more direct proof.

\begin{proposition}\label{p.prolongZQ} $\aut(\hat{Z})^{(1)} =
Q^*.$
\end{proposition}
\begin{proof}
By Theorem \ref{t.HM} (ii) and (iii),  every $A \in
\aut(\hat{Z})^{(1)}$ is determined by an element $\lambda \in W^*
\oplus Q^*$ such that $$2 A_{x, y} = \lambda(x) y + \lambda(y) x$$
for $x \in \hat{Z}$ and $y \in T_x(\hat{Z})$. As $A_x \in
\aut(\hat{Z})$, we can write
$$
A_x = \begin{pmatrix} \mu_x & 0 \\ \phi_x & g_x \end{pmatrix},
\quad \mu_x \in \mathbb{C}, \phi_x \in W \subset \Hom(Q, W), g_x
\in \mathfrak{co}(Q).
$$
If we write $x=(x_1, x_2)$ and $y=(y_1, y_2)$ with $x_1, y_1 \in
W$ and $x_2, y_2 \in Q$, then we have $$2(\mu_x y_1 + \phi_x(y_2),
g_x(y_2)) = 2A_x(y)= 2A_{x,y} = (\lambda(x) y_1, \lambda(x) y_2) +
(\lambda(y)x_1, \lambda(y) x_2).$$
 As this holds
for all $y\in T_x(\hat{Z})$, we have
\begin{equation} \label{equa1}
\mu_x = \lambda(x)/2,  \phi_x(y_2) = (\lambda(y)/2) x_1, g_x(y_2)
= (\lambda(x)/2) y_2 + (\lambda(y)/2) x_2.
\end{equation}
If we take $y_2=0$ in the previous equations, then $\lambda((y_1,
0)) = 0$ for all $y_1 \in W$, which implies $\lambda \in Q^*$.
Conversely, for any $\lambda \in Q^*$, we can use formulae in
\eqref{equa1} to construct $A_x$ and one checks that $A \in
\aut(\hat{Z})^{(1)}$.
\end{proof}

\subsection{Hyperplane section of $\Gr(2, 5)$}\label{e.Gr25}
Let $Q$ be a 5-dimensional orthogonal vector space  and let $W$ be
the 4-dimensional spin representation of $\fso(Q) = \fso(5)$.
There exists a ${\rm Spin}(5)$-stable  $5$-dimensional Fano
manifold $Z$ of Picard number 1 with an embedding $Z \subset \BP
(W \oplus Q).$ In fact, $Z$ is a general hyperplane section of
$\Gr(2,5)$, which is isomorphic to  a symplectic Grassmannian
$\Gr_{\omega}(2,5)$ by Lemma \ref{l.sectionGr}. This can be seen
as follows. As $\fsp(4) \cong \fso(5)$, we can regard $W$ as a
4-dimensional symplectic vector space.  We have an
$\fsp(4)$-module decomposition
$$\wedge^2 W \cong \C \oplus Q.$$
Equip $W \oplus \C$ with the skew symmetric form $\omega$ obtained
from $W$ with $\C$ as its null-space. A natural embedding $Z
\subset \BP(W \oplus Q)$ of the symplectic Grassmannian $Z$ can be
obtained by viewing $Z$ as the hyperplane section of $$\Gr(2, W
\oplus \C) \subset \BP \wedge^2(W \oplus \C) \cong \BP(\wedge^2 W
\oplus W)$$ where the hyperplane is given by the kernel of
$\omega$ in $\wedge^2 W = \C \oplus Q$. From the table in Section
\ref{e.IHSS}, we see that $Z$ is the VMRT at a general point of a
hyperplane section of $\mathbb{S}_5$. As Proposition
\ref{p.autZQ}, the next proposition follows from Theorem 1.11 Case
5 of \cite{Pa}.

\begin{proposition}\label{p.autZQ1} $\aut(\hat{Z}) = \C \oplus W \sd (\fso(Q) \oplus \C).$ \end{proposition}

The next proposition can be proved in the same way as Proposition
\ref{p.prolongZQ}.

\begin{proposition}\label{p.prolongZQ2} $\aut(\hat{Z})^{(1)} =
Q^*.$
\end{proposition}

\section{Prolongation and projection}\label{s.prolong}

Given a linear space $L \subset V$, denote by $p_L: \BP V
\dasharrow \BP (V/L)$ the projection.
 In this section, we study the
prolongation of $p_L(S) \subset \BP (V/L)$ for the examples $S
\subset \BP V$ listed in the previous section for suitable linear
spaces $L$.

Let us recall the following two elementary facts.

\begin{lemma}\label{l.sec}
 Given an irreducible variety $S
\subset \BP V$ and a linear subspace $L \subset V$ with $S
\not\subset \BP L$, the proper image $p_L(S) \subset \BP V/L$ is
well-defined. When $S$ is nonsingular, the restriction $p_L|_S$ is
a morphism sending $S$ biregularly to $p_L(S)$ if and only if $\BP
L \cap \Sec(S) = \emptyset$.
\end{lemma}

\begin{lemma}\label{l.secant}
Let $S \subset \pit V$ be an irreducible closed subvariety and $L
\subset  V$ a linear subspace with $S \not\subset \BP L$. Then
$\Sec(p_L(S)) = p_L(\Sec(S))$.
\end{lemma}

To study the prolongation of $p_L(S)$, it is convenient to
introduce the following.

\begin{definition}\label{d.L} Let $S \subset \BP V$ be an irreducible projective variety and let $L \subset V$ be a
linear subspace. We define two Lie subalgebras of $\aut(\hat{S})$
as follows.
$$\aut(\hat{S}, L):=\{g \in \aut(\hat{S}) | g(L)\subset L\} \subset \mathfrak{gl}(V)$$ $$\aut(\hat{S}, L,
0):=\{g \in \aut(\hat{S}) | g(L)=0\} \subset \mathfrak{gl}(V).$$
\end{definition}

\begin{proposition}\label{p.proj}
Let $S \subset \BP(V)$ be a non-degenerate irreducible subvariety.
Let $L \subset V$ be a linear subspace with $S \not\subset \BP L$.
Assume that
 \begin{enumerate}
 \item[(i)] the natural Lie algebra homomorphism $\aut(\hat{S}, L) \to \aut(\widehat{p_L(S)})$ is an
 isomorphism; and \item[(ii)] for a general $\alpha \in \hat{S}$, $T_{\alpha}(\hat{S}) \cap L = 0$. \end{enumerate}
Then we have an isomorphism of vector spaces
$$
\aut(\widehat{p_L(S)})^{(1)} \cong \aut(\hat{S}, L, 0)^{(1)} = \{
A \in \aut(\hat{S})^{(1)} | A_\alpha(L)=0, \forall \alpha \in
\hat{S} \}.
$$
\end{proposition}
\begin{proof}
 For any element $A \in \aut(\widehat{p_L(S)})^{(1)}
\subset \Hom(V/L, \aut(\widehat{p_L(S)}))$, we define an element
$\tilde{A} \in \Hom(V, \aut(\hat{S}, L))$ by composing $A$ with
the natural projection $V \to V/L$ and the isomorphism
$\aut(\widehat{p_L(S)}) \cong \aut(\hat{S}, L)$ given by the
condition (i).

 For a general element $\alpha \in \hat{S}$, denote by $\bar{\alpha} \in \widehat{p_L(S)}$
 its image in $V/L$.
 By the condition (ii), we have a natural identification $T_{\alpha}(\hat{S}) \cong T_{\bar{\alpha}}(\widehat{p_L(S)})$
 for a general $\alpha \in \hat{S}$. Let $\beta \in \hat{S}$ be a
 general point. As $A_{\bar{\alpha}}(\bar{\beta}) = A_{\bar{\beta}}(\bar{\alpha}) \in T_{\bar{\alpha}}(\widehat{p_L(S)})
 \cap T_{\bar{\beta}}(\widehat{p_L(S)}),$  we have $\tilde{A}_\alpha (\beta) \in
  (T_\alpha(\hat{S}) \cap T_\beta (\hat{S})) \oplus L$. On the
  other hand, as $\tilde{A}_\alpha \in \aut(\hat{S})$, we have
  $\tilde{A}_\alpha(\beta) \in T_\beta(\hat{S})$ by Definition \ref{d.basic} (v).
    As $T_\beta(\hat{S}) \cap L =0$, this implies that
   $\tilde{A}_\alpha(\beta) \in T_\alpha(\hat{S}) \cap T_\beta
  (\hat{S})$. In particular, we have $\Im(\tilde{A}_\alpha) \subset T_\alpha
  (\hat{S})$ and $\tilde{A}_\alpha(\beta)=
\tilde{A}_\beta(\alpha)$ for all general $\alpha, \beta\in
\hat{S}$  because $A_{\bar{\alpha}}(\bar{\beta}) =
A_{\bar{\beta}}(\bar{\alpha}).$ As $S$ is non-degenerate, the
equality $\tilde{A}_\alpha(\beta)= \tilde{A}_\beta(\alpha)$ holds
for all $\alpha, \beta\in V$. Thus $\tilde{A} \in \Hom(\SYM^2V,
V)$.

  As $\tilde{A}_\alpha \in \aut(\hat{S}, L)$, we have
$\tilde{A}_\alpha(L) \subset L$. This implies that
$\tilde{A}_\alpha(L) \subset L \cap T_\alpha (\hat{S}) = 0$ for
general $\alpha \in \hat{S}$ by the condition (ii). Consequently,
$\tilde{A}_\alpha(L)=0$ for a general $\alpha \in \hat{S}$, hence
for any $\alpha \in V$ by the non-degeneracy of $S$. It follows
that
$$\tilde{A} \in \Hom(V, \aut(\hat{S}, L,0)) \cap \Hom(\SYM^2V, V)
=\aut(\hat{S}, L, 0)^{(1)}.$$  Conversely, for any $\tilde{A} \in
\aut(\hat{S}, L, 0)^{(1)}$, it is easy to see that it induces an
element $A \in \aut(\widehat{p_L(S)})^{(1)}$, proving the
proposition.
\end{proof}

\begin{proposition}\label{p.aut}
Let $S \subset \BP V$ be a linearly normal  nonsingular
non-degenerate projective variety. If $L \subset V$ is a subspace
with $\BP L \cap \Sec(S) = \emptyset$, then it satisfies the two
conditions in Proposition \ref{p.proj}. In particular,  $$
\aut(\widehat{p_L(S)})^{(1)} \cong \aut(\hat{S}, L, 0)^{(1)} =  \{
A \in \aut(\hat{S})^{(1)} | A_\alpha(L)=0, \forall \alpha \in
\hat{S} \}.
$$ \end{proposition}

\begin{proof}
Let us identify $S \subset \BP V$ with $S \subset \BP H^0(S,
\sO(1))^*$ for the hyperplane line bundle $\sO(1)$ on $S$. The
condition (ii) is immediate from $\BP T_{\alpha}(\hat{S}) \subset
\Sec(S)$ for any $\alpha \in \hat{S}$. The condition (i) will
follow from the next lemma. \end{proof}

\begin{lemma}\label{l.Aut} Let $S \subset \BP V$ be a linearly normal nonsingular  non-degenerate projective
variety. Let $\BP L_1, \BP L_2 \subset \BP V \setminus {\rm
Sec}(S)$ be two linear subspaces and $p_{L_i}: S \to p_{L_i}(S)
\subset \BP(V/L_i)$ the projection from $\BP L_i$. Suppose that
there exists an isomorphism $\sigma: \BP(V/L_1) \to \BP(V/L_2)$
with $\sigma( p_{L_1}(S)) = p_{L_2}(S)$. Then there exists a
unique isomorphism $\tilde{\sigma}: \BP V \to \BP V$ with $\tilde{
\sigma}(S) = S$ and $\tilde{\sigma}(\BP L_1) = \BP L_2$.
\end{lemma}

\begin{proof}
The restriction $\sigma|_{p_{L_1}(S)}: p_{L_1}(S) \to p_{L_2}(S)$
is an automorphism $\bar{\sigma}$ of $S$ with
$\bar{\sigma}^*\sO(1) \cong \sO(1)$ such that sections of $\sO(1)$
annihilated by $\BP L_2 \subset \BP H^0(S, \sO(1))^*$ correspond
to sections of $\sO(1)$ annihilated by $\BP L_1 \subset \BP H^0(S,
\sO(1))^*$. Thus it induces a homomorphism $\tilde{\sigma}: \BP
H^0(S, \sO(1))^* \to \BP H^0(S, \sO(1))^*$ with
$\tilde{\sigma}(\BP L_1) = \BP L_2$.
\end{proof}

By  Propositions \ref{p.proj} and \ref{p.aut}, studying the
prolongation of $p_L(S)$ under a biregular projection of a
linearly normal $S$ is reduced to the study of $\aut(\hat{S}, L,
0)^{(1)}.$ Let us carry this out for the examples listed in
Section \ref{s.examples}.

For the VMRT of an irreducible Hermitian symmetric space of
compact type, we have the following uniform description.

\begin{proposition}\label{p.projIHSS}
Let $S \subsetneq \BP V$ be the VMRT of an irreducible Hermitian
symmetric space $M$ of compact type.  Recall that from Section
\ref{e.IHSS}, we have a graded Lie algebra structure of
$\fg:=\aut(M)$, $\fg= \fg_{-1} \oplus \fg_0 \oplus \fg_1$ such
that $\fg_{-1} \cong T_o(M)=V$ and $\aut(\hat{S})^{(1)} \cong
\fg_1$. For a subspace $L \subset \fg_{-1}$, we have
$$\aut(\hat{S}, L, 0)^{(1)} = \{ X \in \fg_1 | [X, Z] = 0, \forall
Z \in L \}.$$
\end{proposition}

\begin{proof}
From the definition $$\aut(\hat{S}, L, 0)^{(1)} = \{ A \in
\aut(\hat{S})^{(1)} | A_\alpha(L)=0, \forall \alpha \in \hat{S}
\}$$ and the isomorphism $\aut(\hat{S})^{(1)} \cong \fg_1$ given
in Section \ref{e.IHSS}, we get
\[
\begin{aligned}
\aut(\hat{S}, L, 0)^{(1)} & = \{ X \in \fg_1 | [[X,Y],Z] = 0,
\forall Y \in \fg_{-1}, \forall Z \in L \} \\
&=\{ X \in \fg_1 | [[X,Z],Y] = 0, \forall Y \in \fg_{-1}, \forall
Z \in L \} \\
 & =\{ X
\in \fg_1 | [X, Z] = 0, \forall Z \in L \}.
\end{aligned}
\]
The last equality follows from the fact that if an element $u \in
\fg_{1}=\fg_0^{(1)}$ satisfies $[u, Y] = 0$ for all $Y \in
\fg_{-1}$, then $u=0$ (cf. \cite{Ya} Lemma 3.2). \end{proof}

 The following four propositions give more explicit description of
 Proposition \ref{p.projIHSS} when $\Sec(S) \neq \BP V$. Here we
 will use the data in the table of Section \ref{e.IHSS} freely. Our main interest
 is when $\BP L \cap \Sec(S)= \emptyset.$ But for the classical types, we will treat
 also general $L$, because it will be needed later and requires
 little extra work.

\begin{proposition} \label{p.I}
Let $A$ and $B$ be vector spaces with $a:=\dim A \geq b:=\dim B
\geq 3$. Let $V = \Hom (A, B)$ and let $\hat{S} \subset V$ be the
set of elements of rank $\leq 1$. For a subspace $L \subset
\Hom(A, B)$, we define $\Im(L) \subset B$ as the linear span of
$\{ \Im(\phi) \subset B, \phi \in L \} $ and $\Ker(L) :=
\cap_{\phi \in L} \Ker(\phi)$. Then \begin{enumerate} \item[(i)]
there is a canonical isomorphism of vector spaces
 $$\aut(\hat{S}, L, 0)^{(1)} \cong \Hom(B/\Im(L), \Ker(L));$$ \item[(ii)] for any $\psi \in \Hom(B, A)$
 with $\Im(L) \subset \Ker(\psi)$ and $\Im(\psi) \subset \Ker(L)$,  $L$ is
  contained in $$L(\psi):=\{ \phi \in \Hom(A, B) | \phi \circ \psi=0, \psi \circ \phi=0\}  \cong \Hom(A/\Im(\psi), \Ker(\psi));$$ \item[(iii)]
 if $\BP L \cap \Sec(S) = \emptyset$ and $\aut(\hat{S}, L, 0)^{(1)}$ contains an element of
 rank $r$ in $\Hom(B/\Im(L), \Ker(L))$,
 $$\dim L \leq ab-2(a+b)+4 - r(a+b-r-4).$$ \end{enumerate}
\end{proposition}

\begin{proof}
From \cite{Ya} p.457, the grading on $\fg$ in Proposition
\ref{p.projIHSS} for $M= \Gr(a, a+b)$ can be identified with
$$
\fg_{-1} = \Hom(A, B), \fg_0 = \End(A) \oplus \End(B), \fg_1 =
\Hom(B, A).
$$
The bracket $[\fg_{-1}, \fg_1] \subset \fg_0$ is given by $[\phi,
\psi] = \phi \circ \psi-\psi \circ \phi   \in \fg_0$ for $\phi \in
\fg_{-1}$ and $\psi \in \fg_1$. By Proposition \ref{p.projIHSS},
this gives
\[
\begin{aligned}
\aut(\hat{S}, L, 0)^{(1)} &= \{ \psi \in \Hom(B, A) | \psi \circ \phi =0, \phi \circ \psi=0, \forall \phi \in L\}\\
&= \{ \psi \in \Hom(B, A) | \Im(\psi) \subset \Ker(\phi),
\Im(\phi) \subset \Ker(\psi), \forall \phi \in L
\} \\
&=  \{ \psi \in \Hom(B, A) | \Im(\psi) \subset \Ker(L), \Im(L)
\subset \Ker(\psi) \} \\
&\cong \Hom(B/\Im(L), \Ker(L)),
\end{aligned}
\] proving (i).
For (ii), it is clear from above that   $L$ is contained in
$L(\psi)$.
% := \{ \phi \in \Hom(A, B) | [\phi, \psi]=0\} =\{ \phi
%\in \Hom(A, B) | \Im(\phi) \subset \Ker(\psi),
%        \Im(\psi) \subset \Ker(\phi)\}. $$

Now assume that $\BP L \cap \Sec(S) = \emptyset$ and there is
$\psi$  in (ii) of rank $r$.  Note that $ \Sec(\hat{S}) \subset V
= \Hom(A,B)$ consists of elements of rank $\leq 2$  (e.g.
\cite{HK} p.188 Type I). Let $S_\psi \subset \BP L(\psi) \cong
\BP(\Hom(A/\Im(\psi), \Ker(\psi)))$ be the set of elements of rank
$\leq 1$, then $\Sec(S_\psi)$ consists of elements of rank $\leq
2$, which has dimension $2(a+b-2r)-5$.
 By the assumption, $\BP
L \subset \BP L(\psi)$ is disjoint from $\Sec(S_\psi)$, which
implies that $a-r \geq b-r \geq 3$ and
$$
\dim L \leq (a-r)(b-r) -(\dim \Sec(S_\psi)) -1 = ab-2(a+b)+4 -
r(a+b-r-4).
$$
\end{proof}

\begin{proposition} \label{p.II}
Let $W$ be a vector spaces of dimension $n\geq 6$. For each $\phi
\in \wedge^2 W$, denote by $\phi^{\sharp} \in \Hom( W^*, W)$, the
corresponding element via the natural inclusion $\wedge^2 W
\subset W \otimes W = \Hom(W^*, W)$.  Let $V = \wedge^2 W$ and let
$\hat{S} \subset V$ be the set of elements $\phi$ with
$\phi^{\sharp}$ of rank $\leq 2$. For a subspace $L \subset V$,
define $\Im(L) \subset W$ as the linear span of $\{
\Im(\phi^{\sharp}) \subset W, \phi \in L \}.$ Then
\begin{enumerate} \item[(i)]
there is a canonical vector space isomorphism
 $$\aut(\hat{S}, L, 0)^{(1)} \cong \wedge^2 (W/\Im(L))^*;$$
 \item[(ii)] for each $\psi \in \wedge^2 (W/\Im(L))^* \subset \wedge^2 W^*,$
 denoting by $\psi^{\sharp}$ the corresponding element in $\Hom(W, W^*)$,
  $L$ is contained in $$L(\psi):= \{ \phi \in \wedge^2 W | \Im(\phi^{\sharp}) \subset \Ker(\psi^{\sharp})\}  \cong \wedge^2 \Ker(\psi^{\sharp});$$
  \item[(iii)] if $\BP L \cap \Sec(S) = \emptyset$ and $\aut(\hat{S}, L, 0)^{(1)}$ contains an element of
 rank $r$ in $\wedge^2 (W/\Im(L))^*$ (i.e. the corresponding
 element in $\Hom(W, W^*)$ has rank $r$), $$\dim L
\leq \cfrac{n(n-1)}{2}-4n+10 - \cfrac{r(2n-r-9)}{2}.$$
\end{enumerate}
\end{proposition}

\begin{proof}
From \cite{Ya} pp. 459-461, the grading on $\fg$ in Proposition
\ref{p.projIHSS} for $M= \mathbb{S}_n$ can be identified with
$$
\fg_{-1} = \wedge^2 W, \fg_0 = \End(W), \fg_1 = \wedge^2 W^*.
$$ For each $\psi
\in \wedge^2 W^*$, denote by $\psi^{\sharp} \in \Hom( W, W^*)$,
the corresponding element via the natural inclusion $\wedge^2 W^*
\subset W^* \otimes W^* = \Hom(W, W^*)$.

For any $\phi \in \wedge^2 W, \psi \in \wedge^2 W^*$, the
endomorphism $[\phi, \psi] \in \End(W)$ is given   by
$$
[\phi, \psi] = \phi^{\sharp} \circ \psi^{\sharp}.
$$
Note that we have the following equivalences:
\[
\begin{aligned}
\quad  [\phi, \psi]=0  & \Leftrightarrow \Im(\psi^{\sharp})
\subset \Ker (\phi^{\sharp})
 \Leftrightarrow \psi \in \wedge^2 \Ker(\phi^{\sharp}) \\ & \Leftrightarrow \phi \in \wedge^2 \Ker(\psi^{\sharp})
 \Leftrightarrow \Im(\phi^{\sharp}) \subset \Ker(\psi^{\sharp}).
\end{aligned}
\]
By Proposition \ref{p.projIHSS}, this gives
\[
\begin{aligned}
\aut(\hat{S},L,0)^{(1)} &= \{ \psi \in \wedge^2 W^* |
\Im(\phi^{\sharp}) \subset \Ker(\psi^{\sharp}), \forall \phi \in L
\} \\
&=  \{ \psi \in \wedge^2 W^* | \Im(L) \subset \Ker(\psi^{\sharp})\} \\
&\cong \wedge^2 (W/\Im(L))^*,
\end{aligned}
\] proving (i).
For (ii), it is clear from above that   $L$ is contained in
$L(\psi)$.
%  := \{ \phi \in \wedge^2 W | [\phi, \psi]=0\} = \{ \phi
%\in \wedge^2 W | \Im(\phi^{\sharp}) \subset
%\Ker(\psi^{\sharp})\}.$$

Now assume that $\BP L \cap \Sec(S) = \emptyset$ and there is
$\psi$  in (ii) of rank $r$.  Note that $ \Sec(\hat{S}) \subset V
= \wedge^2 W$ consists of elements with rank $\leq 4$ (e.g.
\cite{HK} p.188 Type II).
 Let
$S_\psi \subset \BP(L(\psi))$ be the variety consisting of
elements of rank $\leq 2$, then we have $\dim \Sec(S_\psi) =
4n-4r-11.$
 By the hypothesis, $\BP
L \subset \BP L(\psi)$ is disjoint from $\Sec(S_\psi)$, which
implies that $n-r \geq 6$ and
$$\dim L
\leq \cfrac{n(n-1)}{2}-4n+10 - \cfrac{r(2n-r-9)}{2}.$$
\end{proof}

\begin{proposition} \label{p.III}
Let $W$ be a vector space of dimension  $n\geq 3$. For each $\phi
\in \Sym^2 W$, denote by $\phi^{\sharp} \in \Hom( W^*, W)$, the
corresponding element via the natural inclusion $\Sym^2 W \subset
W \otimes W = \Hom(W^*, W)$.  Let $V = \Sym^2 W$ and let $\hat{S}
\subset V$ be the set of elements $\phi$ with $\phi^{\sharp}$ of
rank $\leq 1$. For a subspace $L \subset V$, define $\Im(L)
\subset W$ as the linear span of $\{ \Im(\phi^{\sharp}) \subset W,
\phi \in L \}.$ Then
\begin{enumerate} \item[(i)] there is a canonical isomorphism of vector
spaces

 $$\aut(\hat{S}, L, 0)^{(1)} \cong \Sym^2 (W/\Im(L))^*;$$ \item[(ii)] for each
 $\psi \in \Sym^2 (W/\Im(L))^* \subset \Sym^2 W^*,$ denoting by $\psi^{\sharp}$ the
 corresponding element in $\Hom(W, W^*)$,  $L$ is contained in $$L(\psi):= \{ \phi \in \Sym^2 W | \Im(\phi^{\sharp}) \subset \Ker(\psi^{\sharp})\} \cong \Sym^2 \Ker(\psi^{\sharp});$$
    \item[(iii)] if $\BP L \cap \Sec(S) = \emptyset$ and $\aut(\hat{S}, L, 0)^{(1)}$ contains an element of
 rank $r$ in $\Sym^2 (W/\Im(L))^*$ (i.e. the corresponding
 element in $\Hom(W, W^*)$ has rank $r$),
 $$\dim L
\leq \cfrac{n(n+1)}{2}-2n+1 - \cfrac{r(2n-r-3)}{2}.$$
\end{enumerate}
\end{proposition}

\begin{proof}
From \cite{Ya} pp. 458-459, the grading on $\fg$ in Proposition
\ref{p.projIHSS} for $M= \Lag(2n)$ can be identified with
$$
\fg_{-1} = \Sym^2 W, \fg_0=\End(W), \fg_1 = \Sym^2 W^*.
$$
For each $\psi \in \Sym^2 W^*$, denote by $\psi^{\sharp} \in \Hom(
W, W^*)$, the corresponding element via the natural inclusion
$\Sym^2 W^* \subset W^* \otimes W^* = \Hom(W, W^*)$. For any $\phi
\in \Sym^2 W, \psi \in \Sym^2 W^*$,  the endomorphism $[\phi,
\psi] \in \End(W)$ is given by
$$
[\phi, \psi] = \phi^{\sharp} \circ \psi^{\sharp}.
$$

As in the proof of Proposition \ref{p.II}, $[\phi, \psi]=0$ is
equivalent to $\Im(\phi^{\sharp}) \subset \Ker(\psi^{\sharp})$,
which gives, by Proposition \ref{p.projIHSS},
\[
\begin{aligned}
\aut(\hat{S}, L, 0)^{(1)} &= \{ \psi \in \Sym^2 W^* |
\Im(\phi^{\sharp}) \subset \Ker(\psi^{\sharp}), \forall \phi \in L
\} \\
&=  \{ \psi \in \Sym^2 W^* | \Im(L) \subset \Ker(\psi^{\sharp})\} \\
&\cong \Sym^2 (W/\Im(L))^*,
\end{aligned}
\] proving (i).
For (ii), it is clear from above that   $L$ is contained in
$$ L(\psi) = \{ \phi \in \SYM^2 W | [\phi, \psi]=0\} = \{ \phi \in \SYM^2 W | \Im(\phi^{\sharp}) \subset
\Ker(\psi^{\sharp})\}.$$

Now assume that $\BP L \cap \Sec(S) = \emptyset$ and there is
$\psi$  in (ii) of rank $r$. Note that $ \Sec(\hat{S}) \subset V =
\Sym^2 W$ consists of elements with rank $\leq 2$ (e.g. \cite{HK}
p.188 Type III).
 Let
$S_\psi \subset \BP(L(\psi))$ be the variety consisting of
elements of rank $\leq 1$. Then we have $\dim \Sec(S_\psi) =
2(n-r-1).$
 By the hypothesis, $\BP
L \subset \BP L(\psi)$ is disjoint from $\Sec(S_\psi)$, which
implies that $n-r \geq 3$ and
$$\dim L
\leq \cfrac{n(n+1)}{2}-2n+1 - \cfrac{r(2n-r-3)}{2}.$$
\end{proof}

\begin{proposition} \label{p.VI}
Let $ S \subset \BP  V$ be the minimal (Severi) embedding of the
Cayley plane $\mathbb{O}\mathbb{P}^2$ with $\dim V = 27$.   By
\cite{Za} (p. 59-60), $\Sec(S) \subset \BP V$ is a cubic
hypersurface.  For any 1-dimensional subspace $L \subset V$ with
$\BP L \not\in \Sec(S)$, $\aut(\hat{S}, L, 0)^{(1)} =0.$
\end{proposition}

\begin{proof} It is known that $\aut(\hat{S}, L)$ is a simple Lie algebra of type $F_4$ and the
natural representation on $V/L$ is the minimal irreducible
representation of dimension 26 (\cite{Za}, p. 59-60). Let $S'
\subset \BP(V/L)$ be the highest weight variety of this
$F_4$-representation, which is not biregular to the VMRT of an
irreducible Hermitian symmetric space. From Theorem \ref{t.KN}, we
have $\aut(\widehat{p_L(S)})^{(1)} = \aut(\hat{S}')^{(1)} =0.$
Thus by Proposition \ref{p.aut}, $\aut(\hat{S}, L, 0)^{(1)} =
\aut(\widehat{p_L(S)})^{(1)} =0.$
\end{proof}

At this point, we can give the postponed proof of  Proposition
\ref{p.autZ}. We start with examining a special case of
Proposition \ref{p.III}.

\begin{proposition}\label{p.tanSymp}
 Let $W$ and $Q$ be vector spaces of
dimensions $k \geq 2$ and $m$ respectively. Set $L:= \Sym^2Q
\subset V:= \Sym^2(W \oplus Q)$. Let $$ S:= v_2(\pit(W \oplus Q))
\subset \pit(\Sym^2(W \oplus Q))$$ be the second Veronese
embedding of $\pit(W \oplus Q)$. Then for a general point $\alpha
\in \hat{S}$, the tangent space $T_{\alpha}(\hat{S})$ satisfies
$T_{\alpha}(\hat{S}) \cap L =0$. In particular, $\BP L \not\subset
\Sec(S).$ \end{proposition}

\begin{proof}
It suffices to exhibit a point $\alpha \in \hat{S}$ with
$T_{\alpha}(\hat{S}) \cap L =0$. Fix a non-zero $w \in W$ and let
$$\alpha := w^2 \in \Sym^2 W \subset \Sym^2(W \oplus Q).$$
Fix any $w' \in W, q \in Q$.  The arc $$\{ w + t(w'+ q) \in W
\oplus Q | t \in  \C\}$$ in $\hat{S}$ has its tangent vector
$$ \frac{d}{dt}|_{t=0}(w + t(w'+q))^2 =2w(w'+q) \in \Sym^2 W
\oplus (W \otimes Q).$$ Since  such tangent vectors span
$T_{\alpha}(\hat{S})$, $T_{\alpha}(\hat{S})$ intersects $L$ at 0.
\end{proof}

\begin{proposition}\label{p.autSymp}
In the setting of Proposition \ref{p.tanSymp},  let $Z = p_L(S)$
be the proper image of $S$ under $p_L$.  The natural Lie algebra
homomorphism $\aut(\hat{S}, L) \to \aut(\hat{Z})$ is an
 isomorphism, inducing a Lie algebra isomorphism
$$\aut(\hat{Z}) \cong (W^* \otimes Q) \sd (\fgl(W) \oplus
\fgl(Q)).$$ \end{proposition}

\begin{proof} It is clear that $$\aut(\hat{S}, L) \cong
(W^* \otimes Q) \sd (\fgl(W) \oplus \fgl(Q)).$$ The Lie algebra
homomorphism $\aut(\hat{S}, L) \to \aut(\widehat{p_L(S)})$ is
clearly injective. Thus it suffices to show that $\dim
\aut(\hat{Z}) \leq m^2 + km +k^2,$ or equivalently,  $\dim \aut(Z)
\leq m^2 + km +k^2-1.$ From Section \ref{e.SymGrass}, we have a
natural projection $\psi: Z \to \BP W$ realizing $Z$ as the
projectivization of the vector bundle $\sO(-1)^m \oplus \sO(-2)$
on $\BP W \cong \BP^{k-1}.$ From the exact sequence
$$0 \to T^{\psi} \to T(Z) \to \psi^* T(\BP W) \to 0$$ where
$T^{\psi}$ denotes the relative tangent bundle and $$\dim H^0(Z,
T^{\psi}) = \dim H^0(\BP^{k-1}, {\rm End}^0(\sO(-1)^m \oplus
\sO(-2)) = m^2 + km,$$ where ${\rm End}^0$ denotes the traceless
endomorphisms, we have $$\dim \aut(Z) = \dim H^0(Z, T(Z))  \leq
m^2 + km + k^2-1.$$
\end{proof}

\begin{proof}[Proof of Proposition \ref{p.autZ}]
 $\aut(\hat{Z})$ is given by Proposition \ref{p.autSymp}.
From Propositions \ref{p.tanSymp} and \ref{p.autSymp}, we can
apply Proposition \ref{p.proj} to $S$ and $L$. Thus by Proposition
\ref{p.III}, $$\aut(\hat{Z})^{(1)} = \aut(\widehat{p_L(S)})^{(1)}
\cong \aut(\hat{S}, L, 0)^{(1)} \cong \Sym^2 W^*$$ because $\Im(L)
= Q \subset (W \oplus Q)$. \end{proof}

Now we turn to  study the prolongation of the biregular projection
of $Z \subset \BP U$ with $U=(W \otimes Q) \oplus \Sym^2 W$ in
Section \ref{e.SymGrass}. This can be reduced to Proposition
\ref{p.III} by the following.

\begin{proposition}\label{p.double}
Let $S_1 \subset \BP V_1$ be a non-degenerate subvariety.  Let
$L_1 \subset V_1$ be a linear subspace with $S_1 \not\subset \BP
L_1$. Let $V_2:= V_1/L_1$ and let  $S_2 := p_{L_1}(S_1) \subset
\BP V_2$ be the proper image of $S_1$. Let $L_2 \subset V_2$ be a
linear subspace with $S_2 \not\subset \BP L_2$ and let $L_3
\subset V_1$ be the subspace containing $L_1$ with $L_3/L_1 =
L_2$. Suppose that $(S_1, V_1, L_1)$ (resp. $(S_2, V_2, L_2)$)
satisfies the two conditions in Proposition \ref{p.proj} with
$S=S_1, V=V_1, L=L_1$ (resp. with $S=S_2, V= V_2, L=L_2$). Then
$(S_1, V_1, L_3)$ satisfies the two conditions in Proposition
\ref{p.proj} with $S=S_1, V=V_1, L=L_3$. In particular, when the
two conditions are satisfied by $(S_1,V_1,L_1)$ and
$(S_2,V_2,L_2)$, we have
$$\aut(\widehat{p_{L_2}(S_2)})^{(1)} \cong \aut(\hat{S}_1, L_3,
0)^{(1)}.$$ \end{proposition}

\begin{proof}
 For the condition (i) in Proposition \ref{p.proj}, it suffices
to show that the homomorphism $$\aut(\hat{S}_1, L_3) \to
\aut(\widehat{p_{L_2}(S_2)})$$ is surjective.
 But under the isomorphism $\aut(\hat{S}_2) \cong \aut(\hat{S}_1, L_1)$, the subalgebra $\aut(\hat{S}_2, L_2)$ is sent into $$\aut(\hat{S}_1, L_1 \subset L_3):= \{ \sigma \in \aut(\hat{S}_1, L_1), \sigma(L_3) \subset L_3\} \subset \aut(\hat{S}_1, L_3)$$ from which the surjectivity is clear.
Now for the condition (ii) in Proposition \ref{p.proj}, if $v \in
T_{\alpha}(\hat{S}_1) \cap L_3$ for a general $\alpha \in
\hat{S}_1$, then its image $\bar{v} \in V_2$ satisfies $\bar{v}
\in T_{\bar{\alpha}}(\hat{S}_2) \cap L_2$. Thus by condition (ii)
for $(S_2, V_2, L_2)$, we have $\bar{v} =0$, i.e., $v \in L_1$.
Then by condition (ii) for $(S_1, V_1, L_1)$, we get $v=0.$
\end{proof}

We have the following

\begin{corollary}\label{c.Symp}
In the notation of Proposition \ref{p.tanSymp}, let $Z:= p_L(S)
\subset \BP U :=\BP(V/L)$ be the VMRT of a symplectic Grassmannian
as explained in Section \ref{e.SymGrass}. For a subspace $L_2
\subset U$, let $L_3 \subset \Sym^2(W \oplus Q)$ be the inverse
under the projection $\Sym^2(W \oplus Q) \to U$. If $\BP L_2 \cap
\Sec(Z) = \emptyset$, then $$\aut(\widehat{p_{L_2}(Z)})^{(1)}
\cong \aut(\hat{S}, L_3,0)^{(1)}.$$ \end{corollary}

\begin{proof}
We just apply Proposition \ref{p.double} with $S_1= S, V_1= V,
L_1=L,$ together with Propositions \ref{p.aut}, \ref{p.tanSymp}
and \ref{p.autSymp}. \end{proof}

We can make this more explicit as follows.

\begin{proposition}\label{p.Symp}
For $\phi \in \Sym^2(W \oplus Q)$, denote by $\phi^{\sharp}\in
\Hom(W^*\oplus Q^*, W \oplus Q)$ the corresponding homomorphism.
 For $L_2 \subset U$, we denote by $L_3
 \subset \Sym^2(W \oplus Q)$ the subspace satisfying $L_3/(\Sym^2 Q) \cong L_2$ and by
$\Im(L_2)$ the linear span of $\{ \Im (\phi^{\sharp}) \subset
 W \oplus Q, \phi \in L_2\}.$ Define $\Im_W(L_2) : = p_Q( \Im(L_2)) \subset W$,
 where $p_Q: W \oplus Q \to W$ is the  projection to the first factor.
 Then \begin{enumerate} \item[(i)] there is a canonical vector space isomorphism
 $$
\aut(\hat{S}, L_3,0)^{(1)} \cong \Sym^2 (W/\Im_W(L_2))^*;
 $$
 \item[(ii)] for each $\psi \in \Sym^2 (W/\Im_W(L_2))^* \subset \Sym^2 W^*,$ denoting by $\psi^{\sharp}$ the corresponding element in $\Hom(W, W^*)$ and writing $W':= \Ker(\psi^{\sharp}),$ $L_2$ is contained in
     $$L'(\psi):= \{ \phi \in U| \Im(\phi^{\sharp}) \subset \Ker(\psi^{\sharp})\oplus Q\} \cong (W' \otimes Q) \oplus \Sym^2 W';$$
    \item[(iii)] if $\BP L_2 \cap \Sec(Z) = \emptyset$ and $\aut(\hat{S}, L_3, 0)^{(1)}$ contains an element of
 rank $r$ in $\Sym^2 (W/\Im_W(L_2))^*$ (i.e. the corresponding
 element in $\Hom(W, W^*)$ has rank $r$),
$$
\dim L_2 \leq
mk+\cfrac{k(k+1)}{2}-2m-2k+1-\cfrac{r(2m+2k-r-3)}{2}.
$$\end{enumerate} \end{proposition}

The following lemma is immediate from Lemma \ref{l.secant} and the
information on $\Sec(S)$ ($S$ as in Proposition \ref{p.tanSymp})
from the table in section \ref{e.IHSS}.

\begin{lemma}\label{l.Symp}
In the notation of Proposition \ref{p.Symp}, $\dim \Sec(Z) = 2m +
2k -2$ where $m= \dim Q, k = \dim W$. In particular, $\Sec(Z) =
\pit(U)$ if and only if $k=2$. \end{lemma}

\begin{proof}[Proof of Proposition \ref{p.Symp}]
From Proposition \ref{p.III} (i), we have
$$\aut(\hat{S}, L_3, 0)^{(1)} \cong \Sym^2((W \oplus Q)/\Im(L_3))^*.$$ From $L=\Sym^2 Q \subset L_3$
 and $L_3/L = L_2$, we have $\Im(L_3) = Q \oplus \Im_W(L_2)$, proving (i).

Now for $\psi \in \Sym^2 (W/\Im_W(L_2))^*$,   denote by
$\tilde{\psi}^\sharp \in \Hom( W \oplus Q, W^* \oplus Q^*)$  the
element induced by  $\psi^{\sharp} \in \Hom(W, W^*)$ via the
composition
$$ W\oplus Q \stackrel{p_Q}{\longrightarrow} W
\stackrel{\psi^{\sharp}}{\longrightarrow} W^* \hookrightarrow
W^*\oplus Q^*.$$  From Proposition \ref{p.III} (ii), we see that
$L_3$ is contained in
$$L(\psi) := \{ \phi \in \Sym^2(W \oplus Q)| \Im(\phi^{\sharp})
 \subset \Ker(\tilde{\psi}^{\sharp})\}. $$ Clearly, $\Sym^2 Q \subset L(\psi)$
 and the quotient $L(\psi)/(\Sym^2 Q)$ is naturally isomorphic to
 $L'(\psi)$, proving (ii).

For (iii), assume that $\BP L_2 \cap \Sec(Z) = \emptyset$ and
there is $\psi$ in (ii) of rank $r$. Then $\dim W' = k-r.$ Let
$V'= \Sym^2(W'\oplus Q)$ and $S' \subset \BP V'$ be the second
Veronese embedding of $\BP(W' \oplus Q)$. Then  we have
$\dim(\Sec(S')) = 2m +2(k-r)-2$. By the assumption, $\BP L_2
\subset \BP L'(\psi)$ is disjoint from $$ p_{L}(\Sec(S')) \subset
\Sec(p_{L}(S)) = \Sec(Z),$$ which implies that $k-r \geq 3$ by
Lemma \ref{l.Symp} and also $$\dim L_2 \leq
m(k-r)+\frac{(k-r)(k+1-r)}{2}-1-(2m+2(k-r)-2).$$
\end{proof}

Let us derive an important consequence of our study of the
prolongation of biregular projections of the examples in Section
\ref{s.examples},  Theorem \ref{t.mainprolong} below, which is a
key ingredient in the proof of Main Theorem.

\begin{definition}
A linear subspace $\emptyset \neq \BP L \subset \BP V \setminus
\Sec(S)$ is called {\em maximal} if $\Sec(p_L(S)) = \BP(V/L)$. In
this case,  $\dim L = \dim \BP V - \dim \Sec(S)$ from Lemma
\ref{l.secant}.
\end{definition}

\begin{theorem} \label{t.mainprolong}
Let $S \subset \BP V$ be one of the linearly normal varieties
listed in Main Theorem (A1)-(B5) with $\Sec(S) \neq \BP V$. Let
$\BP L \subset \BP V \setminus \Sec(S)$ be a linear space and
$p_L$ the projection along $\BP L$. If $\BP L$ contains a general
point of $\BP V$ or if $\BP L$ is maximal, then
$\aut(\widehat{p_L(S)})^{(1)} = 0$.
\end{theorem}

\begin{proof}
From $\Sec(S) \neq \BP V$, it suffices to check those covered by
Propositions \ref{p.I}, \ref{p.II}, \ref{p.III}, \ref{p.VI} and
\ref{p.Symp}. In fact, the examples in Section \ref{e.S_5} and
Section  \ref{e.Gr25} admit no biregular projections since
$\Sec(S) = \BP V$ by \cite{Za} (Chapter V, Corollary 1.13). There
is nothing to check for the case of Proposition \ref{p.VI}.

In Proposition \ref{p.I}, suppose $L$ contains a general element
$\phi$ of $V$. Then $\phi \in L$ is of maximal rank and
$\Im(L)=B$, proving $\aut(\widehat{p_L(S)})^{(1)} \cong
\aut(\hat{S}, L, 0)^{(1)} =0$ from Proposition \ref{p.I} (i). On
the other hand, if $L$ is maximal, $\dim L = ab -2a-2b +4$. Thus
from Proposition \ref{p.I} (iii), the rank of any element of
$\aut(\hat{S}, L, 0)^{(1)}$ must be zero, i.e., $\aut(\hat{S}, L,
0)^{(1)}=0.$

In Propositions \ref{p.II} (resp. \ref{p.III}), suppose $L$
contains a general element $\phi$ of $ V$. Then $\phi^{\sharp}$ is
of maximal rank and $\Im(L)=W$, proving
$\aut(\widehat{p_L(S)})^{(1)} \cong \aut(\hat{S}, L, 0)^{(1)} =0$
from Proposition \ref{p.II} (i) (resp. Proposition \ref{p.III}
(i)). On the other hand, if $L$ is maximal, $\dim L =
\frac{1}{2}(n^2-n-2) - 4n +11$ (resp. $\frac{1}{2}(n^2 + n -2) -
(2n-2)$). Thus from Proposition \ref{p.II} (iii), (resp.
Proposition \ref{p.III} (iii)), the rank of any element of
$\aut(\hat{S}, L, 0)^{(1)}$ must be zero, i.e., $\aut(\hat{S}, L,
0)^{(1)}=0.$

In Proposition \ref{p.Symp},   suppose $L_2$ contains a general
$\phi\in V_2$. Then $\phi^{\sharp}$ is of maximal rank and
$\Im_W(L_2)=W$, proving $\aut(\widehat{p_{L_2}(Z)})^{(1)} \cong
\aut(\hat{S}, L_3, 0)^{(1)} =0$ from Proposition \ref{p.Symp} (i).
On the other hand,  $\dim(\Sec(Z))  = 2m+2k-2$ by Lemma
\ref{l.secant},
 which implies
that if $L_2$ is maximal, $\dim L_2 = mk +
\frac{k(k+1)}{2}-1-(2m+2k-2)$. Thus from Proposition \ref{p.Symp}
(iii), the rank of any element of $\aut(\hat{S}, L_3, 0)^{(1)}$
must be zero, i.e., $\aut(\widehat{p_{L_2}(Z)})^{(1)}=0.$
 \end{proof}

\section{Cone structure and $G$-structure}\label{s.G-structure}

This section collects some general facts on cone structures and
G-structures. The main theme is to reveal the relationship between
the existence of an Euler vector field, the local flatness of the
cone structure and the prolongation of a linear Lie algebra.
\begin{definition} A {\em cone structure} on a complex manifold $M$ is a
closed analytic subvariety
 $\sC \subset \BP T(M)$  such that the projection
$\pi: \sC \to M$ is  proper, flat and surjective with connected
fibers. A cone structure induces an equivalence relation on $M$:
two points $x, y \in M$ are equivalent if the projective varieties
$\sC_x \subset \BP T_x(M)$ and $\sC_y \subset \BP T_y(M)$ are
projectively equivalent. The equivalence classes define a
holomorphic foliation (possibly with singularity) whose leaves are
maximally connected submanifolds of maximal dimension in $M$
consisting of equivalent points. These leaves will be called the
{\em isotrivial leaves} of the cone structure. The  dimension of
isotrivial leaves is to be denoted by $\delta(\sC)$. If
$\delta(\sC) = \dim M$, i.e., all general fibers of $\sC \to M$
are projectively equivalent, then we say that the cone structure
is {\em isotrivial}. In this case, if we denote by $Z$ the
projective variety $\sC_x \subset \BP T_x M$, then we call $\sC$ a
{\em $Z$-isotrivial cone structure}.
\end{definition}

For an isotrivial cone structure, we can associate to it another
geometric structure: the $G$-structure.

\begin{definition}
Given a complex manifold $M$, fix a vector space $V$ with $\dim V
= \dim M$. The frame bundle $\sF(M)$ has the fiber at $x \in M$,
$$\sF_x(M) := {\rm Isom}(V, T_x(M)).$$ For a closed connected subgroup $G \subset
{\rm GL}(V)$, a $G$-{\em structure} on $M$ is a $G$-subbundle $\sG
\subset \sF(M)$.
 If $G$ contains the scalar group $\C^\times \cdot {\rm
Id} \subset {\rm GL}(V)$, we say that the $G$-structure is {\em of
cone type}. An isotrivial cone structure $\sC \subset \BP T(M)$
induces a $G$-structure $\sG$ of cone type on an open subset $M_o$
of $M$ where each fiber $\sC_x \subset \BP T_x(M), x \in M_o$ is
projectively equivalent  to $\sC_o \subset \BP T_o(M)$ for a base
point $o \in M_o$. In fact, setting $V= T_o(M)$ and $G =
\Aut_0(\hat{\sC}_o),$ the fiber
 $\sG_x \subset \sF_x(M)$ is given by $$\{ \sigma \in {\rm Isom}(V, T_x(M)), \sigma(\hat{\sC}_o) = \hat{\sC}_x \}.$$
\end{definition}

\begin{definition}
A $G$-structure $\sG \subset \sF(M)$ on $M$ and a $G$-structure
$\sG' \subset \sF(M')$ on $M'$ are {\em equivalent} if there
exists a biholomorphic map $\varphi: M \to M'$ such that the
induced map $\varphi_*: \sF(M) \to \sF(M')$ sends $\sG$ to $\sG'$.
 \end{definition}

 \begin{definition} On the vector space $V$ as a complex manifold,
we have a canonical trivialization $\sF(V) = {\rm GL}(V) \times
V.$ For any subgroup $G \subset {\rm GL}(V)$, this induces the
{\em flat $G$-structure}
 on $V$ defined by $$\sG = G \times V \subset {\rm GL}(V) \times V = \sF(V).$$
A $G$-structure $\sG \subset \sF(M)$ is {\em locally flat} if its
restriction to some open subset is equivalent to the restriction
of the flat $G$-structure to some open subset of $V$. An
isotrivial cone structure is {\em locally flat} if its associated
$G$-structure is locally flat. \end{definition}

\begin{definition} Given a cone structure $\sC \subset \BP T(M)$
and a  point $x \in M$, a germ of holomorphic vector field $v$ at
$x$ is said to {\em preserve the cone structure} if the local
1-parameter family of biholomorphisms integrating $v$ lifts to
local biholomorphisms of $\BP T(M)$ preserving $\sC$. The flows of
such a vector field must be tangent to the isotrivial leaves of
$\sC$.  The set of all such germs form a Lie algebra, called the
{\em Lie algebra of infinitesimal automorphisms of the cone
structure} $\sC$ at $x$, to be denoted by $\aut(\sC, x)$.
\end{definition}

\begin{definition}
Let $\sC \subset \BP T(M)$ be a cone structure.
 For a non-negative integer $\ell$, let $\aut(\sC, x)_{\ell}
\subset \aut(\sC,x)$ be the subalgebra of vector fields vanishing
at $x$ to order $\geq \ell + 1$. This gives the structure of a
filtered Lie algebra on $\aut(\sC, x)$, i.e.,
$[\aut(\sC,x)_{\ell}, \aut(\sC,x)_k] \subset \aut(\sC, x)_{\ell +
k}$.
\end{definition}

The following result is Proposition  1.2.1 \cite{HM05}.

\begin{proposition}\label{p.HM05}
For each $k \geq 0$, regard the quotient space
$\aut(\sC,x)_k/\aut(\sC,x)_{k+1}$ as a subspace of
$\Hom(\Sym^{k+1} T_x(M), T_x(M))$ by taking the leading terms of
the Taylor expansion of the vector fields at $x$. Then
$$
\aut(\sC,x)_k/\aut(\sC,x)_{k+1} \subseteq \aut(\hat{\sC}_x)^{(k)}.
$$
If the cone structure $\sC$ is isotrivial such that the associated
$G$-structure is locally flat, then the equality in the previous
inclusion holds for all $k$.
\end{proposition}

The following is a well-known fact in Poincar\'e normal form
theory of ordinary differential equations (e.g. \cite{AI} Sections
3.3.2 and 4.1.2).

\begin{lemma}\label{l.Poincare}
A germ of holomorphic vector fields at $(\C^n, 0)$ of the form
$$\sum_{i=1}^n (z_i + h_i(z)) \frac{\partial}{\partial z_i}$$ with
$h_i \in {\bf m}^2$ where ${\bf m} \subset \sO_{\C^n, 0}$ is the
maximal ideal, can be expressed as $\sum_{i=1}^n w_i
\frac{\partial}{\partial w_i}$ in  a suitable holomorphic coordinate
system $w_i$. \end{lemma}

\begin{definition}\label{d.Euler}
A germ of vector fields of the form in Lemma \ref{l.Poincare} is
called an {\em Euler vector field}. \end{definition}

\begin{proposition}\label{p.dim} Given a cone structure $\sC
\subset \BP T(M)$ and a general point $x \in M$, denote by
$\hat{\sC_x} \subset  T_x(M)$  the affine cone over the fiber
$\sC_x$ at $x$ and let $\aut(\hat{\sC_x}) \subset {\rm End}(T_x(M))$
be the Lie algebra of infinitesimal automorphisms of the affine
cone. Assume that $\aut(\hat{\sC_x})^{(k+1)} =0$ for some $k\geq 0.$
Then
$$ \dim(\aut(\sC, x)) \leq \delta(\sC) + \dim \aut(\hat{\sC_x}) +
\aut(\hat{\sC_x})^{(1)} + \cdots + \aut(\hat{\sC_x})^{(k)}$$ and if
equality holds then there exists an Euler vector field in
$\aut(\sC,x)_0$.
\end{proposition}

\begin{proof} The codimension of $\aut(\sC,
x)_0$ in $\aut(\sC,x)$  is at most $\delta(\sC)$. That the dimension of
$\aut(\sC, x)_0$ is bounded by $$\dim \aut(\hat{\sC_x}) +
\aut(\hat{\sC_x})^{(1)} + \cdots + \aut(\hat{\sC_x})^{(k)}$$
follows from Proposition \ref{p.HM05}, which also shows that the equality
holds only if each element of $\aut(\hat{\sC_x}) \subset {\rm
End}(T_x(M))$ can be realized as the linear part of a vector field
in $\aut(\sC, x)_0$. Thus if the equality holds,
 there exists an Euler vector field in $\aut(\sC, x)$.
\end{proof}

The following is from \cite{Gu} (also see Section 1 of
\cite{HM97}).

\begin{theorem}\label{t.Guillemin}
Given a $G$-structure $\sG \subset \sF(M)$, we can define
vector-valued  functions $c^k$, $k = 0, 1, 2 \ldots$ on $\sG$ with
the following properties.
\begin{enumerate} \item[(1)] $c^k$ is a $H^{k,2}(\fg)$-valued
function, well-defined if $c^{k-1} \equiv 0$. Here $H^{k,2}(\fg)$
is the cohomology of the natural sequence $$\fg^{(k)} \otimes V^*
\to \fg^{(k-1)} \otimes \wedge^2 V^* \to \fg^{(k-2)} \otimes
\wedge^3 V^*$$ and by convention, $c^{-1} \equiv 0$, $\fg^{(-1)} =
V$ and $\fg^{(-2)} = 0.$ \item[(2)] Under the action of $G$ on
$\sG$, the function $c^k$ transforms like the $G$-module
$\fg^{(k-1)} \otimes \wedge^2 V^* \subset \Hom(\Sym^k V, V)
\otimes \wedge^2 V^*$. \item[(3)] If $c^k \equiv 0$ for all
non-negative integers $k$, then $\sG$ is locally flat. \item[(4)]
$c^k$ is an invariant of the $G$-structure, i.e., it is invariant
under an automorphism of the $G$-structure. \end{enumerate}
\end{theorem}

It has the following consequence.

\begin{proposition}\label{p.locflat}
Let $\sC \subset \BP T(M)$ be a cone structure.  Assume that for a
general point $x \in M$, there exists an Euler vector field in
$\aut(\sC,x)_0$. Then the cone structure is isotrivial and locally
flat.
\end{proposition}

\begin{proof}
In a neighborhood of a general point $x$, the isotrivial leaves of
$\sC$ form a regular foliation. Given any vector field $v \in
\aut(\sC,x)_0,$ the flows of $v$ must be tangent to the leaves of
the foliation. But by Lemma \ref{l.Poincare}, each flow of an
Euler vector field $v$  has limit $x.$ Thus $x$ is a singularity
of the foliation, unless there is only one leaf. This shows that
$\sC$ is isotrivial.

To prove the local flatness, by Theorem \ref{t.Guillemin}, it
suffices to show that the functions $c^k$ on the associated
$G$-structure of cone type are identically zero. By induction,
assume that $c^{k-1} \equiv 0$ and $c^k$ is well-defined. Pick a
general point $x \in M$. The subgroup $\C^\times \cdot {\rm Id}
\subset G$ acts on the fiber $\sG_x$ and under this action,  the
characteristic function $c^k$ is multiplied by $t^{-(k+1)} \in
\C^\times$ by Theorem \ref{t.Guillemin} (2). But by integrating
the Euler vector field in $\aut(\sC, x)_0$, we get a 1-parameter
family of automorphisms of the $G$-structure which preserve the
fiber $\sG_x$ and acts by $\C^\times \cdot {\rm Id}$-action on it.
Since this is an automorphism of the $G$-structure, the functions
$c^k$ cannot change under this action by Theorem \ref{t.Guillemin}
(4), a contradiction unless  $c^k$  vanishes on $\sG_x$. Since
this is true for a general $x$, we get $c^k \equiv 0$. Thus by
induction we have $c^k \equiv 0$ for all $k \geq 0$ and the local
flatness of $\sG$ from Theorem \ref{t.Guillemin} (3).
\end{proof}

\begin{corollary}\label{c.dim}
If the equality holds in Proposition \ref{p.dim}, then the cone
structure is locally flat. \end{corollary}

Let us now recall some results on the infinitesimal automorphisms
of a locally flat $G$-structure.
 The following can be seen from Section 2.1 of \cite{Ya}.

\begin{proposition}\label{p.flat}
Assume that $\fg^{(k+1)} = 0$ and that the $G$-structure $\sG$ on
a complex manifold $M$ is locally flat. Then for any point $x \in
M$, $\aut(\sG, x)$, the Lie algebra of germs of holomorphic vector
fields at $x$ preserving the $G$-structure, is isomorphic to the
graded Lie algebra $V \oplus \fg \oplus \fg^{(1)} \oplus \cdots
\oplus \fg^{(k)}$ for a vector space $V$ with $\dim V= \dim M$.
\end{proposition}

\begin{proposition}\label{p.flat1}
Assume that $\fg^{(1)} = 0.$  Then the identity component
of the
 automorphism group of the flat $G$-structure $\sG$ on $V$ is the subgroup $V \sd G$ of the affine group
 $V \sd {\rm GL}(V).$ Moreover for any point $x \in V$, $\aut(\sG, x)$ is isomorphic to the Lie algebra of global holomorphic vector fields
 preserving the $G$-structure, i.e., $\aut(\sG, x) = V \sd \fg.$
 Consequently, given two connected open subsets $W_1, W_2 \subset V$ with
 a biholomorphic map $\varphi: W_1 \to W_2$ such that $\varphi_*: \sF(W_1) \to \sF(W_2)$
 preserves the flat $G$-structure, there exists a global  affine transformation $\tilde{\varphi}: V \to V$
 preserving the flat $G$-structure, extending $\varphi$.  \end{proposition}

\begin{proof}
It is obvious that $V\sd G$ acts on $V$ preserving the flat
$G$-structure and $V \sd \fg \subset \aut(\sG,x)$. As
$\fg^{(1)}=0$, Proposition \ref{p.flat} implies that $\aut(\sG, x)
= V \sd \fg.$ For the next statement, note that $\varphi_*$
induces an isomorphism of $\aut(\sG, x_1)$ and $\aut(\sG, x_2)$
for $x_1 \in W_1$ (resp. $x_2 \in W_2$) which sends  $\aut(\sG,
x_1)_0$ to $\aut(\sG, x_2)_0.$ This induces an isomorphism of the
universal covering of  $V \sd G$ to itself sending the isotropy
subgroup at $x_1$ to the isotropy subgroup at $x_2$. This
isomorphism descends to the desired affine transformation
$\tilde{\varphi}.$
\end{proof}

Proposition \ref{p.flat1} enables us to introduce developing map:
\begin{proposition}\label{p.develop}
Assume that $G \subset {\rm GL}(V)$ satisfies  $\fg^{(1)}=0$. Let
$M$ be a simply connected complex manifold equipped with a locally
flat $G$-structure. Then there exists an unramified holomorphic
map $\delta : M \to V$, called {\em developing map} of the
$G$-structure, such that the induced map on the frame bundles
$\delta_*: \sF(M) \to \sF(V)$ sends the $G$-structure on $M$ to
the flat $G$-structure on $V$. \end{proposition}

\begin{proof}
Fix a point $x\in M$. From the definition of a locally flat
$G$-structure, we have a neighborhood $U$ of $x$ and an unramified
holomorphic map $\delta_U: U \to V$ such that $\delta_{U*}$ sends
the $G$-structure on $U$ to the flat $G$-structure on $V$. By
Proposition \ref{p.flat1} we can extend $\delta_U$ to $\delta: M
\to V$ by analytic continuation as follows. For a given point
$y\in M$, we choose a path $\gamma: [0,1] \to M$ joining $x=
\gamma(0)$ to $y = \gamma(1)$. Then we can find finitely many
points $x=x_0, x_1, \ldots, x_N=y$ on $\gamma([0,1])$ such that
each $x_i$ has a neighborhood $U_i$ with an unramified holomorphic
map $\delta_{U_i}: U_i \to V$ sending the $G$-structure to the
flat $G$-structure. By shrinking $U_i$ if necessary, we may assume
that for each $i, 1 \leq i \leq N$,  $U_i \cap U_{i+1}$ is
connected and
$$(U \cup U_1 \cup \cdots \cup U_i) \cap U_{i+2} = \emptyset.$$ By
Proposition \ref{p.flat1}, we can find an affine automorphism
$\eta_{0,1}$ of $V$ preserving the $G$-structure such that
$\eta_{0,1} \circ \delta_{U_1}$ agrees with $\delta_{U}$ on the
intersection $U \cap U_1$. Replacing $\delta_{U_1}$ by $\eta_{0,1}
\circ \delta_{U_1}$, we can extend $\delta_U$ to the open subset
$U \cup U_1$ to get
$$\delta_{U \cup U_1}: U \cup U_1 \to V$$ such that the induced
map $$\delta_{U \cup U_1*}: \sF(U \cup U_1) \to \sF(V)$$ sends the
$G$-structure to the flat $G$-structure. Now repeating the same
argument, replacing $U$ by $U\cup U_1$ and $U_1$ by $U_2$, we can
extend $\delta_U$ to the open subset $U \cup U_1 \cup U_2$.
Continuing this way, we can extend $\delta_U$ to a neighborhood of
the path $\gamma([0,1])$. This defines the value $\delta(y) \in
V$. Since $M$ is simply connected, the value of $\delta(y)$ does
not depends on the choice of $\gamma$ and we can define the
desired map $\delta: M \to V.$
\end{proof}

\section{Proof of Main Theorem modulo Theorem \ref{t.non1} and Theorem \ref{t.non2}}\label{s.class}

In this section, we prove Main Theorem modulo two technical
results, Theorem \ref{t.non1} and Theorem \ref{t.non2}, the proofs
of which will be postponed to Section \ref{s.non1} and Section
\ref{s.non2}, respectively.

Ionescu and Russo's classification  of conic-connected manifolds
in \cite{IR} will be essential in our proof. To recall their
result, it is convenient to introduce the following definition.

\begin{definition}\label{d.primitive}
A conic-connected manifold $X \subset \BP^N$ is said to be {\em
primitive}, if  $X$ is a Fano manifold with $\Pic(X)$ generated by
$\sO_X(1)$
%\item[(2)] the index of $X$
%satisfies $i(X) \geq (n+1)/2,$ where the index $i(X)$ is the
%integer satisfying $-K_X \cong \sO_X(i(X))$ in $\Pic(X)$; and
%\item[(3)]
and $X$ is covered by lines.
\end{definition}

\begin{theorem}[\cite{IR}, Theorem 2.2]\label{t.IR}
Let $X \subset \BP^N$ be a conic-connected manifold of dimension
$n$. Then either $X$ is primitive or it is projectively equivalent
to one of the following or their biregular projections:

\begin{enumerate} \item[(a1)] The second Veronese embedding of $\BP^n$.

\item[(a2)] The Segre embedding of $\BP^a \times \BP^{n-a}$ for $1
\leq a \leq n-1$. \item[(a3)] The VMRT of the symplectic
Grassmannian $\Gr_\omega(k, k+n+1)$ for $2 \leq k \leq n$.

\item[(a4)] A hyperplane section of the Segre embedding $\BP^a
\times \BP^{n+1-a}$ with $2 \leq a, n+1-a$.

\end{enumerate}
\end{theorem}

We know the prolongations of the varieties (a1), (a2) and (a3) in
Theorem \ref{t.IR} from Section \ref{s.prolong}.  The prolongation
of varieties (a4) in Theorem \ref{t.IR} turns out to be zero:

\begin{proposition} \label{p.hyperplaneSegre}
Let $a \geq b \geq 2$ be two integers.  Let $X = \pit^a \times
\pit^b \hookrightarrow \pit^{ab+a+b}$ be the Segre embedding and
$S= X \cap H$ a nonsingular hyperplane section, which is
conic-connected. Then for the non-degenerate embedding  $S \subset
H$, we have $\aut(\hat{S})^{(1)} = 0$.
\end{proposition}
\begin{proof}
The two projections $X \to \pit^a$ and $X \to \pit^b$ induce
 two fibrations: $\pit^a
\xleftarrow{\pi_1} S \xrightarrow{\pi_2} \pit^b$, with fibers
isomorphic to $\pit^{b-1}$ and $\pit^{a-1}$ respectively (cf. the
proof of Theorem 2.2 in [IR], Case II, subcase (b) ). Let $F
\subset T(S)$ be the distribution spanned by the tangent spaces of
fibers of $\pi_1$ and $\pi_2$, then $F$ has rank $a+b-2$. Note
that the projectivization $\BP F \subset \BP T(S)$ is a cone
structure, which is invariant under $\Aut(S)$.  If
$\aut(\hat{S})^{(1)} \neq 0$, then $F$ is integrable by Theorem
\ref{t.HM} and Proposition \ref{p.locflat}. Thus $F$ is a
foliation with
 leaves of codimension 1. For a general point $x \in S$, let
$R_x$ be the set of points on $S$ which can be connected by a
chain of lines contained in the fiber of $\pi_1$ or $\pi_2$. Then
$R_x$ must agree with the  leaf of $F$ through $x$, and is a
divisor on $S$.    Let $l_i$ be a line contained in the fiber of
$\pi_i$ such that $l_1$ meets $l_2$ at $x$. Then we have $0 = R_x
\cdot l_1 = R_x \cdot  l_2.$ By Lefschetz hyperplane theorem,
$H_2(S, \C) \cong \C^2$ is generated by the classes of $l_1$ and
$l_2$. Thus $R_x$ is a numerically trivial effective divisor, a
contradiction to $H^2(S, \C) \cong H^2(X, \C)$.
\end{proof}

In the setting of Main Theorem,  Theorem \ref{t.IR} has the
following consequence.

\begin{proposition}\label{p.Pic2}
Let $S \subset \BP V$ be a nonsingular non-degenerate variety of
dimension $n$ such that $\aut(\hat{S})^{(1)} \neq 0$. Then either
$S \subset \BP V$ is a  primitive conic-connected manifold or it
is projectively equivalent to one of  the following or their
biregular projections:

\begin{enumerate} \item[(a1)] The second Veronese embedding of $\BP^n$.

\item[(a2)] The Segre embedding of $\BP^a \times \BP^{n-a}$ for $1
\leq a \leq n-1$. \item[(a3)] The VMRT of the symplectic
Grassmannian $\Gr_\omega(k, k+n+1)$ for $2 \leq k \leq n$.
    \end{enumerate}
\end{proposition}

\begin{proof}
Note that by Proposition \ref{p.hyperplaneSegre} and Proposition
\ref{p.aut},  biregular projections of varieties in Proposition
\ref{p.hyperplaneSegre} have no prolongation. By Theorem
\ref{t.HM} (i), $S$ is conic-connected. Applying Theorem
\ref{t.IR}, we get Proposition \ref{p.Pic2}. \end{proof}

The following variation of Proposition \ref{p.Pic2} will be
useful.

\begin{proposition}\label{p.PicSec}
Let $S \subset \BP V$ be a nonsingular non-degenerate variety such
that $\aut(\hat{S})^{(1)} \neq 0$ and $\Sec(S) = \BP V$. Then
either $S \subset \BP V$ is a  primitive conic-connected manifold
or it is projectively equivalent to one of the following:
\begin{enumerate}
\item[(i)] The second Veronese embedding  $v_2(\BP^1) \subset
\BP^2$, i.e. a plane conic. \item[(ii)] The Segre embedding of
$\BP^1 \times \BP^k$ with $k \geq 1$. \item[(iii)] The VMRT of a
symplectic Grassmannian $\Gr_\omega(2, V)$ with $\dim V \geq 5$.
\end{enumerate}
\end{proposition}

\begin{proof}
From the table in Section \ref{e.IHSS},  Lemma \ref{l.secant}, and
Lemma \ref{l.Symp},  a variety in (a1)-(a3) of Proposition
\ref{p.Pic2}, satisfying $\Sec(S) = \BP V$ belongs to the above
list. On the other hand, by Theorem \ref{t.mainprolong}, if $S
\subset \BP V$ is a biregular projection of (a1)-(a3) in
Proposition \ref{p.Pic2} with $\aut(\hat{S})^{(1)} \neq 0$, then
it cannot satisfy $\Sec(S) = \BP V$.
\end{proof}

From Proposition \ref{p.Pic2}, the difficulty in proving Main
Theorem lies in the study of primitive Fano manifolds which has
nonzero prolongation. We will study at first the VMRT of such
varieties.

%An essential tool in the study of Fano manifolds of Picard number
%1 is the deformation theory of rational curves. We will need the
%following well-known result (cf. Section 1 in \cite{Hw01}).

%\begin{proposition}\label{p.standard} Recall that a rational curve $C$ in a
%projective manifold $X$ is standard if the normalization $\nu:
%\BP^1 \to C \subset X$ is an immersion and its normal bundle
%$\nu^*T(X)/T(\BP^1)$ is of the form $\sO(1)^p \oplus \sO^q$ for
%some non-negative integers $p$ and $q$.
%\begin{enumerate} \item[(1)] A standard rational curve does not
%have a non-trivial deformation fixing a point and the tangent
%direction at that point, i.e., if there exists a family $\{ C_t,
%|t|<1\}$ of standard rational curves on a manifold with $x \in
%\cap_{|t|<1} C_t$, then $\cup_{|t|<1} \BP T_x(C_t)$ must have
%positive dimension. \item[(2)] For an irreducible component $\sK$
%of ${\rm RatCurves}^n(X)$, if  a general member of $\sK$ is  not
%standard, it admits a non-trivial deformation fixing two points
%and the deformation degenerates to a connected union of rational
%curves of low degree. In particular, if $\sK$ is a minimal
%rational component, a general member of $\sK$ is standard.
%\end{enumerate} \end{proposition}

\begin{proposition}\label{p.irredu}
 Let $X$ be a Fano manifold of Picard number 1
such that for a general point $x \in X,$ there exists
 a holomorphic vector
field $v_x$ on $X$ which is an Euler vector field at $x$ (in the sense of Definition \ref{d.Euler}) and
generates a $\C^{\times}$-action on $X$. Then for any choice of a
minimal rational component, the associated VMRT at a general point
is irreducible.
\end{proposition}
\begin{proof}   Since $v_x$ is Euler at $x$, it acts on $\BP
T_x(X)$ trivially. The $\C^{\times}$-action generated by $v_x$
sends minimal rational curves through $x$ to minimal rational
curves through $x$ fixing their tangent directions. On the other
hand, the normal bundle (pull-back to the normalization) of a
general minimal rational curve $C$ is of the form $\sO(1)^p \oplus
\sO^q$ for some non-negative integers $p$ and $q$ (cf. Section 1
in \cite{Hw01}). This implies that $C$ does not have a non-trivial
deformation fixing a point and the tangent direction at that
point. Thus the $\C^{\times}$-action must send each minimal
rational curve through $x$ to itself, inducing a non-trivial
$\C^{\times}$-action on each minimal rational curve through $x$.
Denote by $\sN \to \sC$ the normalization of the total variety of
minimal rational tangents $\sC$ and $\sN \to \sM \to X$ the Stein
factorization of $\sN \to X$, where $f: \sM \to X$ is a finite
morphism. As $\sN$ is irreducible, to prove Proposition
\ref{p.irredu}, it suffices to show that $f$ is birational.
Suppose not and let $D \subset X$ be an irreducible component of
the branch locus. Any $\C^{\times}$-action on $X$ lifts to a
$\C^{\times}$-action on $\sM$ because it induces an action on the
space of minimal rational curves.  Let $C$ be a general minimal
rational curve through $x$. We can assume that $C$ intersects $D$
transversally. By Lemma 4.2 of \cite{HM01}, there exists a
component $C'$ of $f^{-1}(C)$ which is not birational to $C$.
After normalizing $C$ and $C'$, the morphism $f|_{C'}: C'\to C$
has non-empty branch points at least at two points $z_1, z_2 \in
C$. By the generality of $x$, $x \neq z_1, z_2$.  But $v_x$
generates a $\C^{\times}$-action on $C$ fixing $D \cap C$ and $x$.
Since $z_1, z_2 \in D \cap C$, we have a non-trivial
$\C^{\times}$-action on $\BP^1$ with at least three fixed points,
a contradiction.
\end{proof}

\begin{proposition} \label{p.VMRTCC}
Let $X \subset \BP^N$ be a primitive conic-connected manifold. Fix
a minimal rational component consisting of lines covering $X$.
Assume that $\aut(\hat{X})^{(1)} \neq 0$. Then for a general point
$x \in X$, the VMRT $\sC_x$ is an irreducible nonsingular and non-degenerate
projective variety satisfying $\Sec(\sC_x) = \BP T_x(X)$.
\end{proposition}

\begin{proof} We know that $\sC_x$ is  irreducible from Theorem \ref{t.HM} (iv) and Proposition \ref{p.irredu}.
It is non-singular from Proposition \ref{p.nonsingular}. It
remains to prove that $\Sec(\sC_x) = \BP T_x(X)$, which implies the
non-degeneracy. By the proof of Theorem 2.2 \cite{IR} (p. 155), if
all conics joining two general points are irreducible, then $X$ is isomorphic to the second
Veronese embedding, a contradiction to the assumption that $X$ is
covered by lines. Thus two general points of $X$ can be joined by
a connected union of two lines. Then by Theorem 3.14 of \cite{HK},
$\Sec(\sC_x) = \BP T_x(X)$. \end{proof}

%it suffices to show that two general points of $X$ can be joined
%by a connected union of two lines. Note that when $n$ is the
%dimension of $X$, a conic $C \subset X$ satisfies $C \cdot
%K^{-1}_X \geq n+1$ from Definition \ref{d.primitive}.

%Suppose a general conic on $X$ is not standard. Then it has
%deformation fixing two points degenerating into a connected union
%of two lines by Proposition \ref{p.standard} (2). Thus any two
%general points of $X$ can be joined by a connected union of lines.

%Thus we can assume that a general conic $C$ on $X$ is standard.
%From $C \cdot K^{-1}_X \geq n+1,$ its normal bundle is
%$\sO(1)^{n-1}$. For a general point $x \in X$, the
%$\C^{\times}$-action in Theorem \ref{t.HM} (iv) must be tangent to
%the standard conics through $x$, because otherwise, the
%$\C^{\times}$ deforms the conics with their tangents at $x$ fixed,
%a contradiction to the form of the normal bundle. Then the
%closures of the $\C^{\times}$-orbits are conics through $x$ such
%that a general member of these conics has no points in common
%other than $x$. Now on each of this conic the $\C^{\times}$-vector
%field vanishes at a point outside $x$. Thus this vector field must
%vanish on a hypersurface in $X$ and the orbits have degree 2 with
%respect to $\sO_X(1)$. This implies that $X$ is biregular to
%projective space by \cite{MS}, and $\sO_X(1)$ is isomorphic to
%$\sO(2)$ of projective space, a contradiction to the fact that $X$
%is covered by lines.

The following theorem called Cartan-Fubini type extension theorem,
was proved in \cite{HM01}.

\begin{theorem}\label{t.Cartan-Fubini}
Let $X, X'$ be two Fano manifolds of Picard number 1 and let $\sC,
\sC'$ be the VMRT's associated to some minimal rational
components. Assume that $\sC_x$ is irreducible and nonsingular for
a general point $x \in X$. Given any connected analytic open
subsets $U \subset X, U' \subset X'$ with a biholomorphic map
$\phi: U \to U'$ such that the differential $\phi_*: \BP T_x(X)
\to \BP T_{\phi(x)} (X')$ sends $\sC_x$ to $\sC'_{x'}$
isomorphically for all $x \in U$, then $\phi$ extends to a
biholomorphic map $\Phi: X \to X'$.
\end{theorem}

An immediate consequence of Theorem \ref{t.Cartan-Fubini} is the
following  which allows us to reconstruct some Fano manifolds of
Picard number 1 from its VMRT.

\begin{corollary}\label{c.reconstruct}
Let $X, X'$ be two Fano manifolds of Picard number 1 with
$Z$-isotrival VMRT for an irreducible nonsingular  projective
variety $Z \subset \BP V$. Assume that the VMRT structures are
both locally flat, then $X$ is biregular to $X'$.
\end{corollary}

Now we continue the study of the VMRT of varieties with
prolongation.
\begin{proposition} \label{p.equiv}
Let $S \subsetneq \BP V$ be an $n$-dimensional  non-degenerate primitive conic-connected manifold with
 $\aut(\hat{S})^{(1)} \neq 0$. Then
(i) the cone structure on a Zariski open subset of $S$ defined by VMRT's  is locally flat, and
 (ii)
$S$ is an equivariant compactification of the affine space $\C^n$.
\end{proposition}
\begin{proof}
(i) follows from Theorem \ref{t.HMk=2}.  To show (ii), recall that
the VMRT of $S$ at a general point $x$ is irreducible, nonsingular
and non-degenerate by Proposition \ref{p.VMRTCC}. Thus by Theorem
\ref{t.HM}, the corresponding $\fg$ of this cone structure
satisfies $\fg^{(2)} =0$. By Proposition \ref{p.flat}, $\aut(\sC,
x)$ contains the abelian subalgebra $\C^n$. Moreover, the induced
local action of $\C^n$ has an open orbit in a neighborhood of $x$.
By Theorem \ref{t.Cartan-Fubini}, $\aut(\sC, x) \cong \aut(S)$,
which implies that $\C^n$ acts on $S$ with an open orbit. By Lemma
\ref{l.algebraic} below, this is an algebraic action of $\C^n$ and
its isotropy subgroup is an algebraic subgroup. But the isotropy
must be discrete because the orbit is open. Thus the isotropy is
trivial and $S$ is an equivariant compactification of $\C^n$.
\end{proof}

 \begin{lemma}\label{l.algebraic} Let $W$ be a vector space and let $\tau: B \to W \sd {\rm GL}(W)$ be an injective
 complex-analytic
 homomorphism of a complex algebraic group $B$ into the affine group  whose image $\tau(B)$ contains $W$. Then $\tau^{-1}(W) \subset B$
 is an algebraic subgroup of $B$. \end{lemma}

 \begin{proof}
 In the affine group $W \sd {\rm GL}(W)$, $W$ is a maximal connected abelian subgroup, i.e., any
 connected abelian subgroup containing $W$ is $W$ itself. If $\tau^{-1}(W)$ is not algebraic, let $W' \subset B$ be the Zariski closure of $\tau^{-1}(W)$. Then $W'$ is a connected abelian subgroup of $B$ strictly containing $W$. It follows that $\tau(W')$ is an abelian subgroup of $W \sd
  {\rm GL}(W)$ strictly containing $W$, a contradiction to the maximality of $W$. \end{proof}

The following theorem enables us to use induction to study prolongation,
and is a crucial step in the proof of Main Theorem.

\begin{theorem}\label{t.induction}
Let $S \subsetneq \BP V$ be a primitive conic-connected manifold
with
 $\aut(\hat{S})^{(1)} \neq 0.$ Then the VMRT $\sC_x \subset \BP
T_x(S)$ at a general point $x \in S$ is an irreducible nonsingular
non-degenerate variety satisfying $\aut(\hat{\sC}_x)^{(1)} \neq 0$
and $\Sec(\sC_x) = \BP T_x(S)$.
\end{theorem}
\begin{proof}
All follow from Proposition \ref{p.VMRTCC}, except $\aut(\hat{\sC}_x)^{(1)} \neq 0$.

 By Proposition
\ref{p.aut=1}, we may assume $\dim \aut(\hat{S})^{(1)} \geq 2$.
 For each $A \in \aut(\hat{S})^{(1)}$, we have an associated element
 $\lambda_A \in H^0(S, \0(1))$ in the sense of Theorem \ref{t.HM}. By
 Proposition \ref{p.equiv} (ii), we have an algebraic action of
 $\C^n$ on $S$ with an open orbit. The complement of the open
 orbit is an irreducible hypersurface from $b_2(S) =1$ (e.g.
 Proposition 1.2 (c) in \cite{PS}).
 Suppose that the hyperplane $(\lambda_A =0)$ intersects the
open $\C^n$-orbit. As $S$ is covered by lines and $S \subsetneq \BP V$, $S$ is not
biregular to a projective space and we can apply  Proposition \ref{p.secondorder} and
Proposition \ref{p.HM05}, to get a point $x \in S \cap (\lambda_A
=0)$ with $\aut(\hat{\sC}_x)^{(1)} \neq 0$ such that $x$ is in the
open orbit of the $\C^n$-action. Thus $\aut(\hat{\sC}_x)^{(1)}
\neq 0$ holds for a general point $x$ of $S$.

So we may assume that the hyperplane $(\lambda_A=0)$ is disjoint
from the open orbit, i.e.,  for each $A \in \aut(\hat{S})^{(1)}$,
the complement of the hyperplane $(\lambda_A =0)$ is an open
$\C^n$-orbit.   Since $\dim \aut(\hat{S})^{(1)} \geq 2$,
 there are at least two distinct $\C^n$-orbits in $S$
whose complements are distinct as  hyperplane sections of $S$. By
$\Pic(S) \cong \Z$, this implies that ${\rm Aut}(S)$ has an open
orbit $M$ whose complement has codimension $\geq 2$. Since $S$ is
simply connected, so is $M$.  Moreover, the VMRT defines a locally
flat cone structure on $M$ by Proposition \ref{p.equiv}. Suppose that $\aut(\hat{\sC}_x)^{(1)}
=0.$ By Proposition \ref{p.develop}, we get a nonconstant holomorphic
mao from $M$ to $\C^n$ and $M$ has non-constant holomorphic
functions, a contradiction.
\end{proof}

The following two theorems will be proved in the next two sections:

\begin{theorem}\label{t.non1}
Let $X$ be a Fano manifold with ${\rm Pic}(X) = \Z \langle
\sO_X(1) \rangle$. Assume that $X$ has minimal rational curves of
degree 1 with respect to $\sO_X(1)$ whose VMRT at a general point
is isomorphic to the VMRT of a symplectic Grassmannian
$\Gr_\omega(2, m+4)$ with $m \geq 2$. Then this cone structure is
not locally flat.
\end{theorem}

\begin{theorem}\label{t.non2}
Let $X$ be a 15-dimensional Fano manifold with ${\rm Pic}(X) = \Z
\langle \sO_X(1) \rangle$. Assume that $X$ has minimal rational
curves of degree 1 with respect to $\sO_X(1)$ whose VMRT at a
general point is isomorphic to a hyperplane section of the
10-dimensional spinor variety. Then this cone structure is not
locally flat.
\end{theorem}

Conjecturally, the Fano manifold in Theorem \ref{t.non1} (resp.
Theorem \ref{t.non2}) is isomorphic to $\Gr_\omega(2, m+4)$ (resp.
a general hyperplane section of the Cayley plane $\BO \BP^2$),
whose VMRT structure is not locally flat.
 Theorem \ref{t.non1} and
Theorem \ref{t.non2} are in contrast with the following. Note that
$\Gr_\omega(2,5)$ is the case $m=1$ in the setting of Theorem
\ref{t.non1}.

\begin{proposition}\label{p.S5flat}
Let $\mathbb{S}_5 \subset \BP^{15}$ be the spinor embedding of the
10-dimensional spinor variety and let $S \subset \BP^{14}$ be a
general hyperplane section. The cone structure on $S$ defined by
the VMRT of lines covering $S$ is locally flat at general points.
\end{proposition}

\begin{proof}
$S$ is of Picard number 1 and covered by lines.
 Since  $\aut(\hat{S})^{(1)} \neq 0$ by Proposition \ref{p.prolongZQ}, the
cone structure on $S$ is locally flat by
Proposition \ref{p.equiv} (i). \end{proof}

\begin{proposition}\label{p.SGr(2,5)}
Let $\Gr(2,5)  \subset \BP^{9}$ be the Pl\"ucker  embedding of the
Grassmannian and let $S:= \Gr_\omega(2,5) \subset \BP^8$ be a
general hyperplane section of $\Gr(2,5).$ The cone structure on
$S$ defined by the VMRT of lines covering $S$ is locally flat at
general points. \end{proposition}

\begin{proof}
This is by the same argument as in the proof of Proposition \ref{p.S5flat}, replacing Proposition \ref{p.prolongZQ} by Proposition \ref{p.prolongZQ2}.
 \end{proof}

We are now ready to prove Main Theorem.

\begin{proof}[Proof of Main Theorem]
Suppose that  $S$ is not primitive, then it belongs to (a1)-(a3)
in
 Proposition \ref{p.Pic2} or their biregular projections.  These varieties  correspond to (A1)-(A3)
 and the first entries in (B1)-(B6) of  Main Theorem, together with the biregular projections of (A1)-(A3).
 Here, note that the first entries in (B1)-(B6) do not
 have biregular projections.

    Now suppose that $S$ is  primitive. By Theorem \ref{t.induction}, the VMRT at a
general point $x\in S$, $\sC_x \subset \BP T_x(S)$ is a
nonsingular non-degenerate variety with $\aut(\hat{\sC}_x)^{(1)}
\neq 0$ and $\Sec(\sC_x) = \BP T_x(S)$. Then we can apply
Proposition \ref{p.PicSec} to $\sC_x \subset \BP T_x(S)$ to
determine $\sC_x$, unless $\sC_x$ is again primitive.
 Repeating this,  we end up with a positive integer $\ell$ and a sequence of
projective varieties
$$S_0 \subset \BP V_0, \; S_1 \subset \BP
V_1, \ldots, S_i \subset \BP V_i, \ldots, S_{\ell} \subset \BP
V_{\ell}, $$ such that
\begin{enumerate} \item[(a)] $S_{\ell}:=S$ and $ V_{\ell} =V$,
\item[(b)] when $x_i$ is a general point of $S_{i}$ for $1 \leq i
\leq \ell$, $S_{i-1} \subset \BP V_{i-1}$ is isomorphic to
$\sC_{x_i} \subset \BP T_{x_i}(S_{i})$ and the cone structure
given by this  VMRT  on $S_{i}$ is locally flat; \item[(c)]
$\aut(S_i)^{(1)} \neq 0$ for each $0 \leq i \leq \ell$ and
$\Sec(S_i) = \BP V_i$ for each $0 \leq i \leq \ell-1$; \item[(d)]
$S_i$ is primitive for each $ 1 \leq i \leq \ell$ and $S_{0}$ is one of the
varieties (i)-(iii) in Proposition \ref{p.PicSec}.
\end{enumerate}

We claim that the sequence of varieties $S_0, \ldots, S_{\ell}$
must be biregular to one of the following sequences of varieties.

\begin{enumerate}
\item[(b1)] $\BP^1, \; \Q^3, \; \Q^5, \;\ldots,\; \Q^{2\ell -1},
\; \Q^{2 \ell +1}.$ \item[(b2)]  $\BP^1 \times \BP^1, \; \Q^4, \;
\ldots, \; \Q^{2\ell}, \; \Q^{2\ell + 2}$. \item[(b3)] $\ell =1$
with  $S_0 \cong \BP^1 \times \BP^k$ and $S_1 \cong \Gr(2, k+3)$
for some $k \geq 3$.

\item[(b4)]  $\ell= 1,2$ or $3$,  with $S_0 \cong \BP^1 \times
\BP^2, \; S_1 \cong \Gr(2,5), \; S_2 \cong  \mathbb{S}_5$ and $S_3
\cong \BO \BP^2$. \item[(b5)] $\ell=1$ or $ 2$, with  $S_0 \cong
(\BP^1 \times \BP^2) \cap H_0, \; S_1 \cong \Gr(2,5) \cap H_1
(\cong \Gr_\omega(2, 5))$ and $S_2 \cong \mathbb{S}_5 \cap H_2$,
where $H_0, H_1, H_2$ are general hyperplanes.
\end{enumerate}

Once the claim is proved, then
by the property (c) of $\{S_i\}$ and Theorem \ref{t.mainprolong},
the embedding $S_i \subset \BP V_i$ for $0\leq i \leq \ell-1,$ is determined by the biregular
type of $S_i$ and is linearly normal, while $S_{\ell} \subset \BP V_{\ell}$ is determined up to biregular projections. Thus the list in (bi), $1 \leq i \leq 5$ give rise to the projective varieties
in (Bi) and (C) of Main Theorem, completing the proof of
Main Theorem.

To prove the claim let us recall that $S_{0} \subset \BP V_{0}$ must
be one of the following from Proposition \ref{p.PicSec}.

\begin{enumerate}
\item[(i)] The second Veronese embedding of $\BP^1 \subset \BP^2$.

\item[(ii)] The Segre embedding of $\BP^1 \times \BP^k$  with $k
\geq 1$.

\item[(iii)] The VMRT of a symplectic Grassmannian $\Gr_\omega(2,
V)$ with $\dim V \geq 5$.

\end{enumerate}

In Case (i), by a successive application of Corollary
\ref{c.reconstruct} combined with the property (b) of the sequence
$\{S_i\}$,  we obtain that $S_i$ is isomorphic to an
odd-dimensional hyperquadric, getting (b1).

In Case (ii),  $S_{1}$ is biregular to $\Gr(2, k+3)$ by Corollary
\ref{c.reconstruct}. If $\ell =1$, we end up with the sequence
(b3).  If $\ell \geq 2$,
 by the property (c) of $S_{1}$ combined with Theorem \ref{t.KN} and Theorem
 \ref{t.mainprolong},
 $S_1 \subset \BP V_1$ must be  the Pl\"ucker embedding of $\Gr(2, k+3)$
 with $k=1$ or $2$. If $k=1$, the Pl\"ucker embedding of $\Gr(2, k+3)$
is  the natural embedding of the 4-dimensional hyperqadric $\Q^4
\subset \BP^5$. By a successive application of Corollary
\ref{c.reconstruct}, we get that $S_i$ is an even-dimensional
hyperquadric, yielding the sequence (b2).  Now assume $k=2$, then
by Corollary \ref{c.reconstruct}, we get that $S_{2}$ is biregular
to
 $\mathbb{S}_5$, i.e.  the 10-dimensional spinor variety.
If $\ell =2$, we stop here, ending up with the case of $\ell =2$
in the sequence (b4). On the other hand, if $\ell \geq 3$,
  the embedding $S_2 \subset \BP V_2$ must be the spinor embedding
  by Theorem \ref{t.KN} and Theorem \ref{t.mainprolong}. Then $S_3$ is
  biregular  to the Cayley plane $\BO \BP^2$ by Corollary \ref{c.reconstruct}, giving $\ell =3$ in (b4).
It remains to show that $\ell \leq 3$. If $\ell \geq 4$, then by
(c) the embedding $ S_3 \subset \BP V_3$ must be the projection
along a general point of the minimal embedding $\BO \BP^2 \subset
\BP^{26}$, which has no prolongation by Proposition \ref{p.VI}, a
contradiction.

In Case (iii), if $\dim V \geq 6$, then  Theorem \ref{t.non1}
contradicts the property (b) of $S_1$. Thus this case,
corresponding to (B6) in Main Theorem, does not give rise to a
sequence with $\ell \geq 1$. Now we consider the case $\dim V=5$.
We want to show that the sequence must be (b5). First, $S_0$ is
isomorphic to a hyperplane section of $\BP^1 \times \BP^2$ under
the Segre embedding, which is the VMRT of $\Gr_{\omega}(2,5)$ at a
general point. Then By Corollary \ref{c.reconstruct} and
Proposition \ref{p.SGr(2,5)}, this implies that $S_{1}$ is
biregular to $\Gr_\omega(2, 5)$ and we are done if $\ell=1$. If
$\ell \geq 2$, then from the property (c) of $S_1$, the embedding
$S_1 \subset \BP V_1$ must be the hyperplane section of $\Gr(2,5)$
under the Pl\"ucker embedding. By Corollary \ref{c.reconstruct}
and Proposition \ref{p.S5flat}, this implies that $S_{2}$ is
biregular to a hyperplane section of the 10-dimensional spinor
variety and we are done if $\ell=2$.
 So it remains to show that $\ell \leq 2$. Suppose  $\ell \geq 3$, then from the property (c) of $S_2$ and Theorem \ref{t.mainprolong}, the embedding $S_2 \subset \BP V_2$ must be the hyperplane section of the spinor embedding $\mathbb{S}_5 \subset \BP^{15}$ and the
cone structure on $S_3$ given by $S_2$ cannot be locally flat by Theorem \ref{t.non2}. This contradicts the property (b) of $S_3$,
completing the proof of the claim.

\end{proof}

\section{Proof of Theorem \ref{t.non1}}\label{s.non1}

This section is devoted to the proof of Theorem \ref{t.non1}.

To start with, let us recall some facts about Grassmannians. Let
$W$ be a complex vector space of dimension 2 and $Q$  a complex
vector space of dimension $m \geq 2$.
 Let ${\rm Gr}(2, W^* \oplus Q)$ be the Grassmannian of
2-dimensional subspaces in $W^* \oplus Q$. There exists a
canonical embedding $$W \otimes Q = {\rm Hom}(W^*, Q) \subset
\Gr(2, W^* \oplus Q)$$ by associating to an element of ${\rm
Hom}(W^*, Q)$ the plane in $W^* \oplus Q$ given by its graph. The
next  proposition is elementary.

\begin{proposition}\label{p.actionGrass} Consider a $\C^{\times}$-action with weight 0 on $W$ and
weight 1 on $Q$. This induces a $\C^{\times}$-action on ${\rm
Gr}(2, W^* \oplus Q)$ whose fixed point set consists of the
following three components:
\begin{itemize} \item[(i)] the isolated point $[W^*]$ corresponding to the
plane $W^* \subset (W^* \oplus Q)$; \item[(ii)] the subvariety
${\rm Gr}(2, Q) \subset {\rm Gr}(2, W^* \oplus Q)$ consisting of
planes of $W^* \oplus Q$ contained in $Q$; \item[(iii)] the
subvariety $\BP W^* \times \BP Q$ consisting of planes which can
be written as the direct sum of a line in $W^*$ and a line in $Q$.
\end{itemize} Moreover under this $\C^{\times}$-action, the orbit $\C^{\times} \cdot z$ of any point $$ z \in
\; {\rm Gr}(2, W^* \oplus Q) \setminus ((W \otimes Q) \cup (\BP
W^* \times \BP Q))$$  has a limit point in ${\rm Gr}(2,
Q)$.\end{proposition}

Next, we need to look at  the geometry of a certain Grassmannian
bundle on a Lagrangian Grassmannian.  Let $\Sigma$ be a symplectic
vector space of dimension $4$ and denote by ${\rm Sp}(\Sigma)$
(resp. $\fsp(\Sigma)$) the Lie group (resp. algebra) of symplectic
automorphisms (resp. endomorphisms) of $\Sigma$. Let ${\rm
Lag}(\Sigma)$ be the Lagrangian Grassmannian, i.e., the space of
Lagrangian subspaces in $\Sigma$. This is homogeneous under ${\rm
Sp}(\Sigma)$ and is biregular to the 3-dimensional hyperquadric
$\Q^3$. Let $\sW$ be the universal quotient bundle on ${\rm
Lag}(\Sigma),$ i.e., the rank $2$ vector bundle satisfying $\Sym^2
\sW = T({\rm Lag}(\Sigma)).$ Its dual bundle is the tautological
bundle $\sW^* \subset \Sigma \times \Lag(\Sigma)$ whose fiber over
the point  $[W^*] \in \Lag(\Sigma)$ corresponding to a Lagrangian
subspace $W^* \subset \Sigma$ is $W^*$ itself. Fix a vector space
$Q$ of dimension $m \geq 2$ and denote by $\sQ$ the trivial vector
bundle on ${\rm Lag}(\Sigma)$ with a fiber $Q$.

\begin{proposition}\label{p.autGrass} Let ${\rm Gr}(2, \sW^* \oplus \sQ)$ be the
Grassmannian bundle  of 2-planes in the vector bundle $\sW^*
\oplus \sQ$. Then the Lie algebra $\fg$ of the automorphism group
of the projective variety ${\rm Gr}(2, \sW^* \oplus \sQ)$ is
isomorphic to $(\Sigma^* \otimes Q) \sd (\fsp(\Sigma) \oplus
\fgl(Q)).$ The vector bundle $\sW \otimes \sQ$ has a natural
embedding into ${\rm Gr}(2, \sW^* \oplus \sQ)$ whose complement
$$D := {\rm Gr}(2, \sW^* \oplus \sQ) \setminus (\sW \otimes \sQ)$$
is a hypersurface consisting of 2-planes in $\sW^* \oplus \sQ$
which have positive-dimensional intersection with $\sQ$.
\end{proposition}

\begin{proof}
The group $\Sigma^* \otimes Q=\Hom(\Sigma, Q)$ acts on the vector
space $\Sigma \oplus Q$ by the following rule: $f \cdot (x, y) =
(x, y+f(x))$ for any $x \in \Sigma, y\in Q$ and $f \in
\Hom(\Sigma, Q)$. This action preserves $\sW^* \oplus \sQ \subset
(\Sigma \oplus Q) \times \Lag(\Sigma)$, inducing an action of
$\Sigma^* \otimes Q$ on $\Gr(2, \sW^* \oplus \sQ)$. From this, we
can see there is a natural inclusion
$$(\Sigma^* \otimes Q) \sd (\fsp(\Sigma) \oplus
\fgl(Q)) \subset \fg.$$ To show that this is an isomorphism, it
suffices to compare their dimensions. Let
$$\psi: {\rm Gr}(2, \sW^* \oplus \sQ) \to {\rm Lag}(\Sigma)$$ be
the natural projection. We have an exact sequence
\begin{equation} \label{eq.sym}
0 \to T^{\psi} \to T({\rm Gr}(2, \sW^* \oplus \sQ)) \to \psi^*
T({\rm Lag}(\Sigma)) \to 0,
\end{equation} where $T^{\psi}$ denotes the relative tangent bundle. We have
$\psi_* T^{\psi} = {\rm End}^0(\sW^* \oplus \sQ)$, the bundle of
traceless endomorphisms, and $R^i\psi_* T^{\psi} =0$ for $i \geq 1$.  Write
$$\psi_* T^{\psi} =  \mathcal{F} \oplus (\sW^* \otimes \sQ^*) \oplus (\sW \otimes
\sQ),$$ where $\mathcal{F}$ is given by the exact sequence
\begin{equation}\label{eq.1}
0 \to \mathcal{O} \to {\rm End}(\sW^*) \oplus {\rm End}(\sQ)  \to
\mathcal{F} \to 0.
\end{equation}
Here the map $\mathcal{O} \to {\rm End}(\sW^*) \oplus {\rm
End}(\sQ) $ is given by $s \mapsto s {\rm Id}_{\sW^*} \oplus s
{\rm Id}_{\sQ}$.
 It is well-known that $$H^0({\rm Lag}(\Sigma), \sW) =
\Sigma^*, \; H^0({\rm Lag}(\Sigma), \sW^*) =0, \; H^0({\rm
Lag}(\Sigma), {\rm End}(\sW^*)) =\C, $$
$$
H^1({\rm Lag}(\Sigma), \sW) =H^1({\rm Lag}(\Sigma), \sW^*)
=H^1({\rm Lag}(\Sigma), {\rm End}(\sW^*)) =0,
$$ $$H^1({\rm Lag}(\Sigma), \sO) = H^2({\rm Lag}(\Sigma), \sO) = 0.$$
 Thus by the long-exact sequence associated to \eqref{eq.1},
 we have $$ H^0({\rm Lag}(\Sigma), \sF) = \fgl(Q), \; H^1({\rm Lag}(\Sigma), \sF) =0$$
  and consequently, $$ H^1({\rm
Gr}(2, \sW^* \oplus \sQ), T^\psi)=0$$  $$H^0({\rm Gr}(2, \sW^*
\oplus \sQ), T^{\psi}) = H^0({\rm Lag}(\Sigma), \psi_* T^{\psi}) =
(\Sigma^* \otimes Q ) \oplus \fgl(Q).  $$ Since $H^0({\rm
Lag}(\Sigma), T({\rm Lag}(\Sigma)) = \fsp(\Sigma)$, the long-exact
sequence associated to \eqref{eq.sym} shows that $$ \dim \fg= \dim
H^0(\Gr(2, \sW^* \oplus \sQ), T( \Gr(2, \sW^* \oplus \sQ)) = \dim
((\Sigma^* \otimes Q) \sd (\fsp(\Sigma) \oplus \fgl(Q))).$$

Now the vector bundle $\sW \otimes \sQ= {\rm Hom}(\sW^*, \sQ)$ can
be regarded as a subset of ${\rm Gr}(2, \sW^* \oplus \sQ)$ by
associating to a homomorphism to its graph. The statement on the
complement $D$ is immediate.
\end{proof}

\begin{proposition}\label{p.modelGrass}
Let $G$ be the simply connected group with Lie algebra $\fg$ of
Proposition \ref{p.autGrass}. The open subset $\sW \otimes \sQ
\subset {\rm Gr}(2, \sW^* \oplus \sQ)$ is homogeneous under the
action of $G$ and
 has a natural
isotrivial cone structure $\sC$ invariant under the $G$-action
such that each fiber $\sC_x \subset \BP T_x(\sW \otimes \sQ)$ is
isomorphic to $Z \subset \BP((W \otimes Q) \oplus \Sym^2 W)$, the
VMRT of the symplectic Grassmanian $\Gr_\omega(2, m+4)$ in the
notation of Section \ref{e.SymGrass}. This cone structure $\sC$ is
locally flat and $\aut(\sC, x) \cong \fg$ for each point $x \in
\sW \otimes \sQ$.\end{proposition}

\begin{proof}
Note that the hypersurface $D$ in Proposition \ref{p.autGrass} is
invariant under the action of $G$, hence $\sW \otimes \sQ$ is also
$G$-invariant. The base ${\rm Lag}(\Sigma)$ is homogeneous under
the action of ${\rm Sp}(\Sigma)$.
 Let $W^* \subset \Sigma$ be a Lagrangian subspace with quotient
$W = \Sigma/W^*.$ The subgroup $\Hom(\Sigma, Q) \subset G$ acts on
the fiber $W \otimes Q = \Hom(W^*, Q)$ of $\sW \otimes \sQ$ over
$[W^*] \in \Lag(\Sigma)$ by translation via the restriction to
$W^* \subset \Sigma$ of the action of $\Hom(\Sigma, Q)$ on $\Sigma \oplus Q$ described at the beginning of the proof of Proposition
\ref{p.autGrass}. This action is transitive on the fiber with
isotropy subgroup $$\{ \kappa \in \Hom(\Sigma, Q) \; | \;
\kappa(W^*) =0 \} \cong \Hom(\Sigma/W^*, Q).$$
This shows that $\sW \otimes \sQ$ is $G$-homogeneous.

 Regard ${\rm
Lag}(\Sigma)$ as a submanifold of $\sW \otimes \sQ \subset {\rm
Gr}(\sW^* \oplus \sQ)$ via the zero section of the vector bundle.
The Lie algebra of the isotropy subgroup in ${\rm Sp}(\Sigma)$ of
$[W^*] \in {\rm Lag}(\Sigma)$ is a parabolic subalgebra
$\fp_{[W^*]} \subset \fsp(\Sigma)$. At the point $[W^*] \in {\rm
Lag}(\Sigma) \subset \sW \otimes \sQ,$ the isotropy subgroup $G_{[W^*]}$ has Lie algebra
$$\fg_{[W^*]}:= (W^* \otimes Q) \sd (\fp_{[W^*]} \oplus \fgl(Q))
\; \subset \; (\Sigma^* \otimes Q) \sd (\fsp(\Sigma) \oplus
\fgl(Q)) = \fg.$$

The tangent space at the point $[W^*] \in {\rm Lag}(\Sigma)
\subset \sW \otimes \sQ$
$$T_{[W^*]}(\sW \otimes \sQ) = (W \otimes Q) \oplus T_{[W^*]}({\rm
Lag}(\Sigma)) = (W \otimes Q) \oplus \Sym^2 W$$ contains the
affine cone $\hat{Z}$ in a natural way.
 The isotropy representation of $\fg_{[W^*]}$ on $T_{[W^*]}(\sW \otimes \sQ) =(W \otimes Q) \oplus \Sym^2 W$  satisfies
\begin{itemize}
\item[(1)]  $(W^* \otimes Q)$-component of $\fg_{[W^*]}$ acts
trivially; \item[(2)]  $\fp_{[W^*]}$ acts naturally as $\fgl(W)$
on $W \otimes Q$ and  on $\Sym^2 W$; \item[(3)] $\fgl(Q)$ acts
naturally on $W \otimes Q$ and trivially on $\Sym^2 W$.
\end{itemize} Thus $\hat{Z}$ is
preserved under the isotropy representation of the isotropy
subgroup $G_{[W^*]}$ and the $G$-action defines a natural
isotrivial cone structure $\sC$ on the open set $\sW \otimes \sQ$
whose fiber is isomorphic to $Z$.
 As $\aut(\hat{Z})^{(2)} = 0$ by Theorem \ref{t.HMk=2}, we
have the following inequalities from Proposition \ref{p.dim} where
$x= [W^*]:$
$$
\dim \fg \leq \dim \aut(\sC, x) \leq \dim (\aut(\hat{Z})^{(1)}
\oplus \aut(\hat{Z}) \oplus (W \otimes Q) \oplus \Sym^2 W).
$$
From Propositions \ref{p.autZ},  $$ \dim
(\aut(\hat{Z})^{(1)} \oplus \aut(\hat{Z}) \oplus (W \otimes Q)
\oplus \Sym^2 W) = m^2+4m+10 = \dim \fg$$ implying $\fg \cong \aut(\sC, x).$
 Now
Corollary \ref{c.dim} gives the locally flatness of the cone
structure $\sC$.
\end{proof}

Now to prove Theorem \ref{t.non1}, we will make the following
assumption and derive a contradiction.

({\em Assumption}) Let $X$ be a Fano manifold with ${\rm Pic}(X) =
\Z \cdot L$ for an ample line bundle $L$. Assume that $X$ has
minimal rational curves of degree 1 with respect to $L$ whose VMRT
at a general point is isomorphic to $Z \subset \BP  ((W \otimes Q)
\oplus \Sym^2 W)$ and the cone structure is locally flat.

 \begin{proposition}\label{p.chi} Under ({\em Assumption}),
the group $G$ in Proposition \ref{p.autGrass} acts on $X$ with an
open orbit $X_o$ such that the complement $X\setminus X_o$ has
codimension $\geq 2$. There exists a $G$-biregular morphism $\chi:
\sW \otimes \sQ \to X_o$, sending the $Z$-isotrivial cone
structure of Proposition \ref{p.modelGrass} to the VMRT-structure on $X_o$, inducing a
fibration $\rho: X_o \to {\rm Lag}(\Sigma).$
\end{proposition}

\begin{proof} Since the isotrivial  cone structure on $X$ is locally flat,
it is locally equivalent to the cone structure $\sC$ of
Proposition \ref{p.modelGrass}. By Theorem \ref{t.Cartan-Fubini},
we have $\aut(X) = \aut(\sC, x) = \fg$ for $x \in \sW \otimes
\sQ$, which implies that the group $G$ acts on $X$ with an open
orbit $X_o$. As $\sW \otimes \sQ$ is simply connected,
 we have a $G$-equivariant
unramified covering morphism $\chi: \sW \otimes \sQ \to X_o$. The
image of the zero-section ${\rm Lag} (\Sigma) \subset (\sW \otimes
\sQ)$ is a positive-dimensional subvariety in $X_o$. Thus the
complement $X \setminus X_o$ must be of codimension $\geq 2$
because $X$ has Picard number 1. In particular, $X_o$ is simply
connected and the morphism $\chi: \sW \otimes \sQ \to X_o$ is
biregular. It  certainly  preserves the cone structure. The
fibration $\sW \otimes \sQ \to {\rm Lag}(\Sigma)$ induces a
fibration $\rho: X_o \to {\rm Lag}(\Sigma).$
\end{proof}

The following lemma is elementary. See  Proposition 4.4 of
\cite{HM01} for a proof.

\begin{lemma}\label{l.HartogsGrass}
Let $Y_1$ be a Fano manifold of Picard number one. Let $Y_2$ be a
compact complex manifold. Assume there exist subsets $E_i \subset
Y_i, i=1,2,$ of codimension $\geq 2$ and a biholomorphic morphism
$\varphi : Y_2 \setminus E_2 \; \to \; Y_1 \setminus E_1.$ Then
$\varphi$ can be extended to a biholomorphic morphism
$\tilde{\varphi}: Y_2 \to Y_1$.
\end{lemma}

The proof of the next proposition is essentially contained in the
proof of Proposition 6.3.3 of \cite{HM05}. We recall the proof for the reader's convenience.

\begin{proposition}\label{p.submanifold} Let $\rho: X_o\to
{\rm Lag}(\Sigma)$ be as in Proposition \ref{p.chi}. Given a point
$[W^*] \in {\rm Lag}(\Sigma)$, the closure in $X$ of the fiber
$\rho^{-1}([W^*])$  is a projective submanifold biregular to the
Grassmannian ${\rm Gr}(2, W^* \oplus Q)$ such that
$\rho^{-1}([W^*]) \subset \Gr(2,W^* \oplus Q)$ is isomorphic to
$W\otimes Q \subset \Gr(2,W^* \oplus Q)$. Consequently, the
biregular morphism
$$\chi: \sW \otimes \sQ = {\rm Gr}(2, \sW^* \oplus \sQ) \setminus D \; \to \;  X_o$$ in Proposition
\ref{p.chi}  can be extended to a morphism $\tilde{\chi}:  {\rm
Gr}(2, \sW^* \oplus \sQ) \to X$.
\end{proposition}

\begin{proof} As in Proposition \ref{p.modelGrass}, regard ${\rm Lag}(\Sigma)$ as
a submanifold of $$\sW \otimes \sQ \subset {\rm Gr}(2, \sW^*
\oplus \sQ).$$  From the description of the isotropy subgroup
$G_{[W^*]}$ in the proof of Proposition \ref{p.modelGrass}, we see
that $G_{[W^*]}$ contains a subgroup isomorphic to ${\rm GL}(W)
\times {\rm GL}(Q)$. Choose $\C^\times \subset {\rm GL}(W) \times
{\rm GL}(Q)$  with weight $1$ on $W$ and weight $-1$ on $Q$. Then
it has weight $0$ on $W\otimes Q$ and weight $2$ on $\Sym^2W$.
 It follows that this $\C^{\times}$ action on ${\rm
Gr}(2, \sW^* \oplus \sQ)$ fixes the point $[W^*]$ and the isotropy
action on
$$T_{[W^*]}({\rm Gr}(2, \sW^* \oplus \sQ)) = (W \otimes Q) \oplus
\Sym^2 W$$ fixes exactly $W\otimes Q$. Thus the fixed point set of
this $\C^{\times}$-action on ${\rm Gr}(2, \sW^* \oplus \sQ)$ has
the fiber ${\rm Gr}(2, W^* \oplus Q) $ as a connected component.
Consequently, the corresponding $\C^{\times}$-action on $X$ has
the closure $S \subset X$ of the fiber $\rho^{-1}([W^*])$ as a
connected component of its fixed point set. Since
 the fixed point set of a $\C^{\times}$-action on the projective
manifold $X$ is nonsingular, the closure $S$ is a projective
submanifold.

To show that this submanifold $S$ is biregular to the
Grassmannian, we need to show that the birational map $\delta:
S\dasharrow {\rm Gr}(2, W^* \oplus Q)$ induced by  $\chi^{-1}: X_o
\to {\rm Gr}(2, \sW^* \oplus \sQ)$ is biholomorphic. This is essentially Lemma 6.3.2 in \cite{HM05}. Let us recall the argument.

 For $Z \subset \BP((W\otimes Q) \oplus \Sym^2 W),$ let $Z' = Z \cap \BP (W \otimes Q)$, which is equivalent to a Segre embedding of $\BP W \times \BP Q$. The
$Z$-isotrivial cone structure on $\sW \otimes \sQ \subset {\rm
Gr}(2, \sW^* \oplus \sQ)$ induces a $Z'$-isotrivial cone structure
on ${\rm Gr}(2, W^* \oplus Q)$. This cone structure is exactly the
VMRT of lines on the Grassmannian. The $Z$-isotrivial cone
structure on $X_o$ also induces a $Z'$-isotrivial cone structure
on $\rho^{-1}([W^*])$. This cone structure is the VMRT of $S$
given by the minimal rational curves of $X$  lying on $S$. The map
$\delta$ induces an isomorphism of these $Z'$-isotrivial cone
structures. Thus $\delta$ sends minimal rational curves of $X$
lying on $S$ to lines in the Grassmannian ${\rm Gr}(2, W^* \otimes
Q)$.

Let $\sH \subset S$ be the union of hypersurfaces where $\delta$
is ramified (note that $\delta$ is always well-defined in
codimension 1).
 Suppose $\sH \neq \emptyset$.  Let $A$ be the proper image
of $\sH$ under $\delta$. Then $A$ is a subset of codimension $\geq
2$ in ${\rm Gr}(2, W^* \oplus Q)$. Choose a family of minimal
rational curves $\{\ell_s\; |\; s \in \Delta\}$ on ${\rm Gr}(2,
W^* \oplus Q)$ such that $\ell_0$ intersects $A$ but is not
contained in $A$; all $\ell_s$ with $s \neq 0$ are disjoint from
$A$ and are the strict images of a family of minimal rational
curves $C_s, s \neq 0,$ on $S$. Then the limit $C_0$ is an
irreducible curve because $C_0$ has degree 1 with respect to the
line bundle $L$ on $X$. This implies that the proper image of
$C_0$ must be $\ell_0$ and $C_0$ intersects $\sH$. But $C_s, s
\neq 0$ is disjoint from $\sH$. Thus we have a family of
irreducible curves $C_s, s \in \Delta,$ on the projective manifold
$S$ and a hypersurface $\sH \subset S$ such that $C_0 \cdot \sH
\neq 0$ but $C_s \cdot \sH =0$ for $s \neq 0$, a contradiction. We
conclude that $\sH = \emptyset.$

Since $\sH = \emptyset$, we see that $\delta$ is unramified
outside a subset $E \subset S$ of codimension $\geq 2$. The image
$\delta(S \setminus E) \subset {\rm Gr}(2, W^* \oplus Q)$ is not
an affine subset, because $S \setminus E$ contains projective
curves (general minimal rational curves of $S$). But $\delta(S
\setminus E)$ contains
 the open subset $W \otimes Q \subset  {\rm Gr}(2, W^* \oplus Q)$
  and its complement ${\rm Gr}(2, W^* \oplus Q) \setminus (W \otimes Q)$
  is an irreducible hypersurface. Thus the complement ${\rm Gr}(2, W^* \oplus Q)
 \setminus \delta(S \setminus E)$ is of codimension $\geq 2$. By Lemma \ref{l.HartogsGrass},
 $\delta$ extends to a biregular morphism $S \to {\rm Gr}(2, W^* \oplus Q).$
\end{proof}

\begin{proposition}\label{p.contract} In the setting of Proposition \ref{p.submanifold},
let ${\rm Gr}(2, \sQ) \subset {\rm Gr}(2, \sW^* \oplus \sQ)$ be
the trivial fiber subbundle whose fiber over $[W^*] \in {\rm
Lag}(\Sigma)$ corresponds to ${\rm Gr}(2, Q) \subset {\rm Gr}(2,
W^*\oplus Q)$ of Proposition \ref{p.actionGrass}. Then the
restriction $\tilde{\chi}|_{{\rm Gr}(2, \sQ)}$ agrees with the
projection
$${\rm Gr}(2, \sQ) = {\rm Gr}(2, Q)
\times {\rm Lag}(\Sigma) \to {\rm Gr}(2, Q).$$
\end{proposition}

\begin{proof}
The $\C^{\times} \subset G$  corresponding to the center of ${\rm
GL}(Q) \subset G$  acts on ${\rm Gr}(2, \sW^* \oplus \sQ)$ such
that on each fiber it induces the $\C^{\times}$-action of
Proposition \ref{p.actionGrass}. From Proposition
\ref{p.actionGrass},   ${\rm Gr}(2, \sQ)$ is a component of the
fixed point set of this action such that all general orbits in the
divisor $D$ have limit points in ${\rm Gr}(2, \sQ)$.

The morphism $\tilde{\chi}: {\rm Gr}(2, \sW^* \oplus \sQ) \to X$
defined in Proposition \ref{p.submanifold} sends the divisor $D$
to $X \setminus X_o,$ a subset of codimension $\geq 2$ in $X$ from
Proposition \ref{p.chi}. Let $A \subset D$ be a general fiber of
the contraction $\tilde{\chi}|_D: D \to X \setminus X_0$. The
limit of $A$ under the $\C^{\times}$-action contains a
positive-dimensional subvariety $A'$ in ${\rm Gr}(2, \sQ)$. By the
$\C^{\times}$-equivariance, $A'$ must be contracted by
$\tilde{\chi}$, too.   But the action of ${\rm GL}(Q) \subset G$
is transitive on ${\rm Gr}(2, \sQ_{[W^*]}) = {\rm Gr}(2, Q)$ for
each $[W^*] \in {\rm Lag}(\Sigma)$. Thus $\tilde{\chi}(\Gr(2,
\sQ))$ has
 dimension strictly less than that of ${\rm Gr}(2, \sQ)$, i.e., ${\rm Gr}(2, \sQ)$ is contracted by
 $\tilde{\chi}$.  By the definition of $\tilde{\chi}$ in
Proposition \ref{p.submanifold}, the line bundle $\tilde{\chi}^*
L$
 is ample on the ${\rm Gr}(2,Q)$-factor of $$ {\rm Gr}(2, \sQ) = {\rm Gr}(2, Q)
\times {\rm Lag}(\Sigma).$$  Thus the fibers of
$\tilde{\chi}|_{{\rm Gr}(2, \sQ)}$ must be contained in the ${\rm
Lag}(\Sigma)$-factor. Since ${\rm Lag}(\Sigma)$ has Picard number
1,
 $\tilde{\chi}$ must contract ${\rm Lag}(\Sigma)$ to one point. \end{proof}

\begin{proposition}\label{p.slice}
Pick a subspace $Q' \subset Q$ of dimension 2, defining a fiber
subbundle ${\rm Gr}(2, \sW^* \oplus \sQ')$ of ${\rm Gr}(2, \sW^*
\oplus \sQ)$. Then the image of ${\rm Gr}(2, \sW^* \oplus \sQ')$
under $\tilde{\chi}$ is a 7-dimensional  projective submanifold
$X' \subset X$ such that the restriction of $\tilde{\chi}$ to
${\rm Gr}(2, \sW^* \oplus \sQ')$
$$\mu: {\rm Gr}(2, \sW^* \oplus \sQ') \to X'$$ sends each fiber of
$\psi: {\rm Gr}(2, \sW^* \oplus \sQ') \to {\rm Lag}(\Sigma)$
isomorphically to a projective submanifold of $X'$ and contracts
the submanifold $$[Q'] \times {\rm Lag}(\Sigma) = {\rm Gr}(2,
\sQ') \subset {\rm Gr}(2, \sW^* \oplus \sQ')$$ to one point in
$X'$.
\end{proposition}

\begin{proof}
From Propositions \ref{p.submanifold} and  \ref{p.contract}, all
are obvious except the nonsingularness of the image $X'$. To see
this, fix a decomposition $Q = Q' \oplus Q^{''}$ and choose a
$\C^{\times} \subset {\rm GL}(Q)$ which acts with weight 0 on $Q'$
and weight 1 on $Q^{''}$. The induced $\C^{\times}$-action on
${\rm Gr}(2, \sW^* \oplus \sQ)$ has ${\rm Gr}(2, \sW^* \oplus
\sQ')$ as a component of its fixed point set. Since the morphism
$\tilde{\chi}$ is equivariant under this $\C^{\times}$-action on
${\rm Gr}(2, \sW^* \oplus \sQ)$ and the corresponding
$\C^{\times}$-action on $X$, the image $X'$ is a component of the
fixed point set of this $\C^{\times}$-action. Thus $X'$ is
nonsingular.
\end{proof}

\begin{proof}[End of the proof of Theorem \ref{t.non1}] Let $\iota \in X'$ be the image
$\mu(\Gr(2, \sQ'))$ in Proposition \ref{p.slice}.  The group ${\rm
Sp}(\Sigma)$, with Lie algebra $\fsp(4)=\fso(5)$ acts on  ${\rm
Gr}(2, \sW^* \oplus \sQ')$ preserving ${\rm Gr}(2, \sQ')$. Thus it
acts on $X'$ with $\iota$ fixed, inducing the isotropy
representation of $\fso(5)$ on $T_{\iota}(X') $. This
representation is non-trivial as a non-trivial action of a
reductive group gives a non-trivial isotropy action on the tangent
space of a fixed point.
 As non-trivial irreducible representations of $\fso(5)$ of dimension
$\leq 7$ can be either of dimension 4 (the spin representation) or
5 (the standard representation), the fixed point set of this
$\fso(5)$-action
 has
a component $E \subset X'$ with $\dim E =  3$ or $2$ through
$\iota$.

For any $[W^*] \in {\rm Lag}({\Sigma})$, the isotropy subgroup in
${\rm Sp}(\Sigma)$ contains the subgroup ${\rm GL}(W^*)$, which
acts in a natural way on
 $\Gr(2, W^* \oplus Q')$.  The
fixed point set  of this  ${\rm GL}(W^*)$-action consists of two
isolated points: $[W^*]$ and $[Q']$. As $\Gr(2, W^* \oplus Q')$ is
mapped isomorphically and equivariantly to a projective submanifold of $X'$, the germ of $E$ at $\iota$  intersects
this image submanifold   only at the point $\iota$. As this is true for
all $[W^*] \in {\rm Lag}(\Sigma)$ and the union of all such images
is $X'$, we deduce that $E = \iota$,
 a
contradiction to the dimension of $E$.
\end{proof}

\section{Proof of Theorem \ref{t.non2}}\label{s.non2}

The section is devoted to the proof of Theorem \ref{t.non2}. The
argument is overall parallel to that of Section \ref{s.non1},
replacing Grassmannians by hyperquadrics. In fact, Proposition 8.i
is a direct analogue of Proposition 7.i, etc.

To start with, let us recall some facts about hyperquadrics.
 By an
orthogonal vector space we mean a vector space $U$ equipped with a
non-degenerate quadratic form $\beta$. Given an orthogonal vector
space, the hyperquadric $\Q(U) \subset \BP U$ is the set of
null-vectors, i.e., its affine cone is
$$\widehat{\Q(U)} := \{ u\in  U \; | \; \beta(u,u) = 0 \}.$$

\begin{lemma}\label{l.actionQ}
Let $S$ be an orthogonal vector space with a quadratic form $\alpha$.
Define a 2-dimensional orthogonal space $(\C \oplus \C, \gamma)$
by the multiplication $\gamma(s,t) = st \in \C$. The direct sum
$(S \oplus (\C \oplus \C), \alpha \oplus \gamma)$ is an orthogonal
space. Consider the hyperquadric $\Q(S \oplus (\C \oplus \C))$ of
this orthogonal space.  There is a natural embedding of $S$ into
$\Q(S \oplus (\C\oplus \C))$ as a Zariski open subset whose
complement $D$ is an irreducible divisor defined by $t=0$. The
divisor $D$ has a unique singular point, to be denoted by
$\Gamma.$  The $\C^{\times}$-action on $S$ given by the scalar
multiplication extends to a $\C^{\times}$-action on $\Q(S \oplus
(\C \oplus \C))$ such that a general orbit in $D$ has $\Gamma$ as
a limit point.
%Moreover, if $E \subset D$ is a submanifold stable under the
%$\C^{\times}$-action with $\Gamma \in E$, then $\dim E \leq
%\frac{\dim S}{2}.$
\end{lemma}

\begin{proof}
Choose coordinates $z_1, \ldots, z_n$ on $S$ with respect to which
the quadratic form $\alpha$ is given by $z_1^2 + \ldots + z_n^2.$
In terms of the homogeneous coordinates $[z_1, \ldots, z_n, s, t]$
on $\BP (S \oplus (\C \oplus \C))$, the hyperquadric $\Q(S \oplus
(\C \oplus \C))$ is defined by
$$z_1^2 + \ldots + z_n^2 + st = 0.$$ The open embedding of $S$
is given by $$t=1, \; s = -z_1^2 - \cdots - z_n^2.$$ Its
complement is the divisor $D$ defined by $t=0$ and $D$ has a
unique singular point $$ \Gamma := (z_1 = \cdots = z_n =t = 0).$$
The $\C^{\times}$-action of scalar multiplication on $S$ is given
in this coordinates  as the action of  $\lambda \in \C^{\times}$
by
$$(z_1, \ldots, z_n, s, t) \mapsto (\lambda z_1, \ldots, \lambda
z_n, \lambda^2 s,t ).$$ This certainly induces a
$\C^{\times}$-action on $\Q(S \oplus (\C \oplus \C))$ preserving
$D$. For any point $(z_1, \ldots, z_n, s, 0) \in D$ with $s \neq
0$, the orbit $$\{ [\lambda z_1, \ldots, \lambda z_n, \lambda^2 s,
0)], \lambda \in \C^{\times} \}  = \{ [\lambda^{-1} z_1, \ldots,
\lambda^{-1} z_n, s, 0], \lambda  \in \C^{\times} \}$$ has
$\Gamma$ as a limit point as $\lambda^{-1}$ approaches $0$.
\end{proof}

Next we need to look at the geometry of a certain hyperquadric
bundle over a 7-dimensional hyperquadric.
  Fix a 9-dimensional orthogonal vector space $U$. The hyperquadric $\Q(U)$
  is a 7-dimensional projective manifold homogeneous under ${\rm SO}(U)$.
  The semi-simple part of the isotropy group at a point of $\Q(U)$ has Lie algebra $\fso(7).$
  The 8-dimensional spin representation $W$ of $\fso(7)$ induces
  a homogeneous vector bundle $\sS^*$ of rank 8 on $\Q(U)$, called
  the {\em dual spinor bundle} and its dual is called the {\em spinor bundle} $\sS$. See \cite{Ot} for details.

  \begin{proposition}\label{p.Ottaviani}
  (i) Denoting by $L$ the ample generator of ${\rm Pic}(\Q(U))$, we
  have $ \sS^* \cong \sS \otimes L$ and $$
  H^0(\Q(U), \sS^*\otimes L^{-1})=H^1(\Q(U), \sS^*\otimes
  L^{-1}) = H^1(\Q(U), \sS^*)=0.$$
(ii) For all $i \geq 0$,  $$H^i(\Q(U), \wedge^2 \sS^* \otimes
L^{-1}) = 0.$$

  (iii) The global sections of $\sS^*$ generate the vector bundle $\sS^*$  and
   $H^0(\Q(U), \sS^*)$
   is the 16-dimensional spin
  representation of $\fso(U) = \fso (9).$\end{proposition}
\begin{proof}
Claim (i) and (iii) follow from Theorem 2.3 and Theorem 2.8 of
\cite{Ot}. We now prove claim (ii). Denoting by $Q$ the
 standard $7$-dimensional representation of $\fso(7)$, we have $\wedge^2 W \cong Q
\oplus \wedge^2 Q$ as representations of $\mathfrak{so}(7)$. Let
$P$ be the isotropy group of a point on $\mathbb{Q}(U)$. As the
center of $P$ acts trivially on $W$, we see that, as
$P$-representations, the highest weight of  $Q$ (resp. $\wedge^2
Q$) is  $\lambda_2$ (resp. $\lambda_3$), where $\lambda_i$ is the
$i$-th fundamental weight of the simple Lie algebra of type $B_4$.
 Note that the line bundle $L^{-1}$ is induced by
the representation of highest weight $-\lambda_1$. This gives that
the bundle $\wedge^2 \mathcal{S}^* \otimes L^{-1}$ is a direct sum
of two equivariant vector bundles with highest weights
$\lambda_2-\lambda_1$ and $\lambda_3-\lambda_1$. Let $\delta$ be
the sum of all fundamental weights, then we see $\delta +
\lambda_2-\lambda_1$ and $\delta + \lambda_3-\lambda_1$ contains
no $\lambda_1$, i.e. these sums are singular weights.
 This implies $H^i(\mathbb{Q}(U), \wedge^2 \mathcal{S}^* \otimes
L^{-1})=0$ for all $i \geq 0$ by Borel-Weil-Bott's theorem.
\end{proof}

Note that the spin representation of $\fso(7)$ carries an
invariant non-degenerate quadratic form (e.g. \cite{FH} Exercise
20.38). Thus there exists a fiberwise non-degenerate quadratic
form on $\sS^*$ with values in a line bundle $M$ on $\Q(U)$, i.e.,
$ \Sym^2 \sS^* \to M$ inducing an isomorphism $(\sS \otimes M)
\cong \sS^*$. From Proposition \ref{p.Ottaviani}, we have $$ \sS
\otimes M \; \cong \; \sS \otimes L,$$ implying $M = L$.
Consequently, we have a fiberwise non-degenerate quadratic form
$$\alpha: \Sym^2(\sS^*) \to L.$$ On the other hand the natural
multiplication $L \otimes \sO \to L$, where $\sO=\sO_{\Q(U)},$
induces a fiberwise non-degenerate quadratic form
$$\gamma: \Sym^2( L \oplus \sO) \to L.$$ Thus the vector bundle
$\sS^* \oplus (L \oplus \sO)$ of rank 10 is equipped with the
fiberwise non-degenerate quadratic form
$$\alpha \oplus \gamma: \; \Sym^2(\sS^*  \oplus (L \oplus \sO)) \to L.$$
The associated hyperquadric bundle $$ \psi : \Q(\sS^* \oplus (L
\oplus \sO)) \; \to \; \Q(U)$$ is a fiber bundle on $\Q(U)$ whose
fiber is an 8-dimensional hyperquadric. We will denote this
projective manifold $\Q(\sS^* \oplus (L \oplus \sO))$ by $Y$.

\begin{proposition}\label{p.autQ} Let $\psi: Y \to \Q(U)$ be the hyperquadric bundle  of the orthogonal vector bundle
 $\sS^*  \oplus (L \oplus \sO)$. Let $\Xi=H^0(\Q(U), \sS^*)$ be the 16-dimensional spin representation of $\fso(U).$
 Then Lie algebra $\fg$ of the
automorphism group of the projective variety $Y$ is isomorphic to
$ \Xi \sd (\fso(U) \oplus \C),$ where $\C$ corresponds to the
scalar multiplication on the vector bundle $\sS^*$. The vector
bundle $\sS^*$ has a natural embedding into $Y= \Q(\sS^* \oplus (L
\oplus \sO))$. Its complement $D$ is an irreducible divisor and
the singular locus of $D$ is a section $\Gamma$ of $\psi$.
\end{proposition}

\begin{proof}
By Proposition \ref{p.Ottaviani} (iii), we have a surjective map
$\Xi \otimes \mathcal{O} \to \sS^*$, which gives for any $x \in
\Q(U)$ a surjective map $\zeta_x: \Xi \to \sS^*_x$. The vector
group $\Xi$ acts on $\sS^* \oplus L \oplus \mathcal{O}$ by the
following rule: for any $(v, s, t) \in \sS^*_x \oplus L_x \oplus
\mathcal{O}_x$ and any $f \in \Xi$,
$$
f \cdot (v, s, t) =  (v + t \zeta_x(f), s-2 \alpha(v,
\zeta_x(f))-t \alpha(\zeta_x(f), \zeta_x(f)), t).
$$
One checks easily that this action preserves the quadric form on
$\sS^* \oplus L \oplus \mathcal{O}$. This induces an action of
$\Xi$ on $Y$. From this, we see that there is a natural inclusion
$$
 \Xi \sd (\fso(U) \oplus \C) \subset \fg.
$$
To show that this is an isomorphism, it suffices to compare their
dimensions.

Let $\psi: Y \to \Q(U)$ be the natural projection. We have an
exact sequence
\begin{equation}\label{eq.F4}
0 \to T^{\psi} \to T(Y) \to \psi^* T(\Q(U)) \to 0,
\end{equation}
 where $T^{\psi}$ denotes the relative tangent
bundle. Recall  that for an orthogonal vector
space $\C^m$, there is a natural identification $$H^0(\Q(\C^m),
T(\Q(\C^m))) = \fso(\C^m) = \wedge^2 \C^m.$$  Translating it into
relative setting, we get $$\psi_* T^{\psi} = \wedge^2(\sS^* \oplus
(L \oplus \sO)) \otimes L^{-1} = ((\wedge^2 \sS^*) \otimes L^{-1})
 \oplus \sO \oplus \sS^* \oplus (\sS^*\otimes L^{-1}). $$ By $R^i \psi_* T^{\psi} =0$ for $i \geq 1$ and  Proposition \ref{p.Ottaviani},
we have $H^1(Y, T^{\psi})=0$ and
$$H^0(Y, T^{\psi}) =H^0(\Q(U), \psi_* T^{\psi}) =  H^0(\Q(U), \sO \oplus \sS^*) = \C \oplus
\Xi.$$  Since $H^0(\Q(U), T(\Q(U)) = \fso(U)$, the long-exact
sequence associated to \eqref{eq.F4} shows that $$ \dim \fg = \dim H^0(Y, T(Y)) = \dim (\Xi
\sd (\fso(U) \oplus \C)).$$

Now the rest of Proposition \ref{p.autQ} is a globalization of
Lemma \ref{l.actionQ}. The hyperquadric bundle $\psi: Y \to
\Q(U)$ has a natural section $\Gamma \subset Y$ over $\Q(U)$
determined by
$$ \Gamma:= \BP \sO \subset \Q(\sS^*  \oplus (L \oplus \sO))=Y$$
because the $\sO$-factor of $\sS^*  \oplus (L \oplus \sO)$ is a
null-vector with respect to the quadratic form $\alpha \oplus
\gamma.$ Given a point $v \in \sS^*$, let $v' \in L$ be the unique
vector defined by
$$\alpha(v,v) + \gamma(v',1) = 0$$ where $1$ denotes the
section of $\sO$ determined by the constant function 1 on $\Q(U)$.
 Then we have a canonical embedding of $\sS^* $ into the hyperquadric
 bundle $Y = \Q(\sS^*  \oplus (L \oplus \sO))$ as a Zariski open subset
 by $$ v \in \sS \mapsto (v, (v', 1)).$$ Its complement is an
 irreducible divisor $D$ determined by the zero section of $\sO$ and $\Gamma$ is the singular locus of $D$,
 which can be seen immediately from Lemma \ref{l.actionQ}.
\end{proof}

\begin{proposition}\label{p.modelQ}
Let $G$ be the simply connected group with Lie algebra $\fg$ of
Proposition \ref{p.autQ}. The open subset $\sS^* \subset Y$
described in Proposition \ref{p.autQ} is $G$-homogeneous and has a
natural isotrivial cone structure $\sC$ invariant under the
$G$-action such that each fiber $\sC_x \subset \BP T_x(\sS^*)$ is
isomorphic to
 $Z \subset \BP (W \oplus Q)$ in the notation of Section \ref{e.S_5}.
 This cone structure $\sC$ is locally
flat and $\aut(\sC, x) \cong \fg$ for each $x \in \sS^*$.
\end{proposition}

\begin{proof}
It is easy to see that the open subset $\sS^*$ is $G$-invariant.
The base $\Q(U)$ is homogenous under the action of ${\rm SO}(U)$.
From the proof of Proposition \ref{p.autQ}, the vector group $\Xi$
acts on the fiber $\sS^*_x$ by translation of images of $\zeta_x$,
thus this action is transitive on the fibers of $\sS^* \to \Q(U)$.
This shows that $\sS^*$ is $G$-homogeneous.

 For a point $z \in \Q(U)$, the Lie algebra of the isotropy subgroup in
${\rm SO}(U)$ of $z$  is a parabolic subalgebra $\fp_{z} \subset
\fso(U)$. It is known that the reductive part of $\fp_z$ is
isomorphic to $\fco(7)$.  Regard $\Q(U)$ as a submanifold of
$\sS^*\subset Y$ via the zero section of the vector bundle.  Let
$\Xi_z : = \Ker(\zeta_z) \subset \Xi$, which corresponds to the
sections of $\sS^*$ vanishing at $z$. At the point $z \in \Q(U)$
 the isotropy subgroup $G_{z}$ has Lie algebra
$$\fg_{z}:= \Xi_z \sd (\fp_{z} \oplus \C)
\; \subset \; \Xi \sd (\fso(U) \oplus \C) = \fg.$$
The tangent space
$$T_{z}(\sS^*) = \sS^*_z \oplus T_{z}(\Q(U)) \cong W \oplus Q $$ contains the affine cone $\hat{Z}$ in a natural way.
 The isotropy representation of $\fg_{z} =\Xi_z \sd (\fp_{z} \oplus \C)$ on $T_{z}(\sS^*) =W \oplus  Q$  satisfies
\begin{itemize} \item[(1)]  $\Xi_z$-component of
$\fg_{z}$ acts trivially; \item[(2)]  $\fp_{z}$-factor acts as
$\fco(7)$ in a natural way
 on $W$ and on $ Q$. \item[(3)] The
$\C$-factor has weight 1 on $W$ and weight 0 on $Q$.
\end{itemize} From Proposition \ref{p.autZQ},  $\hat{Z}$ is
preserved under the isotropy representation of the isotropy
subgroup $G_{z}$. Thus the $G$-action defines a natural
$Z$-isotrivial cone structure on the open set $\sS^*$.

As $\aut(\hat{Z})^{(2)} = 0$ by Theorem \ref{t.HMk=2}, we have the
following inequalities from Proposition \ref{p.dim}
$$
\dim \fg \leq \dim \aut(\sC, x) \leq \dim (\aut(\hat{Z})^{(1)}
\oplus \aut(\hat{Z}) \oplus W \oplus  Q).
$$
From Propositions \ref{p.autZQ} and \ref{p.prolongZQ},  $$ \dim
(\aut(\hat{Z})^{(1)} \oplus \aut(\hat{Z}) \oplus (W \otimes Q)
\oplus \Sym^2 W) = 53= \dim \fg$$ implying $\fg \cong
\aut(\sC,z).$ Then Corollary \ref{c.dim} shows that  the
$Z$-isotrivial cone structure on $\sS^* \subset Y$ is locally
flat.
\end{proof}

Now to prove Theorem \ref{t.non2}, we will make the following
assumption and derive a contradiction.

({\em Assumption}) Let $X$ be a 15-dimensional Fano manifold with
${\rm Pic}(X) = \Z \langle \sO_X(1) \rangle$. Assume that $X$ has
minimal rational curves of degree 1 with respect to $\sO_X(1)$
whose VMRT at a general point is isomorphic to $Z \subset \BP  (W
\oplus Q)$ and the cone structure is locally flat.

 \begin{proposition}\label{p.chiQ} Under ({\em Assumption}),
the group $G$ in Proposition \ref{p.modelQ} acts on $X$ with an
open orbit $X_o$ such that the complement $X\setminus X_o$ has
codimension $\geq 2$. There exists a $G$-biregular morphism $\chi:
\sS^*   \to X_o$, sending the $Z$-isotrivial cone structure of
Proposition \ref{p.modelQ} to the $Z$-isotrivial VMRT cone
structure on $X$. This  induces  a fibration $\rho: X_o \to
\Q(U).$
\end{proposition}

 \begin{proof}
Since the $Z$-isotrivial  cone structure on $X$ is locally flat,
it is locally isomorphic to the cone structure $\sC$ of
Proposition \ref{p.modelQ}. By Theorem \ref{t.Cartan-Fubini}, we
have $\aut(X) = \aut(\sC, x) = \fg$ for $x \in \sS^*$ general,
which implies that the group $G$ acts on $X$ with an open orbit
$X_o$. As $\sS^*$ is simply connected,
 we have a $G$-equivariant unramified covering morphism $\chi:
\sS^* \to X_o$. The image of the zero-section $\Q(U) \subset
\sS^*$ is a positive-dimensional subvariety in $X_o$. Thus the
complement $X \setminus X_o$ must be of codimension $\geq 2$
because $X$ has Picard number 1. In particular, $X_o$ is simply
connected and the morphism $\chi: \sS^* \to X_o$ is biregular. It
certainly preserves the cone structure. The fibration $\psi: \sS^* \to
\Q(U)$ induces a fibration $\rho: X_o \to \Q(U).$
\end{proof}

The proof of the next proposition is essentially the same as
that of Proposition 8.3.4 of \cite{HM05}. We recall the proof
for the reader's convenience.

\begin{proposition}\label{p.submanifoldQ} Let $\rho: X_o\to
\Q(U)$ be as in Proposition \ref{p.chiQ}. Given a point $z \in
\Q(U)$, the closure in $X$ of the fiber $\rho^{-1}(z)$ is a
projective submanifold biregular to the hyperquadric
$\psi^{-1}(z)$ such that $\rho^{-1}(z)$ corresponds to $\sS^*_z
\subset \psi^{-1}(z)$. Consequently, the biregular morphism
$$\chi: \sS^*  = Y\setminus D \; \to \;  X_o$$ in Proposition
\ref{p.chiQ}  can be extended to a morphism $\tilde{\chi}: Y \to
X$.
\end{proposition}

\begin{proof} As in the proof of  Proposition \ref{p.modelQ}, regard $\Q(U)$ as
a submanifold of $\sS^* \subset Y$.  From the description of the
isotropy subgroup $G_{z}$ in the proof of Proposition
\ref{p.modelQ},  we see that $\fg_{z}$ contains a subalgebra
isomorphic to $\fco(Q) \subset \fp_z$, whose center has weight 1
on both $W$ and  $Q$. Also, there is a $\C$-factor in $\fg_z$ with
weight 1 on $W$ and 0 on $Q$. This implies that there exists a
subgroup $\C^{\times} \subset G$ which acts with weight $1$ on $Q$
and weight $0$ on $W$. It follows that this $\C^{\times}$ action
on $Y$ fixes the point $z$ and the isotropy action on
$$T_x(Y) = W \oplus  Q$$ fixes exactly $W$. Thus the fixed point set of this $\C^{\times}$-action on
$Y$ has the fiber $Y_z:=\psi^{-1}(z) $ as a connected component.
Consequently, the corresponding $\C^{\times}$-action on $X$ has
the closure $S_z \subset X$ of the fiber $\rho^{-1}(z)$ as a
connected component of its fixed point set. Since
 the fixed point set of a $\C^{\times}$-action on the projective
manifold $X$ is nonsingular, the closure $S_z$ is a projective
submanifold.

To show that this submanifold $S_z$ is biregular to the
hyperquadric $\psi^{-1}(z)$, we need to show that the birational
map $\delta: S_z \dasharrow \psi^{-1}(z)$ induced by  $\chi^{-1}:
X_o \to Y$ is biholomorphic.

 Recall that $Z' = Z \cap \BP W$ is a 6-dimensional hyperquadric $\Q(W)$ determined
 by the orthogonal structure on the 8-dimensional spin representation $W$. The
$Z$-isotrivial cone structure on $\sS^* \subset Y$ induces a
$Z'$-isotrivial cone structure on the fiber $\sS^*_z$. This cone
structure is exactly the VMRT of lines on the hyperquadric. The
$Z$-isotrivial cone structure on $X_o$ also induces a
$Z'$-isotrivial cone structure on $\rho^{-1}(z)$. This cone
structure is the VMRT of $S_z$ given by the minimal rational
curves of $X$  lying on $S_z$. The map $\delta$ induces an
isomorphism of these $Z'$-isotrivial cone structures. Thus
$\delta$ sends minimal rational curves of $X$ lying on $S_z$ to
lines in the hyperquadric $\psi^{-1}(z)$. Then the same argument
as in the proof of Proposition \ref{p.submanifold},  shows that
 $\delta$ extends to a biregular morphism $S_z \to \psi^{-1}(z).$
\end{proof}

\begin{proposition}\label{p.contractQ} In the setting of Proposition \ref{p.submanifoldQ},
let $\Gamma \subset Y$ be the section of $\psi$ given by the
singular locus of the divisor $D$.  Then $\tilde{\chi}(\Gamma)$ is
one point.
\end{proposition}

\begin{proof}
We can choose a subgroup $\C^{\times} \subset G$ which corresponds
to the scalar multiplication of the vector bundle $\sS^*$,
corresponding to the $\C^{\times}$-action of Lemma \ref{l.actionQ}
on each fiber of $\psi$. From Lemma \ref{l.actionQ}, $\Gamma$ is a
component of the fixed point set of this action such that all
general orbits in the divisor $D$ have limit points in $\Gamma$.

The morphism $\tilde{\chi}: Y \to X$ defined in Proposition
\ref{p.submanifoldQ} sends the divisor $D$ to $X \setminus X_o,$ a
subset of codimension $\geq 2$ in $X$ from Proposition
\ref{p.chiQ}. Let $A \subset D$ be a general fiber of the
contraction $\tilde{\chi}|_D: D \to X \setminus X_0$. The limit of
$A$ under the $\C^{\times}$-action contains a positive-dimensional
subvariety $A'$ in $\Gamma$. By the $\C^{\times}$-equivariance,
$A'$ must be contracted by $\tilde{\chi}$, too.   But $\Gamma$ is
an orbit of the action of a subgroup of $G$ with Lie algebra
$\fso(U) \subset \fg.$ Thus $\Gamma$ is contracted by
 $\tilde{\chi}$. Since   $\Gamma \cong \Q(U)$ has Picard number
1,
 $\tilde{\chi}$ must contract $\Gamma$ to one point. \end{proof}

\begin{proof}[End of the proof of Theorem \ref{t.non2}] Let $\iota \in X$ be the image
$\tilde{\chi}(\Gamma)$ in Proposition \ref{p.contractQ}.  The
group ${\rm Spin}(9) \subset G$  acts on  $Y$ preserving $\Gamma$.
Thus it acts on $X$ with $\iota$ fixed, inducing the isotropy representation of $\fso(9)$  on
$T_{\iota}(X) $. This representation is non-trivial, because a
non-trivial action of a reductive group gives a non-trivial
isotropy action on the tangent space of a fixed point. Since an
irreducible representation of $\fso(9)$ with dimension $\leq 15$
must be the 9-dimensional standard representation,  $T_{\iota}(X)$
decomposes as a $\fso(9)$-module into the sum of the orthogonal
space $U$ and a complementary subspace of dimension 6 where
$\fso(9)$ acts trivially. This implies that the fixed point set of
the ${\rm Spin}(9)$-action on $X$ has a component $E$ of dimension
6 through $\iota$.

For any $z \in \Q(U)$, the stabilizer of ${\rm Spin}(9)$ contains
the subgroup ${\rm Spin}(7)$, which acts in a natural way on the
hyperquadric $Y_z = \psi^{-1}(z)$. This action when restricted to
$\sS^*_z$ is  the spin representation, which has no fixed point in
$\sS^*_z$ except the zero point $z$. This action on $D_z : = D
\cap Y_z$ has only one isolated fixed point, which is its singular
point $\Gamma \cap Y_z$. As $Y_z$ is mapped isomorphically and
equivariantly to its image in $X$, the germ of $E$ at $\iota$ intersects  this image
only at the  point $\iota$. As this holds for all $z$ and the union
of such images covers $X$, we deduce that $E = \iota$, a
contradiction to $\dim E=6$.
\end{proof}

\section{Application to target rigidity}\label{s.target}

In this section, we will give an application of the following
corollary of
 Main Theorem and Theorem \ref{t.mainprolong}.

\begin{corollary} \label{c.projection}
Let $S \subset \BP V$ be a nonsingular  non-degenerate projective
variety with $\Sec(S) \neq \BP V$. Then
$\aut(\widehat{p_x(S)})^{(1)} =0$ for a general point $x \in \BP
V$.
\end{corollary}

\begin{proof}
Suppose that $\aut(\widehat{p_x(S)})^{(1)} \neq 0$.  From Main
Theorem, $p_x(S) \subset \BP (V/\hat{x})$ must be a biregular
projection of the linearly normal embedding $S \subset \BP W, W =
H^0(S, \sO(1))^*, $ of  the varieties in (A1), (A2), (A3) or (B3) in Main
Theorem. Since $x$ is general, it is a biregular projection from a
subspace $L \subset W$ passing through a general point. This
contradicts Theorem \ref{t.mainprolong}. \end{proof}

Our application is via the following cone structure.

\begin{definition}\label{d.S} Let $S \subset \BP V$ be a nonsingular
 projective variety such that ${\rm Sec}(S) \neq \BP V$. On $M := \BP V \setminus {\rm Sec}(S)$,
 we define a cone structure $\sC \subset \BP T(M)$ as follows:
for each $x \in M$, let $\sC_x \subset \BP T_x(M)$ be the union of
tangents to lines joining $x$ to points of $S$. This cone
structure will be called the {\em cone structure induced by } $S$.
For a point $x \in \BP V,$ denote by $\hat{x} \subset V$ the
1-dimensional subspace over $x$ and $p_x: \BP V \setminus \{x\}
\to \BP (V/\hat{x})$ the projection. If $x \in M$, then $p_x|_{S}:
S \to p_x(S)$ is a biregular projection, embedding  $S$  into $\BP
(V/\hat{x})$. The projective variety $\sC_x \subset \BP T_x(M)$ is
isomorphic to $p_x(S) \subset \BP (V/\hat{x})$.
\end{definition}

\begin{theorem}\label{t.1} In Definition \ref{d.S}, suppose that  $S \subset \BP V$ is
non-degenerate and linearly normal.  Then for the cone structure
$\sC \subset \BP T(M)$ induced by $S$ and a general point $x \in
M$, $\aut(\sC,x) \cong \aut(S).$
\end{theorem}

\begin{proof}
By Corollary \ref{c.projection}, we have
$\aut(\widehat{p_x(S)})^{(1)} =0$ for a general point $x \in \BP
V$. From Proposition \ref{p.dim}, this implies that for a general
point $x \in M$,
$$ \dim(\aut(\sC, x)) \leq \delta(\sC) +
\dim \aut(\hat{\sC}_x).$$ By Lemma \ref{l.Aut}, for the cone
structure induced by $S$, the isotrivial leaf through a point $x$
is exactly the orbit of $x$ under the projective automorphism group $\Aut(S)$, which implies that
$$\delta(\sC) = \dim \Aut(S) \cdot x = \dim \Aut(S) - \dim
\Aut(S,x),$$ where $\Aut(S,x) \subset \Aut(S)$ is the isotropy subgroup at $x$.  Since $\dim \aut(\hat{\sC}_x) = \dim \Aut(S, x) +1$
by Lemma \ref{l.Aut}, we have $ \dim(\aut(\sC, x)) \leq \dim
\Aut(S) +1.$ But we have $\aut(S) \subset \aut(\sC, x)$ because
$\Aut(S)$ preserves the cone structure $\sC$ on $M$. Thus $$ \dim
\Aut(S) \leq \dim(\aut(\sC, x)) \leq \dim \Aut(S) +1.$$  To prove
Theorem \ref{t.1}, we may assume $ \dim(\aut(\sC, x)) = \dim
\Aut(S) +1$ and derive a contradiction.

 The assumption means that  the equality holds in
$$ \dim(\aut(\sC, x)) \leq \delta(\sC) +
\dim \aut(\hat{\sC}_x).$$  Then by Corollary \ref{c.dim}, the cone
structure is isotrivial and locally flat, which implies that we
can choose a connected open subset $U \subset M$ and a
biholomorphic map $\psi: U \to W$ into a vector space $W$ of
dimension $\dim W= \dim M$, such that $\psi_* \sC$ is of the form
$\psi(U) \times \sC_{o} \subset \BP T(W)$ for some base point $o$.
In particular, $$\delta(\sC) = \dim M = \dim W = \dim \Aut(S) -
\dim \Aut(S,x),$$ implying that $\Aut(S)$ has an open orbit.
 The
action of $\Aut(S)$ on $M$ preserves the cone structure. Thus by
Proposition \ref{p.flat1}, we have a complex analytic injective
homomorphism induced by $\psi$ $$\tau: \Aut(S) \to
W \sd {\rm Aut}(\widehat{p_x(S)}).$$ As $\dim \Aut(S,x) +1 = \dim
\aut(\widehat{p_x(S)})$, the isotropy representation
$$\iota: \aut(S,x) \to \aut(\widehat{p_x(S)})$$ is not surjective.
From the definition of $\aut(S,x)$ and the commutative diagram
$$\begin{array}{ccc} \aut(S) & \longrightarrow & W \sd
\aut(\widehat{p_x(S)}) \\ \cup & & \downarrow \\ \aut (S, x) &
\stackrel{\iota}{\longrightarrow} & \aut(\widehat{p_x(S)}),
\end{array} $$
 the composition $$\Aut(S) \stackrel{\tau}{\longrightarrow}
W \sd {\rm Aut}(\widehat{p_x(S)}) \longrightarrow
\Aut(\widehat{p_x(S)})$$ is not surjective and must have a kernel
of dimension $\dim W$.  This means that $W \subset \tau(\Aut(S))$.
Then by Lemma \ref{l.algebraic}, $W \subset \Aut(S)$ as an
algebraic subgroup. Moreover, the orbit $W \cdot x$ must have
dimension $\dim(W)$ from the description of the $W$-action on $W$
as translations. Since the stabilizer of $x$ in $W$ is an
algebraic subgroup, it must be the identity, because the only
discrete algebraic subgroup of the vector group is the identity.
It follows that the constructible set $W\cdot x$ is biregular to
the affine space and $\BP V$ is an equivariant compactification of
the vector group. This implies that the complement of $W\cdot x$
is a hyperplane in $\BP V$ (e.g. Satz 3.1 in \cite{Ge}). But our
$S \subset \BP V$ must belong to the complement of $W \cdot x$.
This contradicts the non-degeneracy of $S \subset \BP V$.
\end{proof}

It is worth noticing the following consequence of the proof of
Theorem \ref{t.1}.

\begin{proposition}\label{p.non-degenerate}
Let $S \subset \BP V$ be a linearly normal nonsingular
non-degenerate subvariety such that $\Sec(S) \neq \BP V$. Let
$\sC$ be the cone structure induced by $S$ on $\BP V \setminus
\Sec(S)$. Then the cone structure $\sC$ is isotrivial if and only
if $\Aut(S)$ acts on $\BP V$ with an open orbit. In this case, the
cone structure $\sC$ is never locally flat.
\end{proposition}
\begin{proof}
The first claim follows directly from Lemma \ref{l.Aut}. If $\sC$
is isotrivial and locally flat, then the previous theorem gives
$\aut(S) = \aut(\sC, x)=W \sd \aut(\widehat{p_x(S)})$, which
implies that $\Aut(S)$ contains the vector group $W$. By the same
argument as in the proof of Theorem \ref{t.1}, this implies that
$S$ is degenerated, a contradiction.
\end{proof}

\begin{remark} In the setting of Proposition \ref{p.non-degenerate}, if $S$ is degenerate, i.e., $S$ is
contained in a hyperplane $\BP^{n-1} \subset \BP^n$, then the cone
structure on the complement $\BP^n\setminus \BP^{n-1}$ given by
lines intersecting $S$  is isotrivial and locally flat with
$\sC_x$ projectively equivalent to $S \subset \BP^{n-1}$ (cf.
Example 1.7 of \cite{Hw10}.). \end{remark}

 The cone structure in Definition \ref{d.S}  naturally appears as the VMRT
of a uniruled projective manifold in the following way. The proof
is immediate.

\begin{proposition}\label{p.vmrt} In the setting of Definition \ref{d.S}, the cone structure on $M$ induced by $S$
comes from the  varieties of minimal rational tangents on the
blow-up ${\rm Bl}_S(\BP V)$ of $\BP V$ along the subvariety $S$,
associated to the minimal rational component parametrizing proper
transforms of lines on $\BP V$ intersecting $S$.
\end{proposition}

From Theorem \ref{t.1} and Proposition \ref{p.vmrt}, we can derive
 interesting consequences on  ${\rm Bl}_S(\BP V)$.

\begin{definition}\label{d.Louville}
Let $X$ be a uniruled manifold and $\sK \subset \RatCurves^n(X)$ a
minimal rational component. Let $\sC \subset \BP T(X)$ be the cone
structure associated to the VMRT of $\sK$. $X$ is said to have the {\em
Liouville property with respect to $\sK$} if every infinitesimal
automorphism of $\sC$ at a general point $x \in X$ extends to a
global holomorphic vector field on $X$, i.e.,  $\aut(\sC,
x) \cong \aut(X).$
\end{definition}

\begin{proposition}\label{p.2} Let $S \subset \BP V$ be a linearly normal nonsingular non-degenerate projective
variety such that ${\rm Sec}(S) \neq \BP V$. Then the blow-up
${\rm Bl}_S(\BP V)$ has the Liouville property with respect to the
minimal rational component parametrizing the proper transforms of
lines on $\BP V$ intersecting $S$. \end{proposition}

\begin{proof}
 By Theorem \ref{t.1}, we have
$\aut(\sC,x) \cong \aut(S)$. By Proposition \ref{p.vmrt}, the cone
structure $\sC$ induced by $S$ is the VMRT of a minimal component
on  ${\rm Bl}_S(\BP V)$. As any automorphism of ${\rm Bl}_S(\BP
V)$ preserves the VMRT, we have $\aut({\rm Bl}_S(\BP V)) \subset
\aut(\sC, x)$. On the other hand, as $S$ is $\Aut(S)$-invariant,
we have $\aut(S) \subset \aut({\rm Bl}_S(\BP V))$, which gives
$\aut({\rm Bl}_S(\BP V)) = \aut(\sC, x)$ for $x \in \BP V$
general.
\end{proof}

\begin{definition}\label{d.target}
A projective manifold $X$ is said to  have the {\em target
rigidity property}  if for any surjective morphism $f: Y \to X$,
and its deformation $\{f_t:Y \to X, |t| <1\}$, there exist
automorphisms $\sigma_t: X \to X$ such that $f_t = \sigma_t \circ
f$.
\end{definition}

All projective varieties which are not uniruled have the target
rigidity property (modulo \'etale factorizations)  by \cite{HKP}. All known examples of Fano
manifolds of Picard number 1, except projective space, have the
target rigidity property (cf. \cite{Hw09}). For nonsingular uniruled
projective varieties of higher Picard number, very few cases have
been studied (e.g.\cite{Hw07}). Some
examples can be obtained by the following easy lemma.

\begin{lemma}\label{l.equivariant}
Let $X$ be a nonsingular projective variety with the target
rigidity property. Let $\rho: X'\to X$ be a birational morphism
from a nonsingular projective variety which is equivariant in the
sense that there exists a group homomorphism $\rho^*: \Aut(X) \to
\Aut(X')$ with $$ \rho( \rho^*(g) \cdot y) = g \cdot \rho(y) \;
\mbox{ for any } g\in \Aut(X),  y \in X'.$$ Then $X'$ also has the
target rigidity property. \end{lemma}

\begin{proof} Given a deformation of surjective morphisms $f_t: Y\to X'$, the composition $\rho \circ f_t: Y \to X$ satisfies
$\rho \circ f_t = \sigma_t \circ \rho \circ f_0$ for some
$\sigma_t \in \Aut(X)$ from the target rigidity property of $X$.
Since $f_t = \sigma'_t \circ f_0$ with $\sigma'_t =
\rho^*(\sigma_t) \in \Aut(X'),$ $X'$ has the target rigidity
property.
\end{proof}

 For nonsingular uniruled projective
variety, the target rigidity property can be checked by the
Liouville property:

\begin{proposition}\label{p.trp}
Let $X$ be a nonsingular projective variety which has the
Liouville property with respect to a minimal component $\sK
\subset \RatCurves^n(X)$. Then $X$ has the target rigidity
property.
\end{proposition}
\begin{proof}
By the Stein factorization, it is easy to see that it suffices to
check
 the condition in Definition \ref{d.target} for generically finite surjective morphisms
 $\{ f_t: Y \to X, |t|<1 \}$ (cf. \cite{HKP} Section 2.2 for details). Let $\tau \in H^0(Y, f^*T(X))$ be the
Kodaira-Spencer class of the deformation $f_t$ at $0$. It suffices
to show that $\tau \in f^* H^0(X, T(X))$. As $f_t$ is generically
finite, we can regard $\tau$ as a multi-valued holomorphic vector
field on $X$. Take an analytic open subset near a general point $U
\subset Y$ such that $f_t|_U: U \to f_t(U)$ is biholomorphic for
$|t|<\epsilon$, then $\tau|_U$ can be regarded as a vector field
on $f(U)$. By Proposition 3 in \cite{HM99}, there are countably
many subvarieties $\sD_i \subset \BP T(Y), i=1,2,\ldots,$ (called
varieties of distinguished tangents in \cite{HM99}) such that for
any generically finite morphism $h: Y \to X$ and the dominant
rational map $dh: \BP T(Y) \dasharrow \BP T(X)$ defined by the
differential of $h$, the proper inverse image $d
h^{-1}(\sC|_{h(U)})$ coincide with some $\sD_i$. As the family
$df_t^{-1} (\sC), |t|<\epsilon$ is uncountable,  we have
$(df_t^{-1} (\sC))|_U = (df_0^{-1} (\sC))|_U$ for all $t$ small.
This implies that $\tau|_U$ is an infinitesimal automorphism of
the VMRT structure $\sC$, i.e. its germ at $x \in f_0(U)$ is an
element of $\aut(\sC,x)$. As the Liouville property holds, this
local vector field comes from a global vector field on $X$, which
gives $\tau \in f^* H^0(X, T(X))$.
\end{proof}

The following  is immediate from Proposition \ref{p.2}, Lemma
\ref{l.equivariant} and Proposition \ref{p.trp}.

\begin{corollary}\label{c.comp}
Let $S \subset \BP V$ be an irreducible linearly normal
nonsingular non-degenerate projective variety. Then the blow-up
${\rm Bl}_S(\BP V)$ of $\BP V$  along $S$ has the target rigidity
property. Moreover, if $X \to {\rm Bl}_S(\BP V)$ is the
composition of successive blow-ups along proper transforms of
$\Aut(S)$-invariant subvarieties in $\BP V$ not contained in $S$,
then $X$ satisfies the target rigidity property.
\end{corollary}

\begin{example}
Let $S \subset \BP V$ be the VMRT of an irreducible Hermitian
symmetric space $M$ of compact type such that $\Sec(S) \neq \BP
V$, i.e., those discussed in Propositions \ref{p.I}-\ref{p.VI}.
Let $Z \to {\rm Bl}_S(\BP V)$ be the composition of successive
blowing-ups along proper transforms of higher secant varieties of
$S$, from the smallest to the biggest.  Then
 $Z$ satisfies the target rigidity property by Corollary \ref{c.comp}.
When $M$ is of type I, II or III, this variety $Z$ is studied in
\cite{Th}, where it is  called complete collineations, complete
skew forms and complete quadrics, respectively. \end{example}

\bigskip
Baohua Fu

Institute of Mathematics, AMSS, Chinese Academy of Sciences, 55
ZhongGuanCun

East Road, Beijing, 100190, China and
 Korea
Institute for Advanced Study,

 Hoegiro 87, Seoul, 130-722, Korea

 bhfu@math.ac.cn

\bigskip
Jun-Muk Hwang

 Korea Institute for Advanced Study, Hoegiro 87,

Seoul, 130-722, Korea

jmhwang@kias.re.kr

\begin{thebibliography}{KSWZ}
\bibitem[AI]{AI} Arnold, V. I. and Ilyashenko, Yu. S.:
 Ordinary differential equations. in
Encyclopaedia Math. Sci., 1, Dynamical systems, I, 1-148,
Springer, Berlin, 1988
\bibitem[Ca]{Ca} Cartan, E.:  Les groupes de transformations continus,
infinis, simples. Ann. Sci. \'Ecole Norm. Sup. (3) {\bf 26} (1909)
93-161
\bibitem[Fu]{Fu}
Fu, B.: Inductive characterizations of hyperquadrics. Math. Ann.
{\bf 340} (2008),  185-194
\bibitem[FH]{FH} Fulton, W. and Harris, J.: {\em Representation theory. A first course}.
Graduate Texts in Mathematics. 129, Springer-Verlag, New York, 1991
\bibitem[Ge]{Ge} Gellhaus, C.: \"Aquivariante Kompaktifizierungen
des $\C^n$. Math. Zeit. {\bf 206} (1991) 211-217
\bibitem[Gu]{Gu} Guillemin, V.: The integrability problem for
G-structures. Trans. A.M.S. {\bf 116} (1965) 544-560
\bibitem[Hw01]{Hw01}  Hwang, J.-M.: Geometry of minimal rational curves
on Fano manifolds.   ICTP Lect. Notes {\bf 6} (2001) 335-393
\bibitem[Hw07]{Hw07}
Hwang, J.-M.: Deformation of holomorphic maps onto the blow-up of
the projective plane, Ann. Sci. \'Ecole Norm. Sup. (4) {\bf 40}
(2007),  179-189
\bibitem[Hw09]{Hw09} Hwang, J.-M.: Unobstructedness of deformations of holomorphic maps
 onto Fano manifolds of Picard number 1.  J. reine angew. Math. {\bf 637} (2009) 193-205
 \bibitem[Hw10]{Hw10} Hwang, J.-M.:
Equivalence problem for minimal rational curves with isotrivial
varieties of minimal rational tangents.   Ann. Sci. \'Ecole Norm.
Sup. (4) {\bf 43} (2010) 607-620

\bibitem[HK]{HK}
 Hwang, J.-M. and  Kebekus, S.: Geometry of chains of minimal rational
curves, J. Reine Angew. Math. {\bf 584} (2005) 173-194

\bibitem[HKP]{HKP}
 Hwang, J.-M.,  Kebekus, S. and Peternell, T.:
Holomorphic maps onto varieties of non-negative Kodaira dimension.
J. Alg. Geom. {\bf 15} (2006) 551-561

\bibitem[HM97]{HM97} Hwang, J.-M. and Mok, N.: Uniruled
projective manifolds with irreducible reductive G-structures. J.
reine angew. Math. {\bf 490} (1997) 55-64
%\bibitem[HM98]{HM98} Hwang, J.-M. and Mok, N.: Rigidity of irreducible Hermitian
%symmetric spaces of the compact type under Kahler deformation.
%Invent. Math. {\bf 131} (1998) 393-418
 \bibitem[HM99]{HM99} Hwang, J.-M. and Mok, N.: Holomorphic maps from rational homogeneous
spaces of Picard number 1 onto projective manifolds. Invent. math.
{\bf 136} (1999) 209-231
\bibitem[HM01]{HM01} Hwang, J.-M. and Mok, N.:
Cartan-Fubini type extension of holomorphic maps for Fano
manifolds of Picard number 1. Journal Math. Pures Appl. {\bf 80}
(2001) 563-575
%\bibitem[HM03]{H-M1}
%J.-M. Hwang, N. Mok, \emph{Finite morphisms onto Fano manifolds of
%Picard number 1 which have rational curves with trivial normal
%bundles} J. Algebraic Geom. 12 (2003), no. 4, 627--651
%\bibitem[HM04]{H-M2}
%J.-M. Hwang, N. Mok, \emph{Birationality of the tangent map for
%minimal rational curves} Asian J. Math. 8 (2004), no. 1, 51--63
\bibitem[HM05]{HM05} Hwang, J.-M. and Mok, N.:
Prolongations of infinitesimal linear automorphisms of projective
varieties and rigidity of rational homogeneous spaces of Picard
number 1 under K\"ahler deformation. Invent. math. {\bf 160}
(2005) 591-645
%\bibitem[HP07]{H-P} J.-M. Hwang,  T.
%Peternell, \emph{Holomorphic maps onto K?hler manifolds with
%non-negative Kodaira dimension}, J. Korean Math. Soc. 44 (2007),
%no. 5, 1079--1092.
%\bibitem[Ko]{Ko} Kobayashi, S.:
%{\em Transformation groups in differential geometry},
% Ergebnisse der Mathematik und ihrer Grenzgebiete, Band 70. Springer-Verlag, New York-Heidelberg, 1972.
\bibitem[IR]{IR}
Ionescu, P. and Russo, F.: Conic-connected manifolds, J. reine
angew. Math. {\bf 644} (2010), 145-157
\bibitem[KN]{KN}
Kobayashi, S. and  Nagano, T.: On filtered Lie algebras and
geometric structures I., J. Math. Mech. {\bf 13} (1964) 875--907
%\bibitem[LM]{LM} Landsberg, J. M. and Manivel, L.:
%On the projective geometry of rational homogeneous varieties.
%Comment. Math. Helv. {\bf 78} (2003) 65-100
%\bibitem[MS]{MS} Mori, S. and Sumihiro, H.: On Hartshorne's conjecture. J. Math. Kyoto Univ. {\bf 18} (1978) 523-533
\bibitem[Ot]{Ot} Ottaviani, G.: Spinor bundles on quadrics. Trans.
A. M. S. {\bf 307} (1988) 301-316
\bibitem[Pa]{Pa} Pasquier, B.: On some smooth projective two-orbit
varieties with Picard number 1. Math. Ann. {\bf 344} (2009) 963-987
\bibitem[PS]{PS}
Peternell, T.; Schneider, M.: Compactifications of $\C^n$: a
survey. Several complex variables and complex geometry, Part 2
(Santa Cruz, CA, 1989), 455-466, Proc. Sympos. Pure Math., 52,
Part 2, Amer. Math. Soc., Providence, RI, 1991
\bibitem[Th]{Th}
Thaddeus, M.: Complete collineations revisited. Math. Ann. {\bf
315} (1999), no. 3, 469-495
%\bibitem[Va]{Va}
%Vainsencher, I.: Complete collineations and blowing up determinantal
%ideals. Math. Ann. 267 (1984), no. 3, 417--432
\bibitem[Ya]{Ya} Yamaguchi, K.: Differential systems associated with simple graded Lie algebras.
Adv. Study Pure Math. {\bf 22}  (1993) 413-494
\bibitem[Za]{Za}
Zak, F. L.: {\em Tangents and secants of algebraic varieties}.  Translations of Mathematical Monographs, 127. American Mathematical Society, Providence, RI, 1993
\end{thebibliography}
\end{document}